\documentclass{article}

\usepackage[T1]{fontenc}
\usepackage[utf8]{inputenc}
\usepackage[british]{babel}

\usepackage{lmodern}

\usepackage{amsmath}
\usepackage{amsthm}
\usepackage{amsfonts}
\usepackage{hyperref}
\hypersetup{
	colorlinks = true,
	linkcolor = purple, 
	filecolor = purple, 
	citecolor = purple,
	urlcolor = purple
}
\usepackage{url} 
\usepackage{amssymb}
\usepackage{tikz-cd}
\usetikzlibrary{babel}
\allowdisplaybreaks

\usepackage[backend=biber, style=alphabetic, sorting=nyt, maxcitenames=7,
mincitenames=1, maxnames=7, abbreviate=true, block=space]{biblatex}

\numberwithin{equation}{section}

\newtheorem{lemma}{Lemma}[section]
\newtheorem{theorem}[lemma]{Theorem}
\newtheorem{proposition}[lemma]{Proposition}
\newtheorem{corollary}[lemma]{Corollary}

\newtheorem*{theorem*}{Theorem}
\newtheorem*{proposition*}{Proposition}

\theoremstyle{definition}
\newtheorem{rem}[lemma]{Remark}
\newtheorem{example}[lemma]{Example}
\newtheorem{notation}[lemma]{Notation}
\newtheorem{definition}[lemma]{Definition}

\newcommand{\OO}{\mathcal{O}}
\newcommand{\A}{\mathbb{A}}
\newcommand{\PP}{\mathbb{P}}

\newcommand{\Z}{\mathbb{Z}}

\newcommand{\R}{\mathbb{R}}
\newcommand{\C}{\mathbb{C}}
\newcommand{\m}{\mathfrak{m}}

\newcommand{\G}{\mathbb{G}}

\DeclareMathOperator{\Spec}{Spec}
\DeclareMathOperator{\Proj}{Proj}

\DeclareMathOperator{\id}{id}

\DeclareMathOperator{\rank}{rank}

\DeclareMathOperator{\Sym}{Sym}

\DeclareMathOperator{\ev}{ev}
\DeclareMathOperator{\End}{End}
\DeclareMathOperator{\GW}{GW}

\DeclareMathOperator{\can}{can}

\DeclareMathOperator{\Hom}{Hom}

\DeclareMathOperator{\Tr}{Tr}

\DeclareMathOperator{\sgn}{sgn}
\DeclareMathOperator{\Bl}{Bl}
\DeclareMathOperator{\Rep}{Rep}
\DeclareMathOperator{\gr}{gr}
\DeclareMathOperator{\characteristic}{char}
\DeclareMathOperator{\Stab}{Stab}

\DeclareMathOperator{\KQ}{KQ}
\DeclareMathOperator{\CH}{CH}

\newcommand{\SH}{\text{\textnormal{SH}}}
\newcommand{\Sm}{\text{\textnormal{Sm}}}
\newcommand{\Var}{\text{\textnormal{Var}}}

\newcommand{\rs}{\text{\textnormal{uh}}}
\newcommand{\PPal}{\text{\textnormal{PP}}}
\newcommand{\op}{\text{\textnormal{op}}}

\newcommand{\even}{\text{\textnormal{even}}}
\newcommand{\odd}{\text{\textnormal{odd}}}
\newcommand{\dR}{\text{\textnormal{dR}}}

\newcommand{\Et}{\text{\textnormal{Ét}}}

\makeatletter
\newcommand{\colim@}[2]{%
  \vtop{\m@th\ialign{##\cr
      \hfil$#1\operator@font colim$\hfil\cr
      \noalign{\nointerlineskip\kern1.5\ex@}#2\cr
      \noalign{\nointerlineskip\kern-\ex@}\cr}}%
}
\newcommand{\fcolim}{%
  \mathop{\mathpalette\colim@{\rightarrowfill@\textstyle}}\nmlimits@
}
\newcommand{\colim}{%
  \mathop{\mathpalette\colim@{}}\nmlimits@
}
\makeatother

\title{$\mathbb{A}^1$-Euler Characteristic of Low Symmetric Powers and Split
  Toric Varieties}
\author{Louisa F. Bröring}
\date{}

\begin{document}
\maketitle

\begin{abstract}
  For a smooth, projective scheme $X$ over a field $k$ or any variety $X$ if $k$
  has characteristic zero, we compute
  the compactly supported $\mathbb{A}^1$-Euler characteristic of $\operatorname{Sym}^2(X)$ if
  $\operatorname{char}(k) \ne 2$ and of $\operatorname{Sym}^3(X)$ if
  $\operatorname{char}(k) \ne 2,3$.
  We do so by extending the definition of a $G$-equivariant quadratic Euler
  characteristic first studied by Pajwani-Pál to arbitrary characteristic and
  by studying its relation to the $\mathbb{A}^1$-Euler characteristic of
  quotients.
  As an application, we show that the compactly supported $\mathbb{A}^1$-Euler
  characteristic of $\operatorname{Sym}^n(X)$ agrees with the prediction from
  the power structure constructed by Pajwani-Pál for $n = 2,3$.

  Furthermore, we compute the compactly supported $\mathbb{A}^1$-Euler
  characteristic of split toric varieties and show that the compactly
  supported $\mathbb{A}^1$-Euler characteristic of all of their symmetric
  powers agrees with the prediction from the power structure constructed by
  Pajwani-Pál.
\end{abstract}

\tableofcontents

\section*{Introduction}
\addcontentsline{toc}{section}{Introduction}
\markboth{Introduction}{}
The (compactly supported) $\A^1$-Euler characteristic is a refinement of the Euler characteristic
in topology to one in algebraic geometry, which was first
considered by Hoyois \cite{HoyoisQRGLVTF} in the case of smooth, projective
schemes. Its
construction uses motivic homotopy theory and for arbitrary varieties is due to Arcila-Maya, Bethea, Opie,
Wickelgren, and Zakharevich \cite{Wickelgren} in characteristic zero and to
Levine, Pepin Lehalleur, and Srinivas \cite{LevineLehalleurSrinivas} in
arbitrary characteristic. Let $k$ be a perfect field of
characteristic not two and $X$ a variety over $k$. Then the
$\A^1$-Euler characteristic $\chi_c(X/k)$ of $X$ is not an integer but an
element of the Grothendieck-Witt ring $\GW(k)$ of (virtual, non-degenerate) quadratic
forms over $k$. More specifically, the $\A^1$-Euler characteristic defines a motivic measure $\chi_c\colon K_0(\Var_k) \to \GW(k)$.

The $\A^1$-Euler characteristic carries a lot of information: if $k\subset \R$, and $X$ is a smooth, projective scheme over $k$, then the
rank of $\chi_c(X/k)$ equals the topological Euler characteristic of
$X(\C)$ and the signature of $\chi_c(X/k)$ with respect to the specified
embedding into $\R$ equals the topological Euler characteristic of $X(\R)$, see
\cite[Remark~2.3~(1) and Proposition~2.4~(6)]{LevineA}.
Furthermore, a theorem of Saito \cite[Theorem 2]{Saito} allows one to interpret the discriminant of
$\chi_c(X/k)$ in terms of the determinant of $\ell$-adic cohomology for any
$\ell$ that is coprime to the characteristic of the base-field, see \cite[Theorem
  2.22]{Pajwani-PalYZ}.

$\A^1$-Euler characteristics are often used in the programme of refined
enumerative geometry, which aims at proving counts in enumerative geometry
over more general fields, often in terms of quadratic forms and using motivic homotopy theory. One can often deduce known invariants or results by
taking ranks or signatures from the refined results.

In general, $\A^1$-Euler characteristics are hard to compute and there are
almost no examples of explicit computations of $\chi_c(X/k)$ for $X$ projective
but not smooth. In
particular, not much is known about the $\A^1$-Euler characteristic of
quotients.

\subsection*{The $\A^1$-Euler Characteristic of Low Symmetric Powers and the Equivariant Quadratic Euler Characteristic}
\addcontentsline{toc}{subsection}{The $\A^1$-Euler Characteristic of Low Symmetric Powers and the Equivariant Quadratic Euler Characteristic}
Let $X$ be a smooth, quasi-projective scheme over $k$ and let $S_n$ act on
$X^n$ by permuting the factors. The quotient $\Sym^nX := X^n/S_n$ is called
the $n$-th symmetric power of $X$. In the following, let $\langle a \rangle \in \GW(k)$ denote
the quadratic form $x \mapsto ax^2$ and let $H$ denote the hyperbolic form
$\langle 1\rangle + \langle -1\rangle \in \GW(k)$. The aim of this paper is to prove the
following theorems and some consequences of them.

\begin{theorem*}[see Theorem \ref{thm:qec-sym2}]
  Let $X$ be a connected, smooth, projective scheme over $k$. Express the $\A^1$-Euler
  characteristic of $X$ as $\chi_c(X/k) = \beta -mH$ with
  $\beta = \sum_{i=1}^n\langle \alpha_i\rangle$ and $m\ge 0$.
  Then if $\dim X$ is odd, we have in $\GW(k)$
  \[
    \chi_c(\Sym^2X/k) = n\langle 2\rangle + \sum_{1\le i<j\le n}\langle
    \alpha_i\alpha_j\rangle + m\langle -1\rangle + (m^2 - (n+1)m)\cdot H.
  \]
  If $\dim X$ is even, we have
  \[
    \begin{split}
      \chi_c(\Sym^2X/k) &= n\langle 2\rangle + \sum_{1\le i<j\le n}\langle
                          \alpha_i\alpha_j\rangle +
                          m\langle -1\rangle + (m^2 - (n+1)m)\cdot H\\
      &\quad + (\langle -2\rangle -
      \langle -1 \rangle) \cdot \chi_c(X/k)
    \end{split}
  \]
  in $\GW(k)$.
\end{theorem*}

\begin{theorem*}[see Theorem \ref{thm:sym3} and Remark \ref{rem:sym3-formula-indep}]
  Assume that $\characteristic k \ne 2,3$ and
  let $X$ be a connected, smooth, projective scheme over $k$. Express the $\A^1$-Euler characteristic of $X$ as
  $\chi_c(X/k) = \sum_{i=1}^m\langle \alpha_i\rangle -lH \in \GW(k)$ with
  $l\ge 0$ and write $\beta = \sum_{i=1}^m \langle
  \alpha_i\rangle$. Then
  \begin{align*}
    \chi_c(\Sym^3 X/k)
    &= \langle 6\rangle
      \cdot \Bigg((\langle 1\rangle + (m-1)\langle 3\rangle)\beta +
      \sum_{1\le i<j<k\le m}\langle 6\alpha_i\alpha_j\alpha_k\rangle\\
    &\qquad  + l\langle -6\rangle \beta +
      ml(l-1)H - \left(m+\binom m2\right)lH - \frac 12 \binom{2l}{3}H\Bigg)\\
    &\quad + (\langle 1\rangle - \langle 6\rangle)\cdot
      \chi_c(X/k) + (\langle 1\rangle -\langle 2\rangle)(\chi_c(X/k) - \langle
      1\rangle)\chi_c(X/k)\\
  \end{align*}
  in $\GW(k)$.
\end{theorem*}

The proof of both of these theorems, which can be found in the Sections~\ref{sec:explicit-sym2} and \ref{sec:explicit-sym3}, is based on the following observation: let $G$ be
a finite group of order prime to the characteristic of $k$ acting on a smooth,
projective $k$-scheme $X$. If $X/G$ is smooth, there is an explicit
description of $\chi_c((X/G)/k)$ in terms of a composition of the cup product
and trace on the $G$-invariants of the Hodge cohomology of $X$.

The scheme $\Sym^nX$ is in general not smooth for $n \ge 2$. In order to be
able to apply the above observation, we blow $X^n$ up a number of times (along
smooth centers), until the
quotient $Y/S_n$ of the final blow-up is smooth. Since the Hodge cohomology
changes when blowing up, we have to keep track of the changes. In order to do
so, we introduce the notion of a $G$-equivariant quadratic Euler characteristic
$\chi_G(X/k)$ of a smooth, projective $k$-scheme $X$ with a $G$-action in Section \ref{sec:equivariant-Euler-char}. This Euler characteristic is similar to the notion of $G$-equivariant de Rham Euler
characteristic introduced by Pajwani-Pál
\cite[Definition~3.5]{Pajwani-PalYZ}. This $G$-equivariant quadratic Euler
characteristic $\chi_G(X/k)$ is valued in the Grothendieck-Witt ring $\GW^G(k)$ of
so-called $G$-equivariant non-degenerate symmetric bilinear forms over
$k$. The $G$-equivariant quadratic Euler characteristic is well-behaved with
respect to blowing up and the operation of restricting a $G$-equivariant
non-degenerate symmetric bilinear form to its invariants induces a
group homomorphism $\GW^G(k) \to \GW(k)$. The observation that for smooth
$X/G$ we can compute $\chi_c(X/G)$ via the Hodge cohomology and trace of $X$ translates to the following theorem.

\begin{theorem*}[see Theorem \ref{thm:smooth-quotient-euler-char}]
  Let $X$ be a smooth, projective scheme over $k$ and let $G$ be a group of order prime to $\characteristic k$ acting on $X$. Let $Y = X/G$ be the quotient and $\pi\colon X \to Y$
  the quotient map. Assume that $G$ acts freely on an open, dense subset of $X$ and assume that $Y$ is smooth.
  Then $\chi_c(Y/k) = \langle|G|\rangle \cdot \chi_G(X/k)^G$.
\end{theorem*}

We use $\chi_G(X/k)$ facilitate the computation of $\chi_c(\Sym^nX/k)$ for
$n=2,3$.

The computation of $\chi_c(\Sym^nX/k)$ for smooth,
projective $X$ and $n=2,3$ uses the fact that we
can find a sequence of blow-ups in smooth subschemes of $X^n$ stable under the
$S_n$-action such that the $S_n$-quotient of the final blow-up is smooth. We achieve this using a
smoothness criterion for actions by cyclic groups of prime order, see Theorem~\ref{thm:divisor-smooth-quotient}, and by
filtering $S_n$ by cyclic groups. For a smooth scheme with a $\Z/n\Z$-action,
in general one cannot find a sequence of equivariant blow-up at smooth centres such that
the quotient is smooth. This complication prevents our strategy
from working for $n \ge 4$.

\subsection*{Application: Comparison With the Power-Structure Prediction of Pajwani-Pál}
\addcontentsline{toc}{subsection}{Application: Comparison With the Power-Structure Prediction of Pajwani-Pál}
As an application of our computation of $\chi_c(\Sym^nX/k)$ for $n = 2,3$, we
compare our results to a prediction of
Pajwani-Pál \cite[Corollary~3.27]{Pajwani-PalPS}. A finitely determined  power structure on a ring
$R$ is a collection of functions $a_n\colon R\to R$ for $n\in \Z_{\ge 0}$
subject to some relations mimicking the definition of $f(t)^\alpha$ for $f(t) \in (\Z[[t]])^\times$ and $\alpha \in \Z$, see Definition~\ref{def:finite-ps-char}.
In characteristic zero, the functions $a_n^{\Sym}\colon
K_0(\Var_k)\to K_0(\Var_k), [X] \mapsto [\Sym^nX]$ form a power structure on
$K_0(\Var_k)$. If $k$ has positive characteristic, the same holds after
identifying elements in $K_0(\Var_k)$ linked via a universal homeomorphism, see Definition \ref{def:universal-homeo}.
Pajwani-Pál \cite[Corollary~3.27]{Pajwani-PalPS} construct a
power structure $a_\ast^\PPal$ on $\GW(k)$ with the property that if there is
a power structure $a_\ast$ on $\GW(k)$ such that $\chi_c$ respects the power structures, i.e. $a_n(\chi_c(X/k)) =
\chi_c(\Sym^nX/k)$ for all varieties $X$ over $k$ and all $n \ge 0$, then
$a_\ast$ has to agree with $a_\ast^\PPal$. See Section~\ref{sec:power-structures} for an introduction to power structures.

It is still an open question whether $\chi_c$ respects the power structures in
general. There are, however, certain classes of schemes for which $\chi_c$ respects
the power structures.

The explicit computations of $\chi_c(\Sym^2X/k)$ and $\chi_c(\Sym^3X/k)$ allow
us to check the compatibility with $a_\ast^\PPal$ for $\ast = 2,3$, which is detailed in Section \ref{sec:comparison-sym23}.

\begin{theorem*}[see Theorem \ref{thm:compatibility-a2-smooth-projective}]
  For a smooth, projective scheme $X$ over $k$, we have
  \[
    \chi_c(\Sym^2X/k) = a_2^\PPal(\chi_c(X/k))
  \]
  in $\GW(k)$.
\end{theorem*}

\begin{theorem*}[see Theorem \ref{thm:sym3-comparison}]
  Assume that $\characteristic k \ne 2,3$ and
  let $X$ be a smooth, projective scheme over $k$. Then
  \[
    \chi_c(\Sym^3X/k) = a_3^\PPal(\chi_c(X/k))
  \]
  in $\GW(k)$.
\end{theorem*}

\begin{theorem*}[see Theorem \ref{thm:sym23-compatibility-char-zero}]
  Let $X$ be a variety over $k$ and suppose $k$ has characteristic zero. Then
  \[
    \chi_c(\Sym^2X/k) = a_2^\PPal(\chi_c(X/k)) \quad \text{and}\quad \chi_c(\Sym^3X/k) = a_3^\PPal(\chi_c(X/k))
  \]
  in $\GW(k)$.
\end{theorem*}

For the proof of Theorem \ref{thm:sym23-compatibility-char-zero} one reduces the statement to the Theorems \ref{thm:compatibility-a2-smooth-projective} and \ref{thm:sym3-comparison} using Bittner's presentation \cite[Theorem 3.1]{bittner_universal_2004}.

\subsection*{The $\A^1$-Euler Characteristic of Toric Varieties and Their Symmetric Powers}
\addcontentsline{toc}{subsection}{The $\A^1$-Euler Characteristic of Toric Varieties and Their Symmetric Powers}
In Section~\ref{sec:toric}, we compute the $\A^1$-Euler characteristic of a split toric
variety and verify that all toric varieties satisfy the prediction of
Pajwani-Pál; this relies on the work of Pajwani-Rohrbach-Viergever \cite{PajwaniRohrbachViergever}.

\begin{proposition*}[see Proposition \ref{prop:euler-char-toric}]
  Let $\Sigma$ be a fan in a rank $n$ lattice $N$. Let $X_\Sigma$ be the
  associated split toric variety over $k$ and let $\Sigma_i$ be the collection
  of cones in $\Sigma$ of rank $n-i$. Then we have
  \[
    \chi_c(X_\Sigma/k) = |\Sigma_0|\cdot \langle 1\rangle + \sum_{i\ge 1} |\Sigma_i| \cdot
    (-2)^{i-1}(\langle-1\rangle - \langle 1\rangle) \in \GW(k).
  \]
\end{proposition*}

\begin{proposition*}[see Proposition \ref{prop:split-toric-symmetrisable}]
  Let $X = X_\Sigma$ be a split toric variety over $k$ and $n \in \Z_{\ge
    0}$. Then
  \[
    \chi_c(\Sym^nX/k) = a_n^\PPal(\chi_c(X/k))
  \]
  in $\GW(k)$.
\end{proposition*}

The results presented in this article also form two chapters in the author's PhD-thesis, which was submitted on 22 January 2025.

\subsection*{Related Work}
\addcontentsline{toc}{subsection}{Related Work}
There are some related results concerning the compatibility with the power structure. For summarising these, we say following \cite[Definition~4.1]{PajwaniRohrbachViergever} that a variety $X$ over $k$ is \emph{symmetrisable} if $\chi_c(\Sym^nX/k) = a_n^\PPal(\chi_c(X/k))$ holds for all $n \ge 0$.

Let $K^\rs(\Var_k)$ be the quotient of $K_0(\Var_k)$ identifying universally homeomrphic varieties, see Definition \ref{def:universal-homeo}.
Extending earlier work of Pajwani-Pál \cite[Corollary~4.30 and Lemma~2.9]{Pajwani-PalPS},  Pajwani-Rohrbach-Viergever
\cite[Theorem~4.10 and Corollary~4.3]{PajwaniRohrbachViergever} show that all étale linear
varieties are symmetrisable, that is varieties lying in the subring of
$K_0^\rs(\Var_k)$ generated by $[\A^1]$ and zero-dimensional varieties.
They also show that all symmetrisable varieties form a submodule
over the subgroup of $K_0^\rs(\Var_k)$ generated by zero-dimensional
varieties. Cellular varieties in the sense of \cite[2189]{LevineA} are étale linear.

Together with Anna Viergever \cite[Theorem~27]{bv2024quadratic}, we show that all curves are symmetrisable.

Bejleri-McKean \cite[Theorem~7.10]{mckean2024symmetricpowersnullmotivic}
prove that showing that all varieties are symmetrisable is equivalent to showing that all varieties $X$ with $\chi_c(X/k) = 0$ are symmetrisable. They
also prove \cite[Theorem~9.4]{mckean2024symmetricpowersnullmotivic} that all varieties over $k$ are symmetrisable whenever $\GW(k)$
is torsion free. This
is for example the case for Pythagorean fields such as $\R$ by \cite[Chapter 2, Theorem
4.10 (iii)]{ScharlauQHF}.

In the computations in this work, we follow more the approach of Pajwani-Pál
\cite{Pajwani-PalPS} and Pajwani-Rohrbach-Viergever
\cite{PajwaniRohrbachViergever} to extend the class of varieties that are
symmetrisable over arbitrary fields.

\subsection*{Notation and Conventions}
\addcontentsline{toc}{subsection}{Notation and Conventions}
Throughout, let $k$ be a perfect field that is not of characteristic two. The
assumption that $k$ is perfect is not strictly necessary and we do so to avoid
technical complications. All results on Euler characteristics can be obtained for non-perfect fields from the perfect case via a base-change to a perfect closure using Lemma~\ref{lem:QuotFlatBaseChange}~(\ref{lem:QuotFlatBaseChange:quot}) and \cite[Remarks 2.1 2.]{LevineA}.

For schemes
$X$ and $Y$ over $k$, we denote their fibre product over $k$ by $X\times Y$. For
$k$-vector spaces $V$ and $W$, we denote their tensor product as $k$-vector
spaces by $V\otimes W$. For a
point $x \in X$ on a scheme $X$, we write
$k(x)$ for the residue field at $x$.

For a smooth $k$-scheme $X$, we write $\CH^\ast(X)$ for the Chow-ring of $X$
and $TX$ for the tangent bundle of $X$. For a smooth, closed subscheme
$i\colon Y\subset X$, we write $N_YX = i^\ast TX/TY$ for the normal bundle of
$Y$ in $X$.

We write $\A^1$ for the affine line $\A^1_k$ over $k$ and $\PP^n$ for the
$n$-dimensional projective space $\PP^n_k$ over $k$.
If $X$ is a smooth, projective scheme over $k$, we write $\Omega_{X/k}$ for
the sheaf of Kähler differentials of $X$ over $k$. For $i \in \Z_{\ge 0}$, we write
$\Omega^i_{X/k} := \bigwedge^i\Omega_{X/k}$. For a smooth, projective
$k$-scheme $X$ that is equidimensional of dimension $n$, we write
$\Tr_{X/k}\colon H^n(X,\Omega^n_{X/k})\to k$ for the trace map in induced by
Serre duality. If $X$ is clear from the context, we also write $\Tr$ for
$\Tr_{X/k}$.

For $p,q, p', q' \in \Z_{\ge 0}$, we write $ab := a\cdot b :=
a \cup b \in H^{q+q'}(X,\Omega^{p+p'}_X)$ for the cup product of $a \in
H^q(X,\Omega^p_X)$  and $b \in H^{q'}(X,\Omega^{p'}_X)$.
This cup-product is graded commutative if we
consider $a \in H^q(X,\Omega^p_X)$ to be of degree $q-p$; so,  we have $ab =
(-1)^{(q-p)(q'-p')} ba$ for
$a \in H^q(X,\Omega^p_X)$ and $b \in H^{q'}(X,\Omega^{p'}_X)$.

For a quasi-compact and quasi-separated scheme $X$, let $\SH(X)$ denote the stable
motivic homotopy category of $X$ as described by Hoyois \cite{HoyoisSOEMHT},
and let $1_X\in \SH(X)$ be the
unit object with respect to the smash product in $\SH(X)$.

We denote the group of permutations on the set $\{1, \dots, n\}$ by $S_n$ and
write $\sgn\colon S_n \to \{\pm 1\}$ for the sign homomorphism. We denote the
kernel of $\sgn\colon S_n \to \{\pm 1\}$ by $A_n$. For a
quasi-projective $k$-scheme $X$ and $n \in \Z_{\ge 0}$, we denote the $n$-th symmetric power of $X$
by $\Sym^nX = X^n/S_n$. Here $S_n$ acts by permuting the factors of $X^n$.

\subsection*{Acknowledgements}
\addcontentsline{toc}{subsection}{Acknowledgements}
I would like to extend my deepest gratitude to my adviser,
Marc Levine, for suggesting this project to me and for all
of his invaluable support over the past years.
I am also extremely grateful to Anna Viergever for many insightful
discussions that straightened my thoughts around this project.
I would like to extend my sincere thanks to Sabrina Pauli for many fascinating
discussions and for feedback on my thesis.
Many thanks to Jesse Pajwani for enlightening discussions
about equivariant quadratic Euler characteristics and computations of Euler
characteristics.
I would also like to thank to Dhyan Aranha, Tom Bachmann, Linda Carnevale, Jochen Heinloth, Fabien Morel, Clémentine
Lemarié{-}{-}Rieusset, Herman Rohrbach, Hind Souly, and my second adviser, Johannes Sprang, for interesting discussions and explaining basic concepts
to me. I also want to thank Niklas Müller for organising a seminar on
variation of GIT structure, which led to my computations in Section
\ref{sec:toric}, and Ulrich Görtz for pointing me to the paper \cite{Knighten1973DiifQuotients}.
I am also grateful to Stefan Schreieder and the algebraic geometry group at
the Leibniz Universität Hannover for offering me a place to work when I was
visiting my partner and for the welcoming atmosphere that I encountered there.

This work was funded by Deutsche Forschungsgemeinschaft
(DFG, German Research Foundation) - Research Training Group 2553 -
Projektnummer 412744520.

\section{The $\A^1$-Euler Characteristic}
\label{sec:euler-char}
The (compactly supported) $\A^1$-Euler characteristic is a refinement of the
(compactly supported) topological Euler characteristic to varieties over
$k$. This refinement is a motivic measure valued in the Grothendieck-Witt ring of quadratic
forms. The $\A^1$-Euler characteristic was first studied in the case of smooth
projective schemes by Hoyois \cite{HoyoisQRGLVTF}. The first construction of
the $\A^1$-Euler characteristic as a motivic measure on the Grothendieck ring
of varieties in characteristic zero is due to Arcila-Maya, Bethea, Opie,
Wickelgren, and Zakharevich \cite[Theorem~2.13]{Wickelgren}.  We shall give a
brief account of its construction in arbitrary characteristic due to Levine, Pepin Lehalleur, and Srinivas
\cite[Section~5.1]{LevineLehalleurSrinivas} here.

\begin{definition}
  Let $K_0(\Var_k)$ be the \emph{Grothendieck group of varieties} over
  $k$. That is, $K_0(\Var_k)$ is the free Abelian group generated by isomorphism classes $[X]$ of
  varieties $X$ over $k$ modulo the scissor-relation: for $Z\subset X$ a closed
  subvariety, we have
  \[
    [X] = [Z] + [X\setminus Z]
  \]
  in $K_0(\Var_k)$.

  The product rule $[X][Y] = [X\times_kY]$ induces a ring structure on
  $K_0(\Var_k)$. A ring homomorphism out of $K_0(\Var_k)$ is called
  a \emph{motivic measure}.
\end{definition}

\begin{definition}
  The \emph{Grothendieck-Witt ring $\GW(k)$} over $k$ is defined as
  the group completion of the monoid of isometry classes of
  non-degenerate quadratic forms over $k$ with respect to orthogonal
  direct sums.

  The tensor product of quadratic forms induces a ring structure on the group
  $\GW(k)$.
\end{definition}

\begin{rem}
  \label{rem:gw-rel}
  Over a field $k$ whose characteristic is not $2$, non-degenerate quadratic
  forms are
  the same as non-degenerate symmetric bilinear forms. In this setting, elementary linear
  algebra computations yield the following presentation of $\GW(k)$: the
  quadratic forms $\langle a \rangle: x \mapsto ax^2$ for $x \in k^\times$
  generate $\GW(k)$ and these forms are subject to the following relations for
  $a, b\in k^\times$:
  \begin{enumerate}
  \item $\langle a\rangle \cdot \langle b\rangle = \langle a b \rangle$,
    \label{rem:gw-rel:mult}
  \item
    $\langle a\rangle + \langle b\rangle = \langle a+b\rangle + \langle
    ab(a+b)\rangle$, whenever $a+ b \in k^\times$,
    \label{rem:gw-rel:add}
  \item $\langle ab^2\rangle = \langle a\rangle$, and
    \label{rem:gw-rel:square}
  \item
    $\langle a\rangle + \langle -a\rangle = \langle 1\rangle + \langle
    -1\rangle =: H$.
    \label{rem:gw-rel:hyp}
  \end{enumerate}
  The quadratic form $H = \langle 1\rangle + \langle -1\rangle$ is called the
  \emph{hyperbolic form}.

  The first description of this presentation is due to Witt \cite[Section~1]{WittQF}. The
  presentation presented here, is \cite[Lemma~3.9]{MorelATF}, where the presentation is deduced from
  the corresponding statement for Witt rings, proven by Milnor and Husemoller
  \cite[Lemma~(1.1) in Chapter~4]{MilnorSBF}. Note that the relation (\ref{rem:gw-rel:mult}) is not included in
  \cite[Lemma~3.9]{MorelATF} because the description there provides a
  description of $\GW(k)$ as a group. We include it here to describe the ring
  structure on $\GW(k)$.
\end{rem}

Morel-Voevodsky \cite{MorelVoevodsky99HTS} construct an unstable motivic homotopy
category. From this, Voevodsky \cite{Voevodsky98HT} constructs a stable motivic
homotopy category $\SH(S)$ for a Noetherian base scheme $S$. Another
construction of $\SH(S)$ for $S$ a quasi-compact and quasi-separated scheme is due to Hoyois
\cite{HoyoisSOEMHT}. The category $\SH(S)$ is a symmetric monoidal
category and we shall denote its unit object by $1_S \in SH(S)$. See
\cite{HoyoisSOEMHT} for an introduction to $\SH(S)$. Similar to the
stable homotopy category in topology, there is an infinite suspension functor
$\Sigma^\infty\colon \Sm_S \to \SH(S)$ from the category of smooth schemes
over $S$. There is also a construction of a six functor formalism for $\SH$,
as can be found in \cite{HoyoisSOEMHT}.

For $k$ a perfect field, Morel \cite[Thm.~6.4.1, Rmk.~6.4.2]{MorelIAHT} constructs an isomorphism $\GW(k)
\to \End_{\SH(k)}(1_k)$ from the Grothendieck-Witt ring of quadratic forms
over $k$ to the endomorphisms
of the unit of $\SH(k)$.

If $X$ is a variety over a field $k$ of characteristic zero with structure map
$p_X\colon X\to\Spec(k)$, one can show using resolution of singularities that
$(p_X)_!(1_X)\in \SH(k)$ is a strongly dualisable object, as introduced by Dold-Puppe \cite[Theorem~1.3]{dold80-categoricalEC}. The categorical
Euler characteristic of $(p_X)_!(1_X)$ now computes the
\emph{compactly supported
  $\mathbb{A}^1$-Euler characteristic} $\chi_c(X/k)\in \GW(k)$.

If $X$ is a variety over a field $k$ of
characteristic $p>2$ with structure map $p_X: X\to\Spec(k)$, one can
show using a result of Riou \cite[Corollary~B.2]{Riou} that $(p_X)_!(1_X)\in
\SH(k)[\frac 1p]$ is a strongly
dualisable object after inverting
$p$.
The categorical Euler characteristic of $(p_X)_!(1_X)$ now yields the
Euler characteristic $\chi_c(X/k)\in \GW(k)[\frac{1}{p}]$, where one uses
Morel's isomorphism to identify $\End_{\SH(k)[1/p]}(1_{k})$ with
$\GW(k)[\frac{1}{p}]$. One then shows that
$\chi_c(X/k)$ lies in the image of the injective map
$\GW(k)\to\GW(k)[\frac{1}{p}]$, so that one can define the compactly supported
$\mathbb{A}^1$-Euler characteristic as an element $\chi_c(X/k)\in \GW(k)$.

\begin{definition}
  The \emph{(compactly supported) $\mathbb{A}^1$-Euler characteristic}
  $\chi_c(X/k)\in \GW(k)$ is defined to be the element corresponding with
  $\chi^{\text{cat}}((p_X)_!(1_X)) \in \End_{\SH(k)[\frac{1}{e}]}(1_k)$ under
  Morel's isomorphism. Here $e$ denotes the exponential characteristic of $k$.

  This yields a motivic measure $\chi_c\colon K_0(\Var_k)\to \GW(k)$.
\end{definition}

Since we only consider the compactly supported $\A^1$-Euler characteristic in this article, we shall drop the compactly supported and just write $\A^1$-Euler characteristic for brevity. 

\begin{rem}
  Let $k \subset \R$ and let $X$ be a smooth, projective scheme over $k$.
  Levine \cite[Remark 2.3]{LevineA} shows that the rank of $\chi_c(X/k)$ agrees with
  the topological Euler characteristic of $X(\C)$ and the signature
  agrees with the topological Euler characteristic of $X(\R)$.

  If $k$ is not of characteristic $2$ and $X$ is a smooth, projective scheme
  over $k$, a theorem of
  Saito \cite[Theorem 2]{Saito} can be used to interpret the discriminant of
  $\chi_c(X/k)$ in terms of the determinant of $\ell$-adic cohomology for any
  $\ell$ that is coprime to the characteristic of $k$. A more detailed
  description can be found in \cite[Theorem
  2.22]{Pajwani-PalYZ}.
\end{rem}

\begin{example}
  We have
  \[
    \chi_c(\PP^n(k)) = \sum_{i=0}^n \langle (-1)^i\rangle =
    \begin{cases}
      \langle 1\rangle + \frac n2\cdot H & \text{if $n$ is even},\\
      \frac{n+1}{2} H & \text{if $n$ is odd}.
    \end{cases}
  \]
  See for example \cite[Example 1.7]{HoyoisQRGLVTF} or \cite[Example 2.6 1.]{LevineA} for a proof.
\end{example}

$\A^1$-Euler characteristics are in general hard to compute.
The Motivic Gauß-Bonnet Theorem of Levine-Raksit \cite[Theorem 5.3 and Theorem
8.4]{LevineGB} provides a
way to compute the $\A^1$-Euler characteristic using pushforwards of Euler classes. A more general version is due to
Déglise-Jin-Khan \cite[Theorem 4.6.1]{DegliseFCMHT}.
We content ourselves here
with two non-trivial consequences of the Motivic Gauß-Bonnet Theorem,
which are enough to prove the results in this paper.

\begin{theorem}[{\cite[Theorem 8.6 (3)]{LevineGB}}]
  \label{thm:motivic-gauss-bonnet}
  Let $X$ be an equidimensional, smooth, projective scheme over $k$. Then we have
  \begin{equation}\label{eq:motivic-gauss-bonnet:formula}
    \chi_c(X/k) = \left(\bigoplus_{i,j=0}^{\dim_kX}
      H^i(X,\Omega^j_{X/k})[j-i],q^{\text{Hdg}}\right) \in \GW(k)
  \end{equation}
  under the identification of $\GW(k)$ with the Grothendieck-Witt ring of
  symmetric objects of degree $0$ for the category of
  bounded chain complexes of finite-dimensional $k$-vector spaces, as
  detailed by Panin-Walter \cite[Section 4 and
  5]{PaninWalter2018BO}. Here, $q^{\text{Hdg}}$ is the symmetric bilinear form
  \[
    \bigoplus_{i,j=0}^{\dim_kX}
    H^i(X,\Omega^j_{X/k})[j-i] \otimes \bigoplus_{i,j=0}^{\dim_kX}
    H^i(X,\Omega^j_{X/k})[j-i] \xrightarrow{\cup}
    H^n(X,\Omega^n_{X/k})[0] \xrightarrow{\Tr} k
  \]
  where $\Tr$ is the trace map from Serre duality.

  Alternatively, \eqref{eq:motivic-gauss-bonnet:formula} holds after
  applying a canonical map from $\Z$-graded non-degenerate symmetric bilinear
  forms over $k$ to $\GW(k)$, as detailed in Example
  \ref{ex:graded-euler-char}.
\end{theorem}

\begin{corollary}[{\cite[Corollary 8.7]{LevineGB}}]
  \label{cor:motivic-gb}
  Let $X$ be an $n$-dimensional, smooth, projective scheme over $k$.
  \begin{itemize}
  \item If $n$ is odd, then $\chi_c(X/k) = m\cdot H$ for some integer $m \in \Z$.
  \item If $n$ is even, then $\chi_c (X/k) = Q + m\cdot H$ for some integer $m \in
    \Z$ and $Q$ is the quadratic form corresponding with the symmetric bilinear
    form given by the composition
    \[
      H^{n/2}(X,\Omega^{n/2}_X)\times H^{n/2}(X,\Omega^{n/2}_X) \xrightarrow{\cup}
      H^{n}(X,\Omega_X^{n}) \xrightarrow{\Tr} k.
    \]
    Here, $\cup$ is the cup product and $\Tr$ is the trace map from Serre
    duality.
  \end{itemize}
\end{corollary}


\begin{example}
  \label{ex:chi-affine-line-gm}
  The scissors relation yields $[\PP^n] = [\A^n] + [\PP^{n-1}]$ and $[\G_m] =
  [\A^1] - [\Spec k]$. Hence,
  \[
    \chi_c(\A^n) = \chi_c(\PP^n) - \chi_c(\PP^{n-1}) = \langle (-1)^n\rangle
    \quad\text{and}\quad
    \chi_c(\G_m) = \chi_c(\A^1) - \chi_c(\Spec k) = \langle -1\rangle -
    \langle 1\rangle.
  \]
\end{example}

\begin{rem}
  The compactly supported $\A^1$-Euler characteristic satisfies all relations
  satisfied in $K_0(\Var_k)$. These include:
  \begin{itemize}
  \item For $k$-varieties $X$ and $Y$, we have $\chi_c(X\times Y/k) =
    \chi_c(X/k)\cdot \chi_c(Y/k)$.
  \item If $E\to X$ is a projective bundle of rank $r$ over a $k$-variety $X$,
    then $\chi_c(E/k) = \chi_c(\PP^r/k) \cdot \chi_c(X/k)$.
  \item If $\tilde X = \Bl_ZX$ is the blow-up of a $k$-variety $X$ in a closed
    subvariety $Z$ and if $E \subset \tilde X$ is the exceptional divisor,
    then $\chi_c(\tilde X/k) = \chi_c(X/k) - \chi_c(Z/k) +
    \chi_c(E/k)$.
  \item If a $k$-variety $X = U \cup V$ is the disjoint union of two locally
    closed subvarieties $U,V \subset X$, then $\chi_c(X/k) = \chi_c(U/k) +
    \chi_c(V/k)$.\qedhere
  \end{itemize}
\end{rem}

The $\A^1$-Euler characteristic is not an injective motivic measure. There
are some meaningful ideals lying in the kernel of $\chi_c$. One of these is the
ideal generated by universally homeomorphic varieties, as proven by Bejleri-McKean
\cite[Section~5]{mckean2024symmetricpowersnullmotivic}.

\begin{definition}
  \label{def:universal-homeo}
  A morphism of schemes $f\colon X \to Y$ is a \emph{universal homeomorphism}
  if $X\times_YY' \to Y'$ is a homeomorphism of underlying topological spaces
  for every scheme $Y' \to Y$.

  Let
  \[
    I := ([X]-[Y]; \text{there is a universal homeomorphism $f\colon X \to
      Y$})\subset K_0(\Var_k)
  \]
  be the ideal generated by differences of $k$-varieties
  linked by a universal homeomorphism. We define the \emph{Grothendieck
    ring of $k$-varieties up to universal homeomorphism} as the quotient ring
  $K^{\rs}_0(\Var_k) := K_0(\Var_k)/I$.
\end{definition}

\begin{rem}[{\cite[Chapter~2, Proposition~4.4.6 and
    Corollary~4.4.7]{CLNS2018MotivicIntegration}}]
  We have $K_0(\Var_k) =
  K^\rs_0(\Var_k)$ if $k$ has characteristic zero.
\end{rem}

\begin{proposition}[{\cite[Corollary 5.4]{mckean2024symmetricpowersnullmotivic}}]
  The compactly supported $\A^1$-Euler characteristic factors through
  $K_0^\rs(\Var_k)$.
\end{proposition}

See \cite[Section~2]{LevineA} and \cite[Sections~2--5]{mckean2024symmetricpowersnullmotivic}
for a detailed exposition of the $\A^1$-Euler characteristic and its basic
properties.

\section{Power Structures and Symmetric Powers}
\label{sec:power-structures}
The notion of a power structure was first studied and defined by Gusein-Zade,
Luengo, and Melle-Hernández
\cite{Gusein-Zade-Luengo-Melle-HernandezPS} and later refined in \cite{Gusein-Zade-Luengo-Melle-HernandezHS}. Similar to how $\lambda$-rings
are a generalisation of the notion of exterior power, a power structure is a
generalisation of the notion of taking powers. Throughout this
section, let $R$ be a commutative ring. A power structure on $R$ is a map
$\varphi\colon R\times (1+tR[[t]]) \to (1+tR[[t]]); (r, f(t)) \mapsto f(t)^r$
satisfying certain properties, see
\cite[354]{Gusein-Zade-Luengo-Melle-HernandezHS} for a precise definition. In
this paper, we only consider finitely determined power structures in the sense
of \cite[354]{Gusein-Zade-Luengo-Melle-HernandezHS}. In this case, a power
structure is completely determined by the terms $(1-t)^{-r} =
\sum_{i=0}^{\infty} a_i(r)t^i$ for $r \in R$ subject to some conditions. Since
this is the only case we are interested in, we state the precise definition here.

\begin{definition}[{\cite[Corollary~2.5]{Pajwani-PalPS}}]
  \label{def:finite-ps-char}
  A \emph{power structure} on $R$ is a collection of functions $a_n\colon R \to R$
  for $n \in \Z_{\ge 0}$, denoted $a_\ast$, satisfying the conditions:
  \begin{itemize}
  \item $a_0 = 1$,
  \item $a_1 = \id_R$,
  \item $a_i(0) = 0$ for all $i \ge 1$,
  \item $a_i(1) = 1$ for all $i \ge 1$, and
  \item $a_n(r+s) = \sum_{i=0}^n a_i(r)\cdot a_{n-i}(s)$ for all $r,s \in R$
    and $n \ge 0$.
  \end{itemize}
\end{definition}

\begin{definition}[{\cite[Definition~2.6 and Proposition~2.7]{Pajwani-PalPS}}]
  A ring homomorphism $\varphi\colon R\to S$
  between rings equipped with power structures $a^R_\ast$ and $a^S_\ast$
  \emph{respects the power structures} if 
  we have
  $\varphi(a_n^R(r)) = a_n^S(\varphi(r))$ for all $r
  \in R$ and $n \ge 0$.
\end{definition}

\begin{lemma}
  \label{lem:power-structre-large-sum}
  Let $a_\ast$ be a power structure on a ring $R$ and let $r_1, \dots, r_m \in
  R$. Then we have for $n \ge 0$
  \[
    a_n\left(\sum_{i=1}^m r_i\right) = \sum_{i_1+ \dots+ i_m = n} a_{i_1}(r_1)\cdots a_{i_m}(r_m).
  \]

  \begin{proof}
    Perform an induction on $n$ and use the sum relation $a_n(r+r') =
    \sum_{i=0}^na_i(r)a_{n-i}(r')$.
  \end{proof}
\end{lemma}

\begin{example}
  McGarraghy \cite[Section 4]{McGarraghySP} constructs so-called
  non-factorial
  symmetric powers on $\GW(k)$. These assemble into a power structure $a_{\ast}^{\text{NF}}$ essentially by
  \cite[Propostion 4.23]{McGarraghySP} and this power structure is determined
  by $a_n^{\text{NF}}(\langle \alpha\rangle) = \langle \alpha^n\rangle$.
\end{example}

\begin{example}
  There is a power structure on $K^\rs_0(\Var_k)$ determined by $a_n^{\Sym}([X]) =
  [\Sym^nX]$ for quasi-projective $k$-varieties $X$. Note that the $a_i^{\Sym}$ are the
  coefficients of the Kapranov zeta function, see
  \cite[(1.3)]{kapranov2000ECDT} or \cite[Section 1.4]{MustataZF}. In characteristic zero this
  defines a power structure on $K_0(\Var_k)$ since there $K_0(\Var_k) =
  K^\rs(\Var_k)$. Gusein-Zade,
  Luengo, and Melle-Hernández
  \cite{Gusein-Zade-Luengo-Melle-HernandezPS,Gusein-Zade-Luengo-Melle-HernandezHS}
  first note that this defines a power structure. A more detailed account can be found in
  \cite[Section 6]{mckean2024symmetricpowersnullmotivic}.
\end{example}

Since the operation of taking symmetric powers forms a power structure on $K^\rs(\Var_k)$, and the
Grothendieck-Witt ring also admits power structures, one can wonder
whether the compactly supported $\A^1$-Euler characteristic is compatible with
the power structure $a_\ast^{\Sym}$ on $K_0^\rs(\Var_k)$ and some power structure on
$\GW(k)$. This question was asked by McKean and Nicaise to Pajwani and Pál and
they prove the following result.

\begin{theorem}[{\cite[Corollary 3.27]{Pajwani-PalPS}}]
  \label{thm:pp-ps}
  There is a power structure, $a_\ast^\PPal$, on $\GW(k)$ such that
  \[
    a_n^\PPal(\langle \alpha\rangle) = \langle \alpha^n\rangle +
    \frac{n(n-1)}{2}t_\alpha
  \]
  for $\alpha \in k^\times$. Here $t_\alpha = \langle 2\rangle + \langle
  \alpha\rangle - \langle 1\rangle - \langle 2\alpha\rangle$.
  
  Moreover, if there exists a power structure on $\GW(k)$ such that
  $\chi_c\colon K_0^\rs(\Var_k) \to \GW(k)$ respects the power structures, then
  this power structure is the Pajwani-Pál power structure $a_\ast^\PPal$.
\end{theorem}

\begin{rem}
  In the form presented here,
  Pajwani-Pál only prove Theorem~\ref{thm:pp-ps} in characteristic zero.
  In positive characteristic, they prove the assertion for the restriction of $\chi_c$ to the subring of zero-dimensional varieties $K_0(\Et)$ in $K_0(\Var_k)$, since there $a_\ast^{\Sym}$ forms a power structure without inverting universal homeomorphisms. Since $K_0(\Et) \subset K_0^\rs(\Var_k)$, their statement also holds as stated here.
\end{rem}

\begin{rem}
  This power structure $a^\PPal_\ast$ agrees with McGarraghy's non-factorial
  symmetric power structure $a^{\text{NF}}_\ast$
  if $2$ is a square in $k$. More generally, $a_n^\PPal(\langle \alpha\rangle) =
  a_n^{\text{NF}}(\langle \alpha\rangle)$
  if and only if $[2]\cup [\alpha] = 0 \in H^2(k,\Z/2\Z)$ by \cite[Lemma
  3.29]{Pajwani-PalPS}. That is, $a_\ast^\PPal$ agrees with $a_\ast^{\text{NF}}$ on the subgroup
  generated by all $\langle \alpha
  \rangle$ satisfying $a_n^\PPal(\langle \alpha \rangle) = \langle
  \alpha^n\rangle$. In particular, the
  quadratic forms $\langle 1 \rangle$ and $\langle -1 \rangle$ are in this
  subgroup.
\end{rem}

\begin{rem}[{\cite[Corollary 3.19]{Pajwani-PalPS}}]
  By \cite[Corollary~3.18]{Pajwani-PalPS}, the element $t_\alpha$ is two-torsion in
  $\GW(k)$. Thus, the definition of
  $a_\ast^{\text{PP}}$ can be rewritten as
  \[
    a_n^{\text{PP}}(\langle \alpha\rangle) =
    \begin{cases}
      \langle \alpha^n\rangle, & n \equiv 0,1\mod 4,\\
      \langle \alpha^n\rangle +
      t_\alpha, & n \equiv 2,3 \mod 4.
    \end{cases}
  \]
\end{rem}

It is still an open question whether $\chi_c$ respects the power structures in
general. There are, however, certain classes of schemes for which $\chi_c$ respects
the power structures.

\begin{definition}[{\cite[Definition 4.1]{PajwaniRohrbachViergever}}]
  A variety $X$ is \emph{symmetrisable} if $\chi_c(\Sym^n(X)) =
  a_n^\PPal(\chi_c(X))$ for all $n \ge 0$.
\end{definition}

As mentioned in the introduction there are already some examples of symmetrisable varieties.

\begin{rem}
  Pajwani-Rohrbach-Viergever only use the concrete definition of
  $a_\ast^\PPal$ for checking that $\A^n$ is symmetrisable, for which the term
  $t_\alpha$ vanishes. Since $t_\alpha \in \GW(k)$ is torsion, the term
  $t_\alpha$ also does not show up in the compatibility result of Bejleri-McKean
  \cite[Theorem 9.4]{mckean2024symmetricpowersnullmotivic}. The
  compatibility results proven in this work are the first compatibility result since \cite{Pajwani-PalPS} actually
  calculating with a possibly non-vanishing $t_\alpha$.
\end{rem}

\section{$\A^1$-Euler Characteristic of a Split Toric Variety and its
  Symmetric Powers}
\label{sec:toric}
We shall now compute the $\A^1$-Euler characteristic of a split toric variety
and its symmetric powers. Our method was heavily inspired by the computation of
the topological Euler characteristic of a complex toric variety as found in
\cite{MathOverflowToricEuler}.

Demazure \cite{Demazure1970Toric} first constructs toric varieties. We
briefly recall the construction of split toric varieties following
Elizondo, Lima-Filho, Sottile, and Teitler \cite[Section~2.1]{ELFST2014ArithmToric} in order to state the results.
In the following, let $N$ be a finitely generated free Abelian group of rank
$n$ and let $M :=  \Hom_\Z(N,\Z)$ be its dual.

\begin{definition}[{\cite[218]{ELFST2014ArithmToric}}]
  Let $\sigma$ be a finitely generated
  subsemigroup of $N$. The \emph{polar} $\sigma^\vee$ of $\sigma$ is the
  finitely generated subsemigroup
  \[
    \sigma^\vee = \{u \in M \mid u(v)\ge 0\text{ for all } v\in \sigma\}
  \]
  of $M$.

  A \emph{cone} is a saturated finitely generated subsemigroup $\sigma$, that is, $(\sigma^\vee)^\vee = \sigma$.

  A \emph{face} of a cone $\sigma$ is a subsemigroup of $\sigma$ such that
  there exists $u \in \sigma^\vee$ with
  \[
    \tau = \{v \in \sigma \mid u(v) = 0\}.
  \]

  A cone $\sigma$ is \emph{pointed} if $0$ is a face of $\sigma$.
\end{definition}

\begin{rem}
  Let $\sigma$ be a cone.
  In the following, we will write $\tau \preceq \sigma$ to indicate that
  $\tau$ is a face of $\sigma$.

  Also, note that $\sigma$ is a face of
  $\sigma$ since $\sigma$ is the face associated with
  $u = 0 \in \sigma^\vee$.
\end{rem}

\begin{definition}[{\cite[218]{ELFST2014ArithmToric}}]
  Let $\sigma$ be a pointed cone in $N$. We can associate $\sigma$ with the
  affine scheme $V_\sigma := \Spec k[\sigma^\vee]$ of the semigroup algebra
  over $k$ generated by $\sigma^\vee$.
\end{definition}

\begin{rem}[{\cite[218]{ELFST2014ArithmToric}}]
  Let $\tau\subset \sigma$ be a face of a pointed cone in $N$. Hence, the
  inclusion $\sigma^\vee \subset \tau^\vee$ induces an open immersion $V_\tau
  \hookrightarrow V_\sigma$, since $k[\tau^\vee]$ is a subring of the field of
  fractions of $k[\sigma^\vee]$.
\end{rem}

\begin{definition}[{\cite[218]{ELFST2014ArithmToric}}]
  A \emph{fan} $\Sigma$ of $N$ is a finite collection of pointed cones in $N$ such
  that
  \begin{enumerate}
  \item any face of a cone in $\Sigma$ is also a cone in $\Sigma$ and\label{def:fan:face}
  \item the intersection of any two cones in $\Sigma$ is a common face, which
    by (\ref{def:fan:face}) is again a cone in $\Sigma$.\qedhere
  \end{enumerate}
\end{definition}

\begin{definition}[{\cite[218]{ELFST2014ArithmToric}}]
  Given a fan $\Sigma$ in $N$, we construct its associated \emph{split toric
    variety} $X_\Sigma$ as $X_\Sigma = \colim_{\sigma\in\Sigma}V_\sigma$. That
  is, we glue the affine schemes $V_\sigma$ for
  $\sigma$ a cone in $\Sigma$ along their common subschemes corresponding with
  smaller cones in $\Sigma$.
\end{definition}

Note that since all cones in $\Sigma$ are pointed, they contain the cone $0$
as a face. Thus, every affine $V_\sigma$ contains the affine scheme $V_0$.

\begin{rem}[{\cite[218]{ELFST2014ArithmToric}}]
  We have $0^\vee = M$ and thus, we have $V_0 = \Spec k[M] \cong \G_m^{\rank
    N}$, which is a split torus.
\end{rem}

\begin{rem}[{\cite[218]{ELFST2014ArithmToric}}]
  Let $\Sigma$ be a fan in $N$ and
  let $T = \Spec k[M] = V_0 \cong \G_m^{\rank N}$. The inclusion
  $k[\sigma^\vee] \to k[M]$ for a cone $\sigma \in \Sigma$ induces a faithful
  action
  \[
    T \times_{\Spec k} V_\sigma \to V_\sigma
  \]
  via the map $k[\sigma^\vee] \to k[M]\otimes k[\sigma^\vee]; u \mapsto
  u\otimes u$. This action is compatible with the inclusion $V_\tau\to
  V_\sigma$ for $\tau$ a face of $\sigma$, and thus, the $T$-actions on the
  $V_\sigma$ assemble into a faithful $T$-action on $X_\Sigma$.
\end{rem}

In order to compute the compactly supported $\A^1$-Euler characteristic of a
split toric variety as in \cite{MathOverflowToricEuler}, we need to decompose
a split toric variety into split tori based on the cones in the fan
$\Sigma$. This is generally known as the orbit-cone correspondence. We briefly
recall its statement here for the convenience of the reader. See Fulton
\cite[Section~3.1]{Fulton1993IntroToric} for details.

Let $\sigma$ be a cone in
$\Sigma$. Let $N_\sigma \subset N$ be the sublattice generated by $\sigma$ and
let $N(\sigma) = N/N_\sigma$. Let
\[
  M(\sigma) = \{u \in M \mid u(v) = 0\text{
    for all } v \in \sigma\}
\]
thus, $M(\sigma) = N(\sigma)^\vee$. Define
$O_\sigma = \Spec k[M(\sigma)]$ and $V_\sigma(\tau) = \Spec k[\sigma^\vee\cap
M(\tau)]$ for $\tau \preceq \sigma$ a face.

For $\sigma \in \Sigma$ and $\tau \preceq \sigma$ a face, we have the
surjective ring homomorphisms $i^\ast_{\tau\preceq \sigma}\colon
k[\sigma^\vee] \to k[\sigma^\vee\cap M(\tau)]$ and $i_\tau^\ast\colon k[\tau^\vee] \to
k[M(\tau)]$ induced by the assignments
\[
  i^\ast_{\tau\preceq \sigma}(u) =
  \begin{cases}
    u \in k[\sigma^\vee\cap M(\tau)] & \text{if } u\in M(\tau)\cap \sigma^\vee,\\
    0 \in k[\sigma^\vee\cap M(\tau)] & \text{otherwise},
  \end{cases}
\]
for $u \in \sigma^\vee$ and
\[
  i^\ast_{\tau}(u) =
  \begin{cases}
    u \in k[M(\tau)] & \text{if } u\in M(\tau),\\
    0 \in k[M(\tau)] & \text{otherwise}.
  \end{cases}
\]
for $u \in \tau^\vee$.
These ring homomorphisms define closed immersions $i_{\tau\preceq \sigma}\colon V_\sigma(\tau)
\to V_\sigma$ and $i_\tau\colon O_\tau \to V_\tau$.

Furthermore, the inclusions
$\sigma^\vee\cap M(\tau) \subset M(\tau)$ and $\sigma^\vee \subset \tau^\vee$
induce inclusions $j^\ast_{\tau\preceq \sigma}\colon
k[\sigma^\vee \cap M(\tau)] \to k[M(\tau)]$ and $j^\ast_{\sigma\succeq
  \tau}\colon k[\sigma^\vee]\to k[\tau^\vee]$.
The
induced morphisms of schemes $j_{\tau\preceq\sigma}\colon O_\tau \to
V_\sigma(\tau)$ and $j_{\sigma\succeq\tau}\colon
V_\tau\to V_\sigma$ are open immersions.

These morphisms induce the commutative diagram of schemes
\[
  \begin{tikzcd}
    O_\tau\ar[r, "{j_{\tau\preceq\sigma}}"]\ar[d, "{i_\tau}"]& V_\sigma(\tau)\ar[d, "{i_{\tau\preceq\sigma}}"]\\
    V_\tau\ar[r, "{j_{\sigma\succeq\tau}}"]&V_\sigma.
  \end{tikzcd}
\]
This realises $O_\tau$ as a locally closed
subset of $V_\sigma$. Via the description $X_\Sigma =
\colim_{\sigma\in\Sigma} V_\sigma$ this also represents $O_\sigma$ as a
locally closed subset of $X_\Sigma$. In addition $O_\sigma$ is a $T$-orbit,
isomorphic to $\G_m^{\rank N-\rank\sigma}$.

\begin{proposition}[Orbit-Cone Correspondence {\cite[Proposition on p. 54]{Fulton1993IntroToric}}]
  \label{prop:orbit-cone}
  Let $\Sigma$ be a fan in $N$. Then for every cone $\sigma$ in
  $\Sigma$, we have a decomposition of $V_\sigma$ into locally closed subsets
  \[
    V_\sigma = \coprod_{\tau}O_\tau
  \]
  where $\tau$ ranges the faces of
  $\sigma$.

  In particular, we get a decomposition
  \[
    X_\Sigma = \coprod_{\sigma \in \Sigma}O_{\sigma}
  \]
  of $X_\Sigma$ into locally closed subsets.
\end{proposition}

\begin{rem}
  Note that \cite[Proposition on p.~54]{Fulton1993IntroToric} is only stated
  for $k = \C$, but the proof also works for arbitrary fields.
\end{rem}

\begin{lemma}
  \label{lem:half-binom}
  Let $n \ge 1$ be an integer, then we have
  \[
    \sum_{i=0}^{\lfloor n/2\rfloor} \binom{n}{2i} = \sum_{i=0}^{\lfloor
      n/2\rfloor} \binom{n}{2i+1} = 2^{n-1}.
  \]

  \begin{proof}
    First note that we have
    \[
      2^n = (1 + 1)^n = \sum_{i=0}^n\binom ni = \sum_{i=0}^{\lfloor n/2\rfloor} \binom{n}{2i} + \sum_{i=0}^{\lfloor
        n/2\rfloor} \binom{n}{2i+1}.
    \]
    Thus, it is enough to show one of the equalities.

    We start with $n = 2m+1$ odd. Thus, we have $\binom ni =
    \binom{n}{n-i}$ and
    \[
      \sum_{i=0}^{\lfloor n/2\rfloor} \binom{n}{2i} = \sum_{i=0}^{\lfloor
        n/2\rfloor} \binom{n}{n-2i} = \sum_{i=0}^{\lfloor n/2\rfloor} \binom{n}{2i+1},
    \]
    where we re-indexed the sum in the last equation, replacing $i$ by
    $m-i$ and noting $m = \lfloor n/2\rfloor$. This completes the proof in case $n$ is odd.

    In case $n = 2m$ is even, we have by the odd case above
    \[
      \sum_{i=0}^{\lfloor n/2\rfloor} \binom{n}{2i} = \sum_{i=0}^{\lfloor
        n/2\rfloor} \binom{n-1}{2i-1} + \binom{n-1}{2i} = \underset{=2^{n-2}}{\underbrace{\sum_{i=0}^{\lfloor
            n/2\rfloor} \binom{n-1}{2i-1}}} + \underset{=2^{n-2}}{\underbrace{\sum_{i=0}^{\lfloor
            n/2\rfloor} \binom{n-1}{2i}}} = 2\cdot 2^{n-2} = 2^{n-1}.
    \]
    This completes the proof in the even case.
  \end{proof}
\end{lemma}

\begin{lemma}
  \label{lem:gm-euler-char}
  We have for $n \ge 1$
  \[
    \chi_c(\G_m^n/k) = (-2)^{n-1}(\langle -1\rangle - \langle 1\rangle).
  \]

  \begin{proof}
    By Example \ref{ex:chi-affine-line-gm}, we have
    $\chi_c(\G_m/k) = \langle-1\rangle - \langle 1\rangle$.
    This yields using Lemma \ref{lem:half-binom}
    \begin{align*}
      \chi_c(\G_m^n/k)
      &= (\langle -1\rangle - \langle 1 \rangle)^n\\
      &= \sum_{i = 0}^n \binom ni (-1)^{n-i}\langle(-1)^i\rangle\\
      &= (-1)^{n-1}\underset{ = 2^{n-1}}{\underbrace{\sum_{i=0}^{\lfloor n/2\rfloor} \binom{n}{2i+1}}} \langle
        -1\rangle - (-1)^{n-1}\underset{=2^{n-1}}{\underbrace{\sum_{i=0}^{\lfloor n/2\rfloor} \binom{n}{2i}}} \langle
        1\rangle\\
      &= (-2)^{n-1}(\langle-1\rangle - \langle 1\rangle).\qedhere
    \end{align*}
  \end{proof}
\end{lemma}

\begin{proposition}
  \label{prop:euler-char-toric}
  Let $\Sigma$ be a fan in a rank $n$ lattice $N$. Let $X_\Sigma$ be the
  associated split toric variety over $k$ and let $\Sigma_i$ be the collection
  of cones in $\Sigma$ of rank $n-i$. Then we have
  \[
    \chi_c(X_\Sigma/k) = |\Sigma_0|\cdot \langle 1\rangle + \sum_{i\ge 1} |\Sigma_i| \cdot
    (-2)^{i-1}(\langle-1\rangle - \langle 1\rangle) \in \GW(k).
  \]

  \begin{proof}
    The fan $\Sigma$ yields a decomposition of $X_\Sigma$ into a disjoint
    union of locally closed subschemes
    \[
      X_\Sigma= \coprod_{\sigma\in \Sigma} \G_m^{n-\rank\sigma}
    \]
    by Proposition \ref{prop:orbit-cone}.
    In this decomposition, there are $|\Sigma_i|$ copies of
    $\G_m^i$. Thus, the desired formula follows from the scissors relation and
   Lemma \ref{lem:gm-euler-char}.
  \end{proof}
\end{proposition}

Using Proposition \ref{prop:orbit-cone}, we can also show that toric varieties
are symmetrisable.

\begin{proposition}
  \label{prop:split-toric-symmetrisable}
  Let $X = X_\Sigma$ be a split toric variety over $k$. Then $X$ is
  symmetrisable.

  \begin{proof}
    By Proposition \ref{prop:orbit-cone}, we have
    \[
      [X] = \sum_{i=0}^N \nu_i[\G_m^i] \in K_0(\Var_k)
    \]
    for some $\nu_i \in \Z_{\ge 0}$. Now, $\G_m^i$ is $K_0$-étale linear in the
    sense of \cite[Definition~3.1]{PajwaniRohrbachViergever}, since $[\G_m] = [\A^1]-[\Spec k] \in
    K_0(\Var_k)$. Since $K_0$-étale linear varieties form a ring by \cite[Definition~3.1]{PajwaniRohrbachViergever} and are
    symmetrisable by \cite[Theorem 4.10]{PajwaniRohrbachViergever}, we get that $X$ is symmetrisable.
  \end{proof}
\end{proposition}

\section{Results on the Smoothness of Quotients}
\label{sec:smoothness-results}
Throughout this section, $G$ denotes a finite group that is tame over $k$;
that is, $G$ is a finite group with order prime to
the characteristic of $k$.
For computing the $\A^1$-Euler characteristic of a quotient by a group action,
we need some basic results on quotients and a criterion to detect that the quotient is smooth. In this section, we
collect the necessary results and provide proofs of these well-known theorems
for the convenience of the reader
whenever we were not able to find a good reference. We want to note that
a proof of Proposition
\ref{prop:quotient-proj} has also been given by Viehweg
\cite[Corollary~3.46]{Viehweg1995Moduli} over algebraically closed fields.

\begin{proposition}
  \label{prop:quotient-proj}
  If $G$ is a finite group acting on a quasi-projective $k$-scheme $X$, then
  the quotient $X/G$ exists as a scheme and is again quasi-projective. Furthermore if $X$ is
  projective, then so is $X/G$.

  \begin{proof}
    The quotient $X/G$ exists as a scheme by Serre \cite[Chapter~III, Section~12,
    Example~1 after Proposition~19]{SerreAGCF}. Furthermore,
    the quotient map is finite and $X/G$ is normal. Thus, \cite[Proposition~6.6.1]{EGA2} yields that $X/G$ has an ample line bundle and hence, is
    quasi-projective.

    Note that $X/G$ is separated since $X/G$ is quasi-projective. Hence, $X/G$ is proper by
    \cite[\href{https://stacks.math.columbia.edu/tag/01W6}{Tag
      01W6}]{stacks-project} if $X$ is projective. Since a proper,
    quasi-projective scheme is projective by \cite[\href{https://stacks.math.columbia.edu/tag/0BCL}{Tag 0BCL}]{stacks-project}, we get that $X/G$ is
    projective if $X$ is projective.
  \end{proof}
\end{proposition}

The following is a special case of much more general results; see, for example
\cite[Theorem~4.16]{MoonenAV} for a general version of part
(\ref{lem:QuotFlatBaseChange:quot}).

\begin{lemma}\label{lem:QuotFlatBaseChange}
  Let $k\subset K$ be an extension of fields. Let $G$ be a finite group
  acting on a quasi-projective $k$-scheme $X$, inducing the action of $G$ on
  the  quasi-projective $K$-scheme $X_K:=X\times_{\Spec k}\Spec K$ by
  letting $G$ act trivially on $\Spec K$.
  \begin{enumerate}
  \item There is a canonical isomorphism
    \[
      X_K/G\cong (X/G)_K.
    \]
    \label{lem:QuotFlatBaseChange:quot}
  \item There is a canonical isomorphism
    \[
      (X_K)^G\cong (X^G)_K.
    \]
    \label{lem:QuotFlatBaseChange:fixed}
  \end{enumerate}

  \begin{proof}
    As above, we can reduce to the case that $X = \Spec A$ is affine.
    Let $k_G$ be
    the ring of functions on $G$, considered as the disjoint union
    $\amalg_{g\in G}\Spec k$, that is, $k_G=\prod_{g\in G}k$. The
    multiplication on $G$ equips the commutative ring $k_G$ with the structure
    of a commutative
    Hopf algebra over $k$, by identifying $k_G\otimes_kk_G$ with $k_{G\times
      G}$ via $\prod_{g\in G}a_g\otimes \prod_{g'\in G}b_{g'}\mapsto
    \prod_{(g,g')\in G\times G}a_gb_{g'}$. The action $\rho\colon
    G\times_k\Spec A\to \Spec A$ of $G$ on $\Spec A$ equips $A$ with the structure of a co-module
    over $k_G$, via $\rho^\ast\colon A\to k_G\otimes_k A$. We also have the
    map $i\colon A\to  k_G\otimes_k A, i(a)=1\otimes a$. The sub-$k$-algebra of invariants
    $A^G$ is the equalizer of $\rho^\ast$ and $i$.

    Similarly, we have a description of the sub-$K$-algebra of invariant
    $(A\otimes_kK)^G$ as the equalizer of $\rho_K^\ast\colon A\otimes_kK\to
    K_G\otimes_K(A\otimes_kK)$ and $i_K\colon A\otimes_kK\to
    K_G\otimes_K(A\otimes_kK)$; here $\rho_K\colon G\times_{\Spec K}(\Spec
    A\otimes_kK)\to \Spec A\otimes_kK$ is the base-extension of $\rho$. Via the isomorphism 
    $K_G\otimes_K(A\otimes_kK)\cong (k_G\otimes_k A)\otimes_kK$, we obtain
    $\rho_K^\ast=\rho^\ast\otimes\id_K$ and $i_K=i\otimes\id_K$. Since
    $k\subset K$ is a flat extension, this induces the isomorphism
    $(A\otimes_kK)^G\cong A^G\otimes_kK$. This proves (\ref{lem:QuotFlatBaseChange:quot}).

    The proof of (\ref{lem:QuotFlatBaseChange:fixed}) is similar: we have the Cartesian square
    \[
      \begin{tikzcd}
        X^G\ar[r, "i"]\ar[d, "i"]&X\ar[d, "\rho^\vee"]\\
        X\ar[r, "\delta"]&\prod^k_{g\in G}X
      \end{tikzcd}
    \]
    where  $\prod^k_{g\in G}$ means the product over $\Spec k$,  $\delta$ is
    the diagonal map, and $\rho^\vee$ is the map $\rho^\vee(x)=\prod_{g\in
      G}\rho(g)(x)$. We have the corresponding Cartesian square  describing
    $(X_K)^G$. As we have canonical isomorphisms $(\prod^k_{g\in G}X)_K\cong
    \prod^K_{g\in G}(X_K)$, transforming $\rho^\vee\times_{\Spec k}\id_{\Spec
      K}$ to $\rho_K^\vee$ and similarly for $\delta$, we have 
    \[
      (X^G)_K=(X\times_{\prod^k_{g\in G}X}X)_K\cong X_K\times_{\prod^K_{g\in G}X_K}X_K=(X_K)^G.\qedhere
    \]
  \end{proof}
\end{lemma}

\begin{lemma}
  Let $X$ be a smooth, quasi-projective scheme over $k$ and let $G$ be a
  finite group acting on $X$. Assume that the order of $G$ is prime to
  the characteristic of $k$. Then $X^G$ is again smooth.

  \begin{proof}
    By Lemma~\ref{lem:QuotFlatBaseChange}~(\ref{lem:QuotFlatBaseChange:fixed}),
    we can check smoothness after passing to an algebraic closure. Thus, we
    may assume that $k$ is algebraically closed.
    Let $x \in X^G$. There exists a open affine subscheme of $X$
    that is stable under the $G$-action. Since we can check smoothness
    locally, we can thus assume that $X$
    is affine.

    Note that, Luna \cite[Lemme in Section III.1, p. 96]{LunaSlice} only uses
    that the acting group acts semi-simply on the maximal ideal of interest,
    which is the case for our $G$-action on $X$. Thus, \cite[Lemme in Section
    III.1, p. 96]{LunaSlice} implies that there is a
    $G$-equivariant morphism $\varphi\colon X \to T_xX$, étale at $x$, sending $x$ to $0$
    where we equip $T_xX$ with the induced $G$-action. Since $\varphi$
    restricts to a morphism $X^G \to (T_xX)^G$, étale at $x$, it is enough to show
    that $(T_xX)^G$ is smooth. But the action on $T_xX$ is a $k$-linear action
    on a $k$-vector space.

    Note that if $H$ is a cyclic group acting linearly on $\A^N$, then $H$
    decomposes $\A^N$ via characters into linear subspaces. In particular, $(\A^N)^H$ is the linear
    subspace given by the trivial character. In particular if $G$ acts
    linearly on $\A^N$, then $(\A^N)^G = \bigcap_{g \in G}(\A^N)^{\langle
      g\rangle}$ is an intersection of linear subspaces and hence, smooth. Here, $\langle
    g\rangle$ denotes the cyclic group generated by $g \in G$.

    If we apply this observation to $(T_xX)^G$, we get that $(T_xX)^G$ is smooth.
  \end{proof}
\end{lemma}

\begin{lemma}
  \label{lem:normal-bundle-no-fixed}
  Let $X$ be a smooth scheme over $k$ and let $G$ act on $X$. Then the induced
  $G$-action on the normal bundle $N_{X^G}X$ of $X^G$ in $X$ has no non-trivial fixed
  points.
  \begin{proof}
    Let $i\colon X^G \to X$ be the inclusion.
    By definition, we have $N_{X^G}X = i^\ast TX/T(X^G)$. A tangent vector in a
    point $x \in X$ is the same as a morphism
    $\Spec k(x)[\varepsilon]/(\varepsilon^2) \to X$ with $x$ as its
    image. Thus, we
    see $T(X^G) = (i^\ast TX)^G$. In particular, $(N_{X^G}X)^G = (i^\ast
    TX/(i^\ast TX)^G)^G = 0$.
  \end{proof}
\end{lemma}

\begin{theorem}
  \label{thm:divisor-smooth-quotient}
  Let $X$ be a smooth, quasi-projective scheme over $k$ and let $G = \Z/p\Z$ be a cyclic group of prime
  order $p\ne \characteristic k$. Assume that $D = X^G \subset X$ is a divisor. Then
  $X/G$ is smooth.

  \begin{proof}
    We may check the smoothness of $X/G$ after passing to an algebraic
    closure of $k$. Thus, we may assume that $k$ is algebraically closed. Since
    $G$ has prime order, $G$ acts freely on $X\setminus D$ and thus,
    $(X/G)\setminus D$ is smooth since the projection map $\pi\colon X \to
    X/G$ is étale away from $D$.

    We are therefore left to show that $X/G$ is smooth in every point $x \in
    D$. Since $k$ is algebraically closed and hence, perfect, we thus have to
    show that $\OO_{X/G,x}$ is regular for
    every $x \in D$.

    Let $A = \OO_{X,x}$ and let $\m_{X,x} = (t, x_1, \dots, x_n) \subset A$ be
    the maximal ideal such that $t$ is a local equation for $D$, so
    $\OO_{D,x} = A/(t)$, and such that $\dim A = n+1$. We can find such generators
    since $X$ and $D$ are smooth. Since $G$ acts trivially on $D$, we have that $G$
    acts trivially on $\m_{D,x}/\m_{D,x}^2$ where $\m_{D,x} = (x_1, \dots,
    x_n)\subset A/(t)$ is the maximal ideal. Since taking $G$-invariants is
    exact, we can lift generators of $\m_{D,x}/\m_{D,x}^2$ to
    elements of $(\m_{X,x}/\m_{X,x}^2)^G$. By the same argument, we can lift these
    elements to elements $x_1', \dots, x_n' \in A$ that are invariant under the
    $G$-action and such that $\m_{X,x} = (t, x_1', \dots, x'_n)$.

    The group $G$ acts non-trivially on $(t)/(t^2)$ by Lemma \ref{lem:normal-bundle-no-fixed} since $(t)/(t^2)$ is isomorphic to
    the fibre of the dual of the normal bundle of $D$ in $X$. Thus, a generator $g \in G$
    acts by multiplication by a
    primitive $p$-th root of unity $\zeta \in k$.
    Since $k[G]$ is semisimple, the exact sequence of $k[G]$-modules
    \[
      0 \to (t^2) \to (t) \to (t)/(t^2) \to 0
    \]
    splits, and thus, we can find an element $t' \equiv t \mod (t^2)$ with
    $g\cdot t' = \zeta t'$.

    In particular, we have $\m_{X,x} = (t', x'_1, \dots, x'_n)$. By completing
    $A$, we get $\hat A = k[[t', x'_1, \dots, x'_n]]$ using
    \cite[\href{https://stacks.math.columbia.edu/tag/032A}{Tag 032A} and
    \href{https://stacks.math.columbia.edu/tag/0C0S}{Tag
      0C0S}]{stacks-project}. Furthermore, $g$ acts on $\hat
    A$ by $g\cdot t' = \zeta t'$ and $g\cdot x'_i = x'_i$, and thus,
    \[
      \widehat{(A^G)} = (\hat A)^G = \left\{\sum_{i \ge 0}f_it^{pi}\mid f_i \in
      k[[x'_1, \dots, x'_n]]\right\} \cong k[[t^p, x'_1, \dots, x'_n]]
    \]
    is regular. Hence, $\OO_{X/G,x} = A^G$ is regular by \cite[(24.D),
    p.~175]{MatsumuraCA}.
  \end{proof}
\end{theorem}

\begin{lemma}
  \label{lem:z2-normal-action}
  Let $X$ be a smooth scheme over $k$ and let $G=\Z/2\Z$ act on $X$. Then $G$ acts by
  multiplication with $-1$ on the normal bundle $N_{X^G}X$ of $X^G$ in
  $X$.

  \begin{proof}
    Locally on $X^G$, we get a decomposition $N_{X^G}X = (N_{X^G}X)_1
    \oplus (N_{X^G}X)_{-1}$ into isotypical components. Here $G$ acts by
    multiplication with $\pm 1$ on $(N_{X^G}X)_{\pm 1}$. Since $N_{X^G}X$ has
    trivial $G$-invariants by Lemma \ref{lem:normal-bundle-no-fixed}, we get
    $N_{X^G}X = (N_{X^G}X)_{-1}$, which yields the statement.
  \end{proof}
\end{lemma}

\begin{corollary}
  \label{cor:smooth-z2-blow-up-quotient}
  Let $X$ be a smooth, quasi-projective scheme over $k$ and let $G = \Z/2\Z$ act on $X$. Extend the $G$
  action on $X$ to the blow-up $\Bl_{X^G}X$ of $X$ in $X^G$. Then
  $(\Bl_{X^G}X)/G$ is smooth.

  \begin{proof}
    Let $E \subset \Bl_{X^G}X$ be the exceptional divisor. Hence, we have a
    $G$-equivariant isomorphism $E \cong \PP(N_{X^G}X)$. By Lemma
    \ref{lem:z2-normal-action}, $G$ acts by multiplication with $-1$ on
    $N_{X^G}X$. Thus, $G$ acts trivially on $E = \PP(N_{X^G}X)$, and hence, $E
    \subset (\Bl_{X^G}X)^G$. Since we are
    blowing up the fixed point locus of the $G$-action on $X$, we also have
    $(\Bl_{X^G}X)^G \subset E$ and thus, $E = (\Bl_{X^G}X)^G$. Hence,
    $(\Bl_{X^G}X)/G$ is smooth by Theorem \ref{thm:divisor-smooth-quotient}.
  \end{proof}
\end{corollary}

\section{The Equivariant Quadratic Euler Characteristic}
\label{sec:equivariant-Euler-char}
Throughout this section, $G$ denotes a finite group that is tame over $k$;
that is, $G$ is a finite group with order prime to
the characteristic of $k$.

\subsection{Construction}
Pajwani and Pál \cite[Section~3]{Pajwani-PalYZ} first introduce an equivariant
quadratic Euler characteristic using de Rham
cohomology. We will give a
different and more natural construction here using Hodge cohomology and the language of
Grothendieck-Witt rings of exact categories, following Schlichting
\cite{schlichting_hermitian_2010}. Grothendieck-Witt rings of exact categories
have been first constructed by Quebbemann-Scharlau-Schulte
\cite[pp.~280f.]{QSS1979QuadraticCat}. The construction here is equivalent to
Pajwani-Pál's construction in characteristic zero, see Proposition~\ref{prop:equiv-euler-char-equal-pp}.

An \emph{exact category} is an additive category $\mathcal E$ together with a
collection of short exact sequences
\begin{equation}
  \label{eq:ditinguished-ses}
  \begin{tikzcd}
    0 \ar[r] & A \ar[r, "f"] & B\ar[r, "g"] & C \ar[r] & 0
  \end{tikzcd}
\end{equation}
satisfying certain properties. See Quillen \cite[First definition in Section~2;
p.~100]{Quillen1973HAlgebraicK}  for a precise definition.

A morphism $\varphi\colon A \to B$ in $\mathcal E$ is called an \emph{inflation} if $\varphi$ shows up as
$f$ in \eqref{eq:ditinguished-ses} in some exact sequence in $\mathcal E$. If
$\varphi$ shows up as $g$ in \eqref{eq:ditinguished-ses} in some exact
sequence in $\mathcal E$, then $\varphi$ is called a \emph{conflation}. An
exact category $\mathcal E$ in which every short exact sequence splits is
called a \emph{split exact category}.

\begin{example}
  \label{ex:rep-split-exact}
  Let $\Rep_k(G)$ be the category of $G$-representations in finite-dimensional
  $k$-vector spaces, precisely, $\Rep_k(G)$ is the category of left $k[G]$-modules $M$ whose underlying $k$-vector space is finite-dimensional. Then $\Rep_k(G)$ and the category $\gr\Rep_k(G)$ of
  bounded $\Z$-graded $G$-representations in finite-dimensional $k$-vector
  spaces with the usual notion of short exact sequences are both exact
  categories. Note that both categories are split exact because $k[G]$ is
  semi-simple by Maschke's Theorem.
\end{example}

\begin{definition}[{\cite[Definition 2.1]{schlichting_hermitian_2010}}]
  An \emph{exact category with duality} is a triple $(\mathcal E, \vee, \eta)$,
  where $\mathcal E$ is an exact category, $\vee \colon \mathcal{E}^\op\to
  \mathcal E$ is an exact functor and $\eta\colon \id_{\mathcal E} \to
  \vee^2$ is a natural isomorphism, which consists of morphisms $\eta_A\colon A\to
  A^{\vee\vee}$ for all $A \in \mathcal E$, such that $1_{A^{\vee}} = (\eta_A)^\vee\circ\eta_{A^vee}$ for
  all objects $A$ of $\mathcal E$. We also write $\mathcal E$ in place of
  $(\mathcal E, \vee, \eta)$ if the duality functor $\vee$ and the double dual
  identification $\eta$ are clear from the context.
\end{definition}

In the following, we will consider two exact categories with
duality.

\begin{example}
  \label{example:duality-on-reps}
  This is a continuation of Example \ref{ex:rep-split-exact} and
  both of these examples are essentially variations on \cite[Example~2.2]{schlichting_hermitian_2010}.

  \begin{enumerate}
  \item The category $\Rep_k(G)$ can be made into a split exact category with
    duality by considering the duality functor $\vee\colon \Rep_k(G)^\op \to
    \Rep_k(G)$ defined by $A \mapsto A^\vee := \Hom_{k}(A, k)$. We equip $A^\vee$ with the
    following $G$-action: for $g\in G$ and $f \in A^\vee$, we define $gf$ by
    $(gf)(a) := f(g^{-1}a)$ for $a \in A$. That is, we have the involution
    \[
      \iota\colon k[G]\to k[G]; \sum_{i}\lambda_ig_i \mapsto \sum_i
      \lambda_ig_i^{-1},
    \]
    and we make $ A^\vee$ a left $k[G]$-module by
    $(\alpha\cdot f)(a):=f(\iota(\alpha)\cdot a)$ for $\alpha\in k[G]$, $a\in
    A$ and $f:A\to k$ a $k$-linear map. The evaluation map
    $\operatorname{ev}\colon A \to A^{\vee\vee}$ defines a natural isomorphism
    such that $(\Rep_k(G), \vee, \ev)$ becomes an exact category with
    duality.

    Note that $\vee$ is a $\otimes$-functor and that $\eta$ respects the
    tensor structure; that is, we have $\eta_{A\otimes B} = \eta_A\otimes
    \eta_B$ after identifying $(A\otimes B)^{\vee\vee}$ with
    $A^{\vee\vee}\otimes B^{\vee\vee}$. This is a special case of
    \cite[Example 2.2]{schlichting_hermitian_2010} by considering $k[G]$ with
    the involution $\iota\colon k[G]\to k[G]$ defined above.
  \item We can equip the category $\gr\Rep_k(G)$ with the following
    structure of a split exact category with duality. We let $\vee\colon
    \gr\Rep_k(G)^\op \to \gr\Rep_k(G)$ be the functor sending a graded
    $G$-representation $A^\bullet$ to the graded $G$-representation
    \[
      (A^\bullet)^\vee = \Hom_{\Rep_k(G)}(A^{-\bullet},k).
    \]
    Here, the degree $i$
    part of $(A^\bullet)^\vee$ is given by $\Hom_{\Rep_k(G)}(A^{-i},k)$. We
    define the map $\eta\colon A^\bullet \to (A^\bullet)^{\vee\vee}$ to be the map
    given in degree $i$ by
    \[
      (-1)^i\operatorname{ev}_{A^i}\colon A^i \to
      \Hom_{\Rep_k(G)}(\Hom_{\Rep_k(G)}(A^i,k),k),
    \]
    that is, $a\in A^i$ is sent
    to the homomorphism sending $f \in \Hom_{\Rep_k(G)}(A^i,k)$ to
    $(-1)^if(a)$.

    Note that $\vee$ is again a $\otimes$-functor and that $\eta$ respects the
    tensor structure.\qedhere
  \end{enumerate}
\end{example}

\begin{definition}[{\cite[Definition 2.4]{schlichting_hermitian_2010}}]
  A \emph{symmetric form} on an object $X$ in an exact category with duality
  $(\mathcal E, \vee, \eta)$ is a pair $(X,\varphi)$ where $\varphi\colon X\to
  X^\vee$ is a morphism in $\mathcal E$ satisfying $\varphi^\vee\eta_X =
  \varphi$. If $\varphi$ is an isomorphism, the form is called
  \emph{non-singular} or \emph{non-degenerate}; otherwise \emph{singular} or
  \emph{degenerate}.

  Let $(X,\varphi)$ be a symmetric form in $\mathcal E$ and let $f\colon Y \to X$
  be a morphism. Then $(Y, \varphi|_Y)$ with $\varphi|_Y = f^\vee
  \varphi f$ is a symmetric form on $Y$ called the
  \emph{restriction} via $f$ of $\varphi$. A \emph{map of symmetric forms} $f\colon (Y,\psi) \to
  (X,\varphi)$ is a morphism $f\colon Y \to X$ in $\mathcal E$ such that $\psi
  = \varphi|_Y$. The map $f$ is an \emph{isometry} if $f$ is an
  isomorphism. Composition of morphisms in $\mathcal E$ defines the
  composition of maps between symmetric forms.
\end{definition}

\begin{definition}[{\cite[Definition 2.4]{schlichting_hermitian_2010}}]
  The \emph{orthogonal sum} $(X, \varphi)\perp (Y,\psi)$ of two symmetric
  forms $(X,\varphi)$ and $(Y, \psi)$ is the symmetric form $(X\oplus Y,
  \varphi\oplus \psi)$.
\end{definition}

\begin{definition}
  If $\mathcal E$ is symmetric monoidal, $\vee$ is a
  $\otimes$-functor, and $\eta$ respects the tensor structure, then the \emph{product} $(X, \varphi)\otimes (Y,\psi)$ of two
  symmetric forms is $(X\otimes Y, \varphi\otimes
  \psi)$.
\end{definition}

If $(X,\varphi)$ and $(Y, \psi)$ are non-degenerate symmetric forms, then so
are their orthogonal sum and product. If $\mathcal E$ is essentially small,
which is the case for $\Rep_k(G)$ and $\gr\Rep_k(G)$, then the set of isometry
classes of non-degenerate symmetric forms is a monoid under orthogonal sum and a
semi-ring under orthogonal sum and product.

\begin{definition}[{\cite[Definition 2.5]{schlichting_hermitian_2010}}]
  Let $(X,\varphi)$ be a non-degenerate symmetric form in $\mathcal E$. A
  \emph{totally isotropic subspace} of $X$ is an inflation $i \colon L
  \hookrightarrow X$ with $\varphi|_L = i^\vee\varphi i = 0$ such that the
  induced map $L \to L^\perp := \ker(i^\vee\varphi)$ is also an inflation.

  A \emph{Lagrangian} of $X$ is a morphism $i\colon L \to X$
  such that
  \[
    \begin{tikzcd}
      0\ar[r] & L \ar[r,"i"] & X \ar[r, "i^\vee\varphi"] & L^\vee\ar[r] & 0
    \end{tikzcd}
  \]
  is a short exact sequence. That is, $L$ is a totally isotropic subspace with
  $L = L^\perp$.

  A non-degenerate symmetric form $(X,\varphi)$ is called \emph{metabolic} if
  $(X,\varphi)$ has a Lagrangian. If $L$ is any object of $\mathcal E$, then
  \[
    H(L) = \left(L\oplus L^\vee,
    \begin{pmatrix}
      0 & 1\\
      \eta & 0
    \end{pmatrix}
    \right)
  \]
  is a non-degenerate symmetric form called the \emph{hyperbolic form} of
  $L$. The form $H(L)$ is metabolic with Lagrangian the inclusion $L \to L\oplus L^\vee$
  into the first summand. Note that the symmetric forms $H(L\oplus M)$ and
  $H(L)\perp H(M)$ are isometric.
\end{definition}

\begin{example}
  \label{ex:symmetric-form-char}
  This is a continuation of Example \ref{example:duality-on-reps}.

  \begin{enumerate}
  \item Let $(V,\varphi)$ be a symmetric form in $\Rep_k(G)$. Thus, $\varphi$
    is a $G$-equivariant linear map $\varphi \colon V \to V^\vee$. Such a
    linear map is the same as a bilinear form
    \[
      \beta \colon V \otimes V
      \to k; \beta(v,w) = \varphi(v)(w)
    \]
    and the $G$-equivariance condition
    translates to the condition that $\beta(gv,w) = \beta(v, g^{-1}w)$ for all
    $v,w \in V$ and $g \in G$. This is equivalent to the condition that
    $\beta(gv,gw) = \beta(v,w)$ for all $v,w \in V$ and $g \in G$. The requirement $\varphi^\vee\eta_X =  \varphi$ translates
    to the condition that $\beta(v,w) = \beta(w,v)$. Thus, symmetric forms on
    $\Rep_k(G)$ are the same as symmetric bilinear forms $\beta\colon V\otimes
    V \to k$ satisfying the additional condition that $\beta(gv, gw) =
    \beta(v,w)$ for all $v,w \in V$ and $g \in G$. We call such a bilinear
    form a \emph{$G$-equivariant symmetric bilinear form}. Such a form is
    non-degenerate if the underlying symmetric bilinear form is non-degenerate.
  \item Let $(V^\bullet, \varphi)$ be a symmetric form in $\gr\Rep_k(G)$. Hence,
    the degree $i$-part $\varphi^i$ of $\varphi$ is a $G$-equivariant linear map
    $\varphi^i\colon V^i \to (V^{-i})^\vee$. Such a linear map is the same
    as a bilinear form $\beta^i\colon V^i\otimes V^{-i} \to k$ and the
    $G$-equivariance condition again translates to the condition $\beta(gv,w)
    = \beta(v,g^{-1}w)$ for all $v \in V^i, w\in V^{-i}$ and $g \in G$, which
    is again equivalent to the condition $\beta(gv,gw) = \beta(v,w)$ for all $v
    \in V^i, w\in V^{-i}$ and $g \in G$.
    The requirement $\varphi^\vee\eta_{V^\bullet} =  \varphi$ translates to the
    following condition: the bilinear forms
    \[
      \beta^0\colon V^0\otimes V^0 \to
      k\quad \text{and}\quad \beta^{2i}\perp\beta^{-2i}\colon (V^{2i}\oplus V^{-2i})\otimes (V^{2i}\oplus V^{-2i}) \to
      k
    \]
    are symmetric for all $i \ne 0$ and the bilinear forms
    \[
      \beta^{2i+1}\perp
      \beta^{-2i-1}\colon (V^{2i+1}\oplus V^{-2i-1})\otimes (V^{2i+1}\oplus
      V^{-2i-1}) \to k
    \]
    are alternating for all $i\in \Z$. The form $(V^\bullet, \varphi)$ is non-degenerate,
    if $\varphi^i\colon V^i \to (V^{-i})^\vee$ is an isomorphism for all
    $i$.

    We call a collection $\beta^\bullet$ of $k$-bilinear forms $\beta^i\colon V^i\otimes
    V^{-i}\to k$ for all $i \in \Z$ satisfying the conditions outlined above a
    \emph{$\Z$-graded $G$-equivariant symmetric bilinear} form. We call $\beta^\bullet$
    non-degenerate if the corresponding morphism $\varphi^\bullet\colon
    V^\bullet \to (V^\bullet)^\vee$ is an
    isomorphism.\label{ex:symmetric-form-char:Z-grading}\qedhere
  \end{enumerate}
\end{example}

\begin{rem}
  The definition of a $G$-equivariant symmetric bilinear form in Example~\ref{ex:symmetric-form-char} is equivalent to the notion of a
  $G$-equivariant quadratic form in \cite[Definition~3.1]{Pajwani-PalYZ} under
  the correspondence between symmetric bilinear forms and quadratic forms.
\end{rem}

\begin{notation}
  Let $\beta^\bullet$ be a $\Z$-graded non-degenerate $G$-equivariant
  symmetric bilinear form on $V^\bullet$. We write $V^\even := \bigoplus_{i\in \Z}V^{2i}$ and $V^\odd :=
  \bigoplus_{i\in \Z}V^{2i+1}$. We denote the non-degenerate symmetric
  bilinear form $\perp_{i\in \Z}\beta^{2i}$ on $V^\even$ by $\beta^\even$ and the non-degenerate
  alternating bilinear form $\perp_{i\in \Z}\beta^{2i+1}$ on $V^\odd$ by
  $\beta^\odd$. Note that $V^\odd$ is even-dimensional because $\beta^\odd$ is
  non-degenerate and alternating.
\end{notation}

\begin{definition}[{\cite[p.~7 and Corollary~2.10]{schlichting_hermitian_2010}}]
  \label{def:GW-category}
  The \emph{Grothendieck-Witt group} $\GW(\mathcal E)$ of a split exact
  category with duality is the Grothendieck group of the Abelian monoid of
  isometry classes $[X,\varphi]$ of non-degenerate symmetric forms
  $(X,\varphi)$ in $\mathcal E$ with respect to orthogonal sums.
\end{definition}

\begin{rem}
  If $\mathcal E$ is symmetric monoidal, $\vee$ is a
  $\otimes$-functor, and $\eta$ respects the tensor structure, then the product of symmetric forms
  induces a ring structure on $\GW(\mathcal E)$. In this case, we will also
  call $\GW(\mathcal E)$ the \emph{Grothendieck-Witt ring} of $\mathcal E$.
\end{rem}

\begin{notation}
  For the exact categories with duality $\Rep_k(G)$ and $\gr\Rep_k(G)$
  we fix the following notation: write $\GW^G(k)$ for $\GW(\Rep_k(G))$
  and $\gr\GW^G(k)$ for $\GW(\gr\Rep_k(G))$. We call $\GW^G(k)$ the
  \emph{$G$-equivariant Grothendieck-Witt ring} of $k$.
\end{notation}

\begin{rem}
  The definition of $\GW^G(k)$ using Definition \ref{def:GW-category} agrees with the definition of the
  $G$-equivariant Grothendieck-Witt ring in \cite[Definition
  3.1]{Pajwani-PalYZ}. Furthermore, if $G$ is the trivial group, then
  $\GW^G(k)$ agrees with the Grothendieck-Witt ring $\GW(k)$ as introduced in Section~\ref{sec:euler-char}.
\end{rem}

\begin{rem}\label{rem:z-graded-nonzero-hyperbolic}
  Let $\beta^\bullet$ be a $\Z$-graded non-degenerate $G$-equivariant
  symmetric bilinear form on $V^\bullet$. Since $\beta^\bullet$ is non-degenerate, it induces
  an isomorphism $V^{2i}\cong V^{-2i}$. Thus, $\bigoplus_{i\ge 1}V^{2i}$ is a
  Lagrangian for the $\Z$-graded non-degenerate $G$-equivariant symmetric
  bilinear form constructed from $\beta$ by removing all odd degree parts and
  the degree zero part of $\beta$. Thus,  \cite[Lemma
  2.8]{schlichting_hermitian_2010} yields $\perp_{i\in
    \Z\setminus\{0\}} \beta^{2i} \cong H(\bigoplus_{i\ge 1} V^{2i})$,
  which is hyperbolic. In
  particular, we have $\beta^\even = \beta^0 \perp H(\bigoplus_{i\ge 1}
  V^{2i})$ and thus, $\beta^\even = \beta^0 + (\sum_{i\ge 1}\dim_kV^{2i})\cdot H \in \GW(k)$.
\end{rem}

\begin{rem}
  \label{rem:even-odd-mult}
  Let $\beta^\bullet$ be a $\Z$-graded $G$-equivariant non-degenerate
  symmetric bilinear form on $V^\bullet$ and let $\gamma^\bullet$ be a
  $\Z$-graded $G$-equivariant non-degenerate symmetric bilinear form on $W^\bullet$. Then we
  can describe $(\beta\cdot \gamma)^\even$ and $(\beta\cdot \gamma)^\odd$
  explicitly: the grading on $V\otimes W$ yields
  \[
    (V\otimes W)^\even = V^\even \otimes W^\even \oplus V^\odd \otimes W^\odd
  \]
  and
  \[
    (V\otimes W)^\odd = V^\even \otimes W^\odd \oplus V^\odd \otimes W^\even.
  \]
  This decomposition extends to $\beta\cdot \gamma$.
  That is, we have
  \[
    (\beta\cdot \gamma)^\even = \beta^\even \otimes
    \gamma^\even \perp \beta^\odd \otimes \gamma^\odd
  \]
  as $G$-equivariant symmetric bilinear forms and
  \[
    (\beta\cdot
    \gamma)^\odd = \beta^\even \otimes \gamma^\odd \perp \beta^\odd\otimes
    \gamma^\even
  \]
  as $G$-equivariant bilinear forms.
  We will often write these decompositions as $(\beta\cdot \gamma)^\even = \beta^\even \cdot
    \gamma^\even + \beta^\odd \cdot \gamma^\odd$ and $(\beta\cdot
  \gamma)^\odd = \beta^\even \cdot \gamma^\odd + \beta^\odd\cdot
  \gamma^\even$.
\end{rem}

We are now turning towards defining a $G$-equivariant quadratic Euler
characteristic.

\begin{proposition}
  \label{prop:canonical-map-Z}
  There is a canonical ring homomorphism $\can_\Z\colon \gr\GW^G(k) \to \GW^G(k)$ given by sending a
  $\Z$-graded non-degenerate $G$-equivariant symmetric bilinear form
  $\beta^\bullet$ on $V^\bullet$ in $\gr\Rep_k(G)$ to
  \[
    \beta^\even - [\perp_{i \ge 0} H(V^{2i+1})] \in \GW^G(k).
  \]
  In particular, if $G$ is the trivial group, we have
  \[
    \can_\Z(\beta) = \beta^\even - \frac
    12\dim_kV^\odd \cdot H \in \GW(k).
  \]

  \begin{proof}
    First note that this defines a well-defined element in $\GW^G(k)$: since isometry
    of $\Z$-graded $G$-equivariant symmetric bilinear forms respects degrees,
    the definition of $\can_\Z$ does not depend on the chosen representative
    of an isometry class of a symmetric bilinear form. Also note that $\bigoplus_{i\ge 0}V^{2i+1}$ defines a Lagrangian of
    the odd-degree part of $\beta^\bullet$. Furthermore, $\can_\Z$ is canonical because the class of
    $[\perp_{i \ge 0} H(V^{2i+1})] \in \GW^G(k)$ does not depend on the chosen
    Lagrangian by \cite[Lemma~2.8~(c)]{schlichting_hermitian_2010}.

    The definition on $\Z$-graded non-degenerate $G$-equivariant
    symmetric bilinear form extends to a group homomorphism $\can_\Z\colon
    \gr\GW^G(k) \to \GW^G(k)$. Indeed, suppose we have $\Z$-graded non-degenerate $G$-equivariant
    symmetric bilinear forms $\beta^\bullet$
    on $V^\bullet$ and $\gamma^\bullet$ on $W^\bullet$. Now,
    $(\beta\perp\gamma)^i$ is the $k$-bilinear map $(V^i\oplus W^i)\otimes
    (V^{-i}\oplus W^{-i}) \to k$ given by sending $(v_i+w_i) \otimes (v_{-i}+
    w_{-i})$ to $\beta^i(v_i,v_{-i}) + \gamma(w_i,w_{-i})$ for $v_i \in V^i,
    w_i\in W^i, v_{-i}\in V^{-i}$, and $w_{-i} \in W^{-i}$. Thus, we have
    \begin{align*}
      \can_\Z(\beta^\bullet + \gamma^\bullet)
      &= [\perp_{i\in \Z}\underset{=\beta^{2i}\perp
        +\gamma^{2i}}{\underbrace{(\beta^\bullet
          +\gamma^\bullet)^{2i}}}] - [\perp_{i \ge 0} \underset{=H(V^{2i+1})\perp
        H(W^{2i+1})}{\underbrace{H(V^{2i+1}\oplus W^{2i+1})}}]\\
      &= \beta^\even  +
      \gamma^\even -
        [\perp_{i \ge 0} H(V^{2i+1})] - [\perp_{i \ge 0} H(W^{2i+1})]\\
      &= \can_\Z(\beta^\bullet) + \can_\Z(\gamma^\bullet).
    \end{align*}

    Similarly, One can check similarly that $\can_\Z(\beta^\bullet\cdot \gamma^\bullet)
    = \can_\Z(\beta^\bullet)\cdot \can_\Z(\gamma^\bullet)$ for $\Z$-graded non-degenerate $G$-equivariant
    symmetric bilinear forms $\beta^\bullet$
    and $\gamma^\bullet$, which implies that $\can_\Z$ is a ring homomorphism.
  \end{proof}
\end{proposition}

\begin{example}
  \label{ex:graded-euler-char}
  Let $X$ be an equidimensional, smooth, projective scheme over $k$ of
  dimension $n$ with a $G$-action. Hence, the Hodge
  cohomology groups $H^q(X,\Omega^p_{X/k})$ carry a natural $G$-action. The
  composition of cup-product and trace
  \[
    \begin{tikzcd}
      H^q(X,\Omega^p_{X/k})\times H^{n-q}(X,\Omega^{n-p}_{X/k}) \ar[r, "\cup"]
      & H^n(X,\Omega^n_{X/k}) \ar[r, "\Tr"] & k
    \end{tikzcd}
  \]
  defines a $\Z$-graded $G$-equivariant non-degenerate symmetric bilinear form
  \[
    \left(\bigoplus_{i,j=0}^{\dim_kX}
      H^i(X,\Omega^j_{X/k})[j-i],\Tr\right)
  \]
  in $\gr\GW^G(k)$. If we take $G$ to be the trivial group, then the image of
  this symmetric form in $\GW^G(k) = \GW(k)$ computes the $\A^1$-Euler
  characteristic of $X$ over $k$.
\end{example}

\begin{definition}
  \label{def:equiv-quadratic-euler-char}
  Let $X$ be an equidimensional smooth, projective scheme over $k$ with a $G$-action. We define the
  \emph{$G$-equivariant quadratic Euler characteristic} $\chi_G(X/k)$ of $X$
  over $k$ as the image of the non-degenerate symmetric form $\left(\bigoplus_{i,j=0}^{\dim_kX}
    H^i(X,\Omega^j_{X/k})[j-i],\Tr\right)$ in $\GW^G(k)$.

  If $X$ is an arbitrary smooth, projective scheme over $k$ with a $G$-action,
  we write $X = \coprod_{i=1}^N X_i$ as a disjoint union of equidimensional open and
  closed subschemes of $X$ closed under the $G$-action. Then we define the
  \emph{$G$-equivariant quadratic Euler characteristic} $\chi_G(X/k)$ of $X$
  over $k$ as $\chi_G(X/k) = \sum_{i=1}^N \chi_G(X_i/k)$.
\end{definition}

\begin{proposition}
  \label{prop:equiv-euler-char-equal-pp}
  Let $k$ be a field of characteristic zero, and let $X$ be a smooth,
  projective scheme over $k$ with a $G$-action. Then $\chi_G(X/k)$ is equal to the
  $G$-equivariant de Rham Euler characteristic of $X$ as
  defined in \cite[Definition 3.5]{Pajwani-PalYZ}.

  \begin{proof}
    Both Euler characteristics are additive in disjoint unions, so we may
    assume that $X$ is equidimensional of dimension $n$.
    We start by recalling the definition of the $G$-equivariant de Rham Euler
    characteristic of $X$ from \cite[Definition 3.5]{Pajwani-PalYZ}. Denote the $i$-th de
    Rham cohomology of $X$ by $H^i_{\dR}(X)$. Define
    \[
      \beta^i\colon (H^i_\dR(X)\oplus H_\dR^{2n-i})\times (H^i_\dR(X)\oplus
      H^{2n-i}_\dR) \xrightarrow{\cup} H^{2n}_\dR(X) \xrightarrow{\Tr} k
    \]
    for $i \in \{0, \dots, n-1\}$ even and define
    \[
      \beta^n\colon H^n_\dR(X) \times H^n_\dR(X) \xrightarrow{\cup}
      H^{2n}_\dR(X) \xrightarrow{\Tr} k.
    \]
    If $n$ is odd, there exists a Lagrangian $L\subset H^n_\dR(X)$ for $\beta$
    closed under the $G$-action by \cite[Lemma 3.4]{Pajwani-PalYZ}.
    Now, the $G$-equivariant de Rham Euler characteristic of $X$ is defined as
    \[
      \chi_G^\dR(X) =
      \begin{cases}
        \sum_{i=0}^{n/2}\beta^{2i} - \sum_{i=1}^{n/2}H(H^{2i-1}_\dR(X))
        & \text{if } n \text{ even}\\
        \sum_{i=0}^{(n-1)/2}\beta^{2i} - \sum_{i=1}^{(n-1)/2}H(H^{2i-1}_\dR(X)) - H(L)
        & \text{if } n \text{ odd}.
      \end{cases}
    \]
    in $\GW^G(k)$. Furthermore, the Hodge to de Rham spectral sequence
    induces an isomorphism $H^{2d}_\dR(X) \cong
    H^n(X,\Omega^n_{X/k})$ by \cite[\href{https://stacks.math.columbia.edu/tag/0FW6}{Tag~0FW6}]{stacks-project}. Under this isomorphism, we obtain an
    identification of the trace map in de Rham cohomology and the trace map in
    Hodge cohomology by \cite[Proof of
    \href{https://stacks.math.columbia.edu/tag/0FW7}{Tag
      0FW7}]{stacks-project}.

    Note that every short exact sequence of $k[G]$-modules splits.
    Since $X$ is smooth and projective, the Hodge to de Rham spectral sequence
    degenerates on the first page and is compatible with the
    $G$-action by \cite[Théorème~5.5~(ii)]{Deligne1968HodgeSS}. Thus, the filtration
    of $H^\ast_\dR(X)$ induced by the Hodge to de Rham spectral sequence
    guarantees the existence of a $G$-equivariant isomorphism
    \[
      H^i_\dR(X) \cong \bigoplus_{p+q = i}H^q(X,\Omega^p_{X/k}).
    \]
    Hence, we get
    \begin{equation}
      \label{eq:equiv-euler-char-equal-pp:hyperbolic}
      H(H^i_\dR(X)) = \sum_{p+q = i} H(H^q(X,\Omega^p_{X/k})) \in \GW^G(k).
    \end{equation}
    for $i= 1, \dots, n-1$. Furthermore, note that we have
    \begin{equation}
      \label{eq:equiv-euler-char-equal-pp:non-middle}
      \beta^{2i} = H(H^i_\dR(X)) \in \GW^G(k)
    \end{equation}
    for $0 \le i < n/2$ by Poincaré Duality, see \cite[\href{https://stacks.math.columbia.edu/tag/0FW7}{Tag 0FW7}]{stacks-project}.

    Now, let $n = 2m$ be even. Hence, we have a filtration $F^{n+1}\subset \dots
    \subset F^{0}= H^n_\dR(X)$ such that $F^{i}/F^{i+1} =
    H^i(X,\Omega^{n-i}_{X/k})$ as $k[G]$-modules. Since the filtration is multiplicative by
    \cite[\href{https://stacks.math.columbia.edu/tag/0FM7}{Tag 0FM7}]{stacks-project},
    we see that $F^{m+1}$ is a totally isotropic subspace of
    $H^n_\dR(X)$ with respect to $\beta^n$ and that $F^m\subset
    (F^{m+1})^\perp$. Furthermore, a dimension count using Serre duality
    yields $F^m = (F^{m+1})^\perp$. Thus, we get by \cite[Lemma 2.8
    (c)]{schlichting_hermitian_2010} that
    \[
      \beta^n = \gamma + H(F^{m+1}),
    \]
    where $\gamma$ is the symmetric bilinear form
    \[
      H^m(X, \Omega^m_{X/k}) \times H^m(X, \Omega^m_{X/k}) \xrightarrow{\cup}
      H^{n}(X, \Omega^n) \xrightarrow{\Tr} k.
    \]
    The argument to obtain
    \eqref{eq:equiv-euler-char-equal-pp:hyperbolic} also shows
    \[
      H(F^{m+1}) = \sum_{i=m+1}^{n} H(H^i(X,\Omega^{n-i}_{X/k})).
    \]
    By combining this with \eqref{eq:equiv-euler-char-equal-pp:hyperbolic}
    and \eqref{eq:equiv-euler-char-equal-pp:non-middle}, we get
    \begin{align*}
      \chi_G^\dR(X)
      &= \sum_{i=0}^{n/2}\beta^{2i} - \sum_{i=1}^{n/2}H(H^{2i-1}_\dR(X))\\
      &= \gamma + \sum_{i=m+1}^{n} H(H^i(X,\Omega^{n-i}_{X/k})) +
        \sum_{i=0}^{m-1} \sum_{p+q = 2i} H(H^q(X,\Omega^p_{X/k}))\\
      &\quad -
        \sum_{i=1}^{n/2}\sum_{p+q = 2i-1} H(H^q(X,\Omega^p_{X/k}))\\
      &= \chi_G(X/k).  
    \end{align*}
    The case $n$ odd works similarly. The only difference is that $\gamma$
    does not show up and that we identify $L$ with $F^{(n+1)/2}$.
  \end{proof}
\end{proposition}

\begin{rem}[{\cite[15]{Pajwani-PalYZ}}]
  \label{rem:gw-forget-action}
  There is a canonical ring homomorphism $\GW^G(k) \to \GW(k)$ given by
  forgetting the group action. This map sends $\chi_G(X/k)$ to the
  $\A^1$-Euler characteristic $\chi_c(X/k)$ as defined in Section
  \ref{sec:euler-char}.
\end{rem}

\begin{rem}
  There also is a canonical ring homomorphism $\varphi\colon \GW(k) \to \GW^G(k)$ induced by
  equipping a quadratic form with the trivial $G$-action. The ring
  homomorphism $\varphi$ is a right inverse to the ring homomorphism $\GW^G(k)
  \to \GW(k)$ forgetting the group action from Remark \ref{rem:gw-forget-action}. Moreover, $\varphi$
  induces a $\GW(k)$-module structure on $\GW^G(k)$.

  Similarly, we have a canonical ring homomorphism $\gr\GW(k) \to
  \gr\GW^G(k)$ equipping $\beta^\bullet \in \gr\GW(k)$ with the trivial group
  action. This equips $\gr\GW^G(k)$ with the structure of a $\gr\GW(k)$-module.
\end{rem}

\begin{proposition}
  Let $H \subset G$ be a normal subgroup.
  The restriction map sending a $G$-equivariant symmetric bilinear form $\beta$ on $V$
  to the restriction $\beta^H$ on the $H$-invariants $V^H$ induces a
  well-defined $\GW(k)$-linear map $(-)^H\colon \GW^G(k) \to \GW^{G/H}(k)$. We write
  this map as $\beta \mapsto \beta^H$.

  We also have a similar $\gr\GW(k)$-linear map $(-)^H\colon \gr\GW^G(k) \to
  \gr\GW^{G/H}(k)$.

  \begin{proof}
    Let $\beta$ be a non-degenerate $G$-equivariant symmetric bilinear form on
    the $G$-rep\-re\-sen\-ta\-tion $V$ and consider $V$ as an $H$-representation. Write $V$ as a sum of its isotypical
    components,  $V = \bigoplus_{\rho}V^\rho$, where the sum is over the
    irreducible $H$-representations of $V$. Since $\beta$ is
    non-degenerate, we find that the induced linear map $V\to V^\vee; v\mapsto
    \beta(v,-)$ is an $H$-equivariant isomorphism. Since the isomorphism is
    $H$-equivariant, it restricts to an isomorphism $V^\rho \to
    (V^\rho)^\vee$, where $\rho$ denotes the trivial $H$-representation. In particular,
    the restriction of $\beta$ to $V^\rho$ endowed with the induced $G/H$-action, which is $\beta^H$, is again a
    non-degenerate symmetric bilinear form. That is, $\beta^H$ is an element of
    $\GW^{G/H}(k)$.

    Since taking $H$-invariants is an additive functor, we get that
    $(\beta+\gamma)^H = \beta^H + \gamma^H$ for all non-degenerate
    $G$-equivariant symmetric bilinear forms $\beta$ and $\gamma$. Thus, the
    definition of $\beta^H$ induces a well-defined group homomorphism
    $(-)^H\colon \GW^G(k) \to \GW^{G/H}(k)$.

    It remains to show that $(-)^H$ is compatible with the $\GW(k)$-module structure. For this, note that if
    $V$ is a $k$-vector space equipped with the
    trivial $G$-action and $W$ is a $G$-representation, then we have $(V\otimes W)^H = V\otimes W^H =
    V^H\otimes W^H$. This in
    turn implies $(\beta\cdot \gamma)^H = \beta^H\cdot\gamma^H = \beta
    \cdot \gamma^H$ for all non-degenerate
    symmetric bilinear forms $\beta$ and non-degenerate $G$-equivariant
    symmetric bilinear forms $\gamma$. This proves that
    $(-)^H$ is $\GW(k)$-linear.

    One analogously obtains the statement for $(-)^H\colon \gr\GW^G(k) \to \gr\GW^{G/H}(k)$.
  \end{proof}
\end{proposition}

\begin{rem}
  The map $(-)^H \colon \GW^G(k)\to \GW^{G/H}(k)$ is not a ring homomorphism. That
  is $(-)^H$ is not multiplicative: indeed if $G = \Z/2$, consider the
  $G$-equivariant symmetric bilinear form $\beta\colon k\times k \to k;
  (x,y)\mapsto xy$ where we act on $k$ by multiplication with $-1$. Then we have $\beta^G
  =0$ and $(\beta\cdot \beta)^G = \langle 1\rangle \ne 0 = \beta^G \cdot \beta^G$.
\end{rem}

\begin{rem}
  If $H \subset G$ is a subgroup, we can canonically consider an element $\beta
  \in \GW^G(k)$ as an element of $\GW^H(k)$ by restricting the $G$-action to
  an $H$ action. In particular, we can always compute
  $\beta^H \in \GW(k)$ for all $\beta\in \GW^G(k)$ by considering $\beta$ as
  an element of $\GW^H(k)$. In the following, we implicitly consider
  $\beta\in \GW^G(k)$ as an element of $\GW^H(k)$ when computing the
  $H$-invariants of $\beta$ in $\GW(k)$.

  Since this discussion extends to $\gr\GW^G(k)$ and $\gr\GW^H(k)$ for $H
  \subset G$, we can also take the $H$-invariants of $\beta^\bullet \in
  \gr\GW^G(k)$ in $\gr\GW(k)$.
\end{rem}

\begin{proposition}
  For $H\subset G$ a normal subgroup, the following diagram commutes
  \[
    \begin{tikzcd}
      \gr\GW^G(k)\ar[d, "(-)^H"]\ar[r, "\can_{\Z}"] &
      \GW^G(k)\ar[d, "(-)^H"]\\
      \gr\GW^{G/H}(k)\ar[r, "\can_\Z"] & \GW^{G/H}(k).
    \end{tikzcd}
  \]

  \begin{proof}
    This can be checked on $\Z$-graded non-degenerate $G$-equivariant
    symmetric bilinear forms, where it is a straightforward
    computation.
  \end{proof}
\end{proposition}

\subsection{Properties}
The $G$-equivariant quadratic Euler characteristic shares many formal
properties with the $\A^1$-Euler characteristic introduced in Section
\ref{sec:euler-char}.

\begin{rem}
  Pajwani and Pál \cite[Theorem 3.6]{Pajwani-PalYZ} show that the
  $G$-equivariant de Rham Euler characteristic
  can be extended to a ring homomorphism $K^G_0(\Var_k) \to \GW^G(k)$ in
  characteristic zero. Here $K^G_0(\Var_k)$ denotes the Grothendieck ring of
  varieties with a $G$-action. Many of the properties that we show about the
  $G$-equivariant quadratic Euler characteristic here can be deduced directly
  from this description in characteristic zero because the $G$-equivariant de Rham Euler
  characteristic agrees with the $G$-equivariant quadratic Euler
  characteristic by Proposition~\ref{prop:equiv-euler-char-equal-pp}.
\end{rem}

\begin{rem}
  \label{rem:equiv-gauss-bonnet-idea}
  Hoyois \cite{HoyoisSOEMHT} constructs a $G$-equivariant stable homotopy category
  $\SH^G(X)$ for $X$ a scheme with a $G$-action. By \cite[Corollary~6.13]{HoyoisSOEMHT}, $p_\#(1_X) \in
  \SH^G(k)$ is again strongly dualisable for a smooth, projective $k$-scheme $p\colon X \to \Spec k$
  equipped with a $G$-action, and it therefore
  possesses a categorical Euler characteristic. One can now map this to the
  $(0,0)$-th homotopy group of the $G$-equivariant Hermitian K-theory
  spectrum $\KQ^{0,0}(\Spec k,G)$
  using its unit map $u \colon \End_{\SH^G(k)}(1_k) \to \KQ^{0,0}(\Spec k, G)$. By
  this, one would obtain another invariant $u(\chi^{\text{cat}}(p_\#1_X)) \in \KQ^{0,0}(\Spec
  k,G)$ qualifying for the name of $G$-equivariant quadratic Euler
  characteristic.

  One can now try to identify $\KQ^{0,0}(\Spec k, G)$ with $\GW^G(k)$ and then try
  to reprove a Motivic Gauß-Bonnet Theorem in this setting, that is, showing
  that this construction agrees with the construction of $G$-equivariant
  quadratic Euler characteristic from Definition
  \ref{def:equiv-quadratic-euler-char}. With such a theorem, one would
  probably be able
  to prove the results in this section using motivic homotopy
  theory in a more elegant way.
\end{rem}

Since there currently is no Motivic Gauß-Bonnet Theorem in the $G$-equivariant
setting, as outlined in Remark
\ref{rem:equiv-gauss-bonnet-idea}, we are going to prove the properties in
this section using Hodge cohomology directly. For this, we shall first prove
two statements about Hodge cohomology, Proposition \ref{prop:hodge-pbf} and
Proposition \ref{prop:hodge-blowup}. These statements are well-known to the experts.
We include their proofs here for the convenience of the reader since proofs
over fields different from the field of complex numbers are hard to find.

\begin{proposition}
  \label{prop:hodge-pbf}
  Let $V$ be a vector bundle of rank $r+1$ over a smooth, projective scheme
  $X$ over $k$. Assume that $G$ acts on $X$ and equivariantly on $V$. This equips $H^\ast(X,\Omega^\ast_{X/k})$ with
  a $G$-action, compatible with the natural $k$-algebra structure, and equips
  $H^\ast(\PP(V),\Omega^\ast_{\PP(V)/k})$ with a $G$-action,
  compatible with its natural $H^\ast(X,\Omega^\ast_{X/k})$-algebra
  structure. Let $\xi
  \in H^1(\PP(V),\Omega^1_{\PP(V)/k})$ be the first Chern class of
  $\OO_V(1)$. Then $\xi$ is fixed under the induced $G$-action on
  $H^1(\PP(V),\Omega^1_{\PP(V)/k})$ and we have an isomorphism of
  $H^\ast(X,\Omega^\ast_{X/k})$-algebras, compatible with the respective
  $G$-actions, given by
  \[
    \frac{H^\ast(X,\Omega^{\ast}_X)[\xi]}{(\xi^{r+1}+\sum_{i=1}^{r}(-1)^ic_i(V)\xi^{r+1-i})}
    \to H^\ast(\PP(V),\Omega^\ast_{\PP(V)/k})
  \]
  sending $\alpha\cdot \xi^i$ with $\alpha \in H^{q-i}(X,\Omega^{p-i})$ to
  $p^\ast(\alpha)\cup \xi^i$. Here $p \colon \PP(V)\to X$ is the projection,
  $c_i(V) \in H^i(X,\Omega^i_{X/k})$ for $i = 0, \dots, r$ is the $i$-th Chern
  class of $V$ in Hodge cohomology, and we consider $\xi \in H^\ast(X,\Omega^{\ast}_X)[\xi]$ to be
  of bi-degree $(1,1)$.

  \begin{proof}
    The $G$-action on $V$ induces a $G$-action $G\times_k \PP(V)\to \PP(V)$ on
    $\PP(V)$, making the projection $p\colon \PP(V)\to X$ a $G$-equivariant morphism,
    and provides a
    $G$-linearization of the tautological invertible quotient sheaf
    $\OO_V(1)$. In particular, for $g\in G$, the action
    map $\psi_g\colon \PP(V)\to \PP(V)$ satisfies $\psi_g^\ast(\OO_V(1)) =
    \OO_V(1)$, hence, $\psi_g^\ast(\xi) = \xi$, that is, $\xi$ is
    $G$-invariant.

    To prove the desired isomorphism, we follow Rao-Yang-Yang-Yu \cite[Step 3 in Section 4,
    p. 3017]{rao_hodge_2023}. Let $p\colon \PP(V) \to X$ denote the
    projection. By \cite[Exposé XI, Théorème 1.1]{SGA7-2}, the maps
    \[
      p^\ast(-)\cup\xi^i\colon H^{q-i}(X,\Omega^{p-i}_{X/k})\to H^q(\PP(V), \Omega_{\PP(V)/k}^p), \ \alpha\mapsto p^\ast(\alpha)\cup\xi^i,
    \]
    induce a $G$-equivariant isomorphism of bi-graded
    $H^\ast(X,\Omega^\ast_{X/k})$-modules
    \[
      \bigoplus_{i=0}^rp^\ast(-)\cup\xi^i\colon \bigoplus_{i=0}^r
      H^\ast(X,\Omega^\ast_{X/k})\xi^i \to H^\ast(\PP(V), \Omega_{\PP(V)/k}^\ast).
    \]
    Furthermore, the morphism
    \[
      p^\ast(-)\cup\xi^\ast\colon H^\ast(X,\Omega^\ast_{X/k})[\xi]\to H^\ast(\PP(V), \Omega_{\PP(V)/k}^\ast).
    \]
    sending $\alpha\cdot x^i$ to $p^\ast(\alpha)\xi^i$ is a homomorphism of
    bi-graded rings.
    Here, we give $\xi$ bi-degree $(1,1)$.
    By \cite[Theorem 1 (g) and Theorem 2]{Srinivas}
    and \cite[Théorème 1]{Grothendieck1958Chern}, the kernel of
    $\sum_{i=0}^{r-1}p^\ast(-)\cup\xi^i$ is the ideal
    \[
      \left(\xi^{r+1}+\sum_{i=1}^{r}(-1)^ic_i(V)\xi^{r+1-i}\right).
    \]
    This yields the desired isomorphism.
  \end{proof}
\end{proposition}

\begin{proposition}
  \label{prop:hodge-blowup}
  Let $G$ act on a smooth, projective scheme $X$ and let $Z \subset X$ be a
  smooth, closed subscheme that is stable under the $G$-action. Let
  $\tilde X := \Bl_ZX$ be the blow-up of $X$ along $Z$ with projection
  $\pi\colon \tilde X\to X$ and let
  $i_E\colon E \to \tilde X$ be the inclusion of the exceptional divisor
  $E$. Equip $E$ and $\tilde X$ with the induced $G$-action. Then the sequence
  \begin{equation}
    \label{eq:hodge-blowup:sequence}
    \begin{tikzcd}
      0 \ar[r] & H^q(X,\Omega^p_{X/k}) \ar[r, "\pi^\ast"] & H^q(\tilde X,
      \Omega_{\tilde X/k}^p) \ar[l, "\pi_\ast", bend left] \ar[r, "i_E^\ast"] &
      \frac{H^q(E,\Omega_{E/k}^p)}{\pi^\ast H^q(Z, \Omega^p_{Z/k})}\ar[r] & 0
    \end{tikzcd}
  \end{equation}
  is split-exact and $G$-equivariant and split by $\pi_\ast$. Furthermore, the induced splitting
  \[
    H^q(\tilde X, \Omega_{\tilde X/k}^p) \cong H^q(X, \Omega^p_{X/k})\oplus
    \frac{H^q(E,\Omega_{E/k}^p)}{\pi^\ast H^q(Z, \Omega^p_{Z/k})}
  \]
  is orthogonal with respect to cup products and compatible with the
  $G$-action. Here, $\pi_\ast$ denotes the pushforward
  constructed by Srinivas \cite[Theorem 1]{Srinivas}.

  \begin{proof}
    The sequence \eqref{eq:hodge-blowup:sequence} is $G$-equivariant since
    $\pi$ and $i_E$ are $G$-equivariant.
    Note that $\pi_\ast\pi^\ast$ is the identity by \cite[Theorem 1 (d) and
    (f)]{Srinivas} and that $\pi_\ast$ is $G$-equivariant by \cite[Theorem 1
    (c)]{Srinivas}. In particular, we get that $\pi^\ast$ is a split
    injection. Rao-Yang-Yang-Yu \cite[Theorem 1.2 and its
    proof]{rao_hodge_2023} prove that the sequence \eqref{eq:hodge-blowup:sequence} is exact
    if $k$ is algebraically closed. Since the base-change to an algebraic
    closure of $k$ is faithfully flat, we get that the sequence
    \eqref{eq:hodge-blowup:sequence} is exact over a general $k$.

    Consider $x \in H^q(X,\Omega^p_{X/k})$ and $y \in
    H^{q'}(E,\Omega_{E/k}^{p'})/\pi^\ast H^{q'}(Z, \Omega^{p'}_{Z/k})$. Pick the
    lift $\tilde y \in H^{q'}(\tilde X, \Omega_{\tilde X}^{p'})$  of $y$ with
    $\pi_\ast(\tilde y) =0$. In order to see that the induced splitting is
    orthogonal with respect to cup products, we
    have to show that $\pi^\ast(x) \cdot\tilde y = 0$. By multiplicativity of
    pullbacks, we have $i_E^\ast(\pi^\ast(x) \cdot \tilde y) = i_E^\ast(\pi^\ast(x))
    \cdot i_E^\ast(\tilde y) = 0 \cdot y = 0$ because \eqref{eq:hodge-blowup:sequence} is
    exact. Furthermore, we have by the projection formula \cite[Theorem 1
    (d)]{Srinivas} that
    $\pi_\ast(\pi^\ast(x)\cdot \tilde y) = x \cdot \pi_\ast(\tilde y) = x
    \cdot 0 = 0$. Since $i_E^\ast$ and $\pi_\ast$ induce the splitting, we get
    $\pi^\ast(x)\cdot \tilde y = 0$. Thus, the splitting is orthogonal with respect to cup
    products.
  \end{proof}
\end{proposition}

\begin{proposition}
  \label{prop:eqec-procut-formula}
  Let $X$ and $Y$ be smooth, projective schemes over $k$ equipped with a
  $G$-action. Then we have
  \[
    \chi_G(X\times Y/k) = \chi_G(X/k) \cdot \chi_G(Y/k).
  \]

  \begin{proof}
    Both sides of the above equation are bilinear with respect to disjoint
    unions. Hence, we may assume that $X$ and $Y$ are both equidimensional.
    Note that the isomorphism on Hodge cohomology induced by the Künneth Theorem
    \cite[\href{https://stacks.math.columbia.edu/tag/0BED}{Tag
      0BED}]{stacks-project}
    is compatible with the $G$-action. Let $p_X\colon X\times Y \to X$ and
    $p_Y\colon X\times Y \to Y$ be the projection morphisms. The isomorphism $p_X^\ast
    \omega_{X/k}\otimes p_Y^\ast\omega_Y \cong \omega_{X\times
      Y/k}$ induces an isomorphism
    \begin{equation}
      \label{eq:eqec-product-formula:topcohom}
      H^{\dim X + \dim Y}(X\times Y, \omega_{X\times Y/k}) \cong H^{\dim
        X}(X, \omega_{X/k}) \otimes_k H^{\dim Y}(Y,\omega_{Y/k})
    \end{equation}
    using the Künneth Theorem.
    We now want to compute the trace $\Tr_{X\times Y/k}\colon H^{\dim X+ \dim
      Y}(X\times Y, \omega_{X\times Y/k}) \to k$. Let $x \in H^{\dim
      X}(X, \omega_{X/k})$ and $y\in H^{\dim Y}(Y,\omega_{Y/k})$. Let
    $\pi_X\colon X \to \Spec k$ and $\pi_Y\colon Y\to \Spec k$ be the
    structure maps. We
    have using \cite[Theorem~1]{Srinivas}
    \begin{align*}
      \Tr_{X\times Y/k}(p_X^\ast x\cdot p_Y^\ast y)
      &= \pi_{X\ast}p_{X\ast}(p_X^\ast x\cdot p_Y^\ast y)\\
      &= \pi_{X\ast}(x\cdot p_{X\ast}p_Y^\ast y)\\
      &= \pi_{X\ast}(x\cdot \pi_X^\ast\pi_{Y\ast}y)\\
      &= \pi_{X\ast}(x)\cdot \pi_{Y\ast}(y)\\
      &= \Tr_{X/k}(x) \cdot \Tr_{Y/k}(y).
    \end{align*}
    More specifically, we use \cite[Theorem~1 (b) and (f)]{Srinivas} for the
    first and the last equality, \cite[Theorem~1 (d)]{Srinivas} for the second
    and fourth equality, and \cite[Theorem~1 (c)]{Srinivas} for the third
    equality. Hence, we get $\Tr_{X\times Y/k} = \Tr_{X/k}\cdot
    \Tr_{Y/k}$ under the identification
    \eqref{eq:eqec-product-formula:topcohom}.
    Combining the Künneth isomorphism with this trace identification yields the
    claim.
  \end{proof}
\end{proposition}

\begin{lemma}
  \label{lem:pbf-trace}
  Let $V$ be a rank $r+1$ vector bundle over a connected, smooth, projective $k$-scheme
  $X$. Let $\xi
  \in H^1(\PP(V),\Omega^1_{\PP(V)/k})$ be the first Chern class of
  $\OO_V(1)$ and let $\alpha \in H^\ast(X,\Omega^\ast_{X/k})$. Then
  \[
    \Tr_{\PP(V)/k}(p^\ast(\alpha)\cdot \xi^i) =
    \begin{cases}
      \Tr_{X/k}(\alpha) & \text{if } i = r,\\
      0 & \text{if } 0 \le i \le r-1,
    \end{cases}
  \]
  where $p\colon \PP(V)\to X$ denotes the projection.

  \begin{proof}
    Let $\pi\colon X\to \Spec k$ be the structure map and let $n = \dim X$. For degree reasons, we
    can assume that $\alpha\in H^n(X,\Omega^n_{X/k})$ and $i = r$ because otherwise both
    sides of the desired equality are zero. In this situation, we have
    \[
      \Tr_{\PP(V)/k}(p^\ast(\alpha)\cdot \xi^r) = \pi_\ast
      p_\ast(p^\ast(\alpha)\cdot \xi^r) = \pi_\ast(\alpha\cdot
      p_\ast(\xi^r)) = \Tr_{X/k}(\alpha \cdot p_\ast(\xi^r)),
    \]
    where we use \cite[Theorem~1 (b) and (f)]{Srinivas} for the first and last
    equality, and \cite[Theorem~1 (d)]{Srinivas} for the second
    equality. Thus, we are left to prove $p_\ast(\xi^r) = 1 \in H^0(X,\Omega^0_{X/k})$.

    We first prove that $p_\ast(\xi^r) = 1$ in $\CH^0(X)$. Since $X$ is
    irreducible, the localisation sequence yields $\CH^0(X) = \CH^0(U)$ for
    every non-empty open subscheme $U \subset X$. By restricting to some open $U \subset X$ where
    $V$ trivialises, we may thus assume that $V \cong X\times \A^{r+1}$
    and hence, $\PP(V) = \PP^r_X$. Now, $\xi$ is the first Chern class of
    $\OO_{\PP^r_X}(1)$ and thus, corresponds with the class of a
    hyperplane. Consequently, $\xi^r$ corresponds withe the intersection of $r$
    hyperplanes in general position, which is a single $X$-point with
    multiplicity one. Thus, we get $p_\ast(\xi^r) = 1 \in \CH^0(X)$. By
    applying \cite[Theorem~1 (g) and Theorem~2]{Srinivas}, we now also obtain
    $p_\ast(\xi^r) = 1 \in H^0(X,\Omega^0_{X/k})$.
  \end{proof}
\end{lemma}

\begin{proposition}
  \label{prop:eqec-pbf}
  Let $V$ be a rank $r+1$ vector bundle over a smooth, projective $k$-scheme $X$. Let $G$ act
  on $X$ and equivariantly on $V$. Then we have
  \[
    \chi_G(\PP(V)/k) = \chi_c(\PP^{r}/k)\cdot \chi_G(X/k) \in \GW^G(k).
  \]

  \begin{proof}
    Since the equivariant quadratic Euler characteristic is additive in
    disjoint unions, we may assume that $X$ is equidimensional.
    We will analyse the Hodge cohomology and trace using Proposition~\ref{prop:hodge-pbf}. The decomposition in Proposition~\ref{prop:hodge-pbf}, yields an isomorphism of
    $H^\ast(X,\Omega^\ast_{X/k})$-modules
    \begin{equation}
      \label{eq:eqec-pbf:decomp}
      H^\ast(\PP(V), \Omega^\ast_{\PP(V)/k})=\bigoplus_{i=0}^{r-1} H^{\ast-i}(X,
      \Omega^{\ast-i}_{X/k})\xi^i = H^\ast(X,\Omega^\ast_{X/k})\otimes
      H^\ast(\PP^{r}_k,\Omega^\ast_{\PP^{r-1}_k/k})
    \end{equation}
    compatible with the grading.
    Here, the isomorphism on the right hand side is induced by sending the
    first Chern class of $\OO_{\PP^{r-1}}(1)$ to $\xi$. Furthermore,
    both isomorphisms are $G$-equivariant if we equip $\PP_k^{r}$ with the
    trivial $G$-action. We equip the right hand side of
    \eqref{eq:eqec-pbf:decomp} with the $\Z$-graded $G$-equivariant non-degenerate
    quadratic form by $\beta_R$ defined as the tensor product of quadratic forms on
    $H^\ast(X,\Omega^\ast_{X/k})$ and
    $H^\ast(\PP^r_k,\Omega^\ast_{\PP^r_k/k})$ induced by the composition of
    cup product and trace.
    We have to show that $\beta_R$ and the composition of cup product and trace on
    the left hand side of \eqref{eq:eqec-pbf:decomp}, denoted $\beta_L$, induce the same element in $\GW^G(k)$.

    Let $n = \dim X$.
    By Remark \ref{rem:z-graded-nonzero-hyperbolic} and Serre duality, the contribution of
    \[
      V_L := \bigoplus_{(p,q) \ne ((n+r)/2,(n+r)/2)} H^q(\PP(V), \Omega^p_{\PP(V)/k})
    \]
    to $\chi_G(\PP(V)/k) = \can_\Z(\beta_L)$ is hyperbolic with Lagrangian
    \[
      L_L := \bigoplus_{\substack{(p,q) \ne ((n+r)/2,(n+r)/2)\\ q-p \ge
          1\text{ or } p=q>(n+r)/2}} H^q(\PP(V), \Omega^p_{\PP(V)/k}).
    \]
    Thus, this contribution is completely determined by the
    $G$-representation $L_L$. By \eqref{eq:eqec-pbf:decomp}, $L_L$ is
    $G$-equivariantly isomorphic to
    \[
      L_R:=\bigoplus_{\substack{(p,q) \ne ((n+r)/2,(n+r)/2)\\ q-p\ge 1 \text{ or } p=q>(n+r)/2}}
      (H^\ast(X,\Omega^\ast_{X/k})\otimes
      H^\ast(\PP^{r}_k,\Omega^\ast_{\PP^{r}_k/k}))^{p,q}.
    \]
    Let $x \in (H^\ast(X,\Omega^\ast_{X/k})\otimes
    H^\ast(\PP^{r}_k,\Omega^\ast_{\PP^{r}_k/k}))^{p,q}$ and $y \in (H^\ast(X,\Omega^\ast_{X/k})\otimes
    H^\ast(\PP^{r}_k,\Omega^\ast_{\PP^{r}_k/k}))^{p',q'}$ with $(p',q')\neq
    (n+r-p, n+r-q)$. So, we have $xy \in (H^\ast(X,\Omega^\ast_{X/k})\otimes
    H^\ast(\PP^{r}_k,\Omega^\ast_{\PP^{r}_k/k}))^{p+p',q+q'}$ and
    $(p+p',q+q') \ne (n+r,n+r)$. Since the
    traces involved to define $\beta_R$ send everything that is not in
    bi-degree $(n+r,n+r)$ to zero, we have $\beta_R(x,y) = 0$. Thus, the
    $(p,q)$ and $(p',q')$ bi-graded pieces on the right are orthogonal with
    respect to $\beta_R$ if $(p',q')\neq (n+r-p, n+r-q)$.
    In particular, $L_R$ is a Lagrangian for the contribution of
    \[
      V_R :=\bigoplus_{(p,q) \ne ((n+r)/2,(n+r)/2)}
      (H^\ast(X,\Omega^\ast_{X/k})\otimes
      H^\ast(\PP^{r}_k,\Omega^\ast_{\PP^{r}_k/k}))^{p,q}
    \]
    to $\can_\Z(\beta_R) \in \GW^G(k)$.
    Hence, the contributions of $V_R$ to $\can_\Z(\beta_R)$ and $V_L$ to $\can_\Z(\beta_L)$ agree.

    Thus, it is enough to show that the composition of cup product and trace
    $\beta_V \in \GW^G(k)$ on $H^{(n+r)/2}(\PP(V),
    \Omega^{(n+r)/2}_{\PP(V)/k})$ and $\beta_X \in \GW^G(k)$ on $(H^\ast(X,\Omega^\ast_{X/k})\otimes
    H^\ast(\PP^{r}_k,\Omega^\ast_{\PP^{r-1}_k/k}))^{(n+r)/2,(n+r)/2}$ induce
    the same form on $\GW^G(k)$, whenever $n+r$ is even.

    We first show this for $n = 2m$ and $r = 2s$ even. In this case, the
    isomorphism \eqref{eq:eqec-pbf:decomp} of $H^\ast(X,\Omega^\ast_{X/k})$-modules induces
    a $G$-equivariant isomorphism
    \[
      H^{r+s}(\PP(V),\Omega^{r+s}_{\PP(V)/k)}) \cong
      \bigoplus_{i=1}^sH^{m-i}(X,\Omega^{m-i}_{X/k})\xi^{s+i} \oplus
      H^m(X,\Omega^m_{X/k})\xi^s \oplus \underset{=:L}{\underbrace{\bigoplus_{i=1}^s H^{m+i}(X,\Omega^{m+i}_{X/k})\xi^{s-i}}}.
    \]
    The algebra structure from Proposition \ref{prop:hodge-pbf} implies that
    the subspace $L$ is a $G$-invariant totally isotropic
    subspace with
    orthogonal complement $L^\perp = H^m(X,\Omega^m_{X/k}) \oplus L$. Indeed,
    Lemma \ref{lem:pbf-trace} yields
    \[
      \Tr_{\PP(V)/k}(p^\ast(\alpha)\cdot \xi^i) =
      \begin{cases}
        \Tr_{X/k}(\alpha) & \text{if } i = r\\
        0 & \text{if } 0 \le i \le r-1.
      \end{cases}
    \]
    Thus, $L$ is totally isotropic and we have $H^m(X,\Omega^m_{X/k}) \oplus
    L\subset L^\perp$. Since
    \[
      \dim \left(\bigoplus_{i=1}^sH^{m-i}(X,\Omega^{m-i}_{X/k})\xi^{s+i}\right) = \dim L,
    \]
    we get
    $L^\perp = H^m(X,\Omega^m_{X/k}) \oplus L$ by a dimension count.

    Thus, we get by \cite[Lemma 2.8 (c)]{schlichting_hermitian_2010} that
    $\beta_V = \gamma + H(L) \in \GW^G(k)$ where $\gamma \in \GW^G(k)$ is the composition
    of cup product and trace
    \[
      H^m(X,\Omega^m_{X/k})\times H^m(X,\Omega^m_{X/k}) \xrightarrow{\cup}
      H^n(X,\Omega^n_{X/k}) \xrightarrow{\Tr_{X/k}} k.
    \]
    Similarly, we get $\beta_X = \gamma + H(L) \in\GW^G(k)$ and thus, $\beta_X
    =\beta_V$.

    The proof in case $n$ and $r$ are both odd is similar, with the only
    difference that the summand $H^m(X,\Omega^m_{X/k})$ does not appear and that $L$ is a
    Lagrangian.
  \end{proof}
\end{proposition}

\begin{corollary}
  \label{cor:g-eq-normal-bundle-fixed}
  Let $X$ be a smooth, projective scheme over $k$ with a $G$-action. Let $Y \subset X$
  be a smooth closed subscheme fixed by $G$. Then we have
  \[
    \chi_G(\PP(N_YX)/k) = \chi_c(\PP(N_YX)/k) \in \GW^G(k).
  \]

  \begin{proof}
    Since the equivariant quadratic Euler characteristic is additive in
    disjoint unions, we can assume that $Y$ is connected. Hence, $N_YX$ is a
    vector bundle of fixed rank $r$.
    Because the $G$-action on $Y$ is trivial, we have $\chi_G(Y/k)^G =
    \chi_c(Y/k)$. Now, Proposition~\ref{prop:eqec-pbf} yields
    \[
      \chi_G(\PP(N_YX)/k) = \chi_c(\PP^{r-1}/k)\cdot \chi_G(Y/k) =
      \chi_c(\PP^{r-1}/k)\cdot \chi_c(Y/k) = \chi_c(\PP(N_YX)/k).\qedhere
    \]
  \end{proof}
\end{corollary}

\begin{proposition}
  \label{prop:g-eq-blow-up}
  Let $X$ be a smooth, projective scheme over $k$ with a $G$-action. Let $Y \subset X$
  be a smooth, closed subscheme that is stable under the $G$-action. Then we
  have
  \[
    \chi_G(\Bl_YX/k) = \chi_G(X/k) - \chi_G(Y/k) + \chi_G(\PP(N_YX)/k) \in \GW^G(k).
  \]

  \begin{proof}
    Since the equivariant quadratic Euler characteristic and the exceptional
    divisor of blow-ups are additive in
    disjoint unions, we may assume that $X$ and $Y$ are equidimensional.
    Let $E\subset \Bl_YX/k$ be the exceptional divisor and denote the
    morphisms in the blow-up square by
    \[
      \begin{tikzcd}
        E \ar[r, "i_E"] \ar[d, "\pi"] & \Bl_YX\ar[d, "\pi"]\\
        Y \ar[r, "i_Y"] & X.
      \end{tikzcd}
    \]
    Thus, we have
    $E \cong \PP(N_YX)$ and hence, $\chi_G(E/k) = \chi_G(\PP(N_YX)/k)$.
    By Proposition \ref{prop:hodge-blowup}, we have
    \[
      \chi_G(\Bl_YX/k) = \chi_G(X/k) + q,
    \]
    where $q\in \GW^G(k)$ is the $G$-equivariant quadratic form induced by
    cup-product and trace on the $H^q(E,\Omega_{E/k}^p)/\pi^\ast H^q(Y,
    \Omega^p_{Y/k})$ part of the splitting. Thus, it remains to show that $q =
    \chi_G(E/k) - \chi_G(Y/k)$.

    By the projective bundle formula in Proposition \ref{prop:hodge-pbf} for
    $V = N_YX$ and $E = \PP(V)$, we have
    \[
      H^\ast(E,\Omega^\ast_{E/k}) = \bigoplus_{i=0}^{r}H^{\ast-i}(Y,\Omega^{\ast-i}_{Y/k})\xi^i,
    \]
    where $r+1$ denotes the codimension of $Y$ in $X$ and $\xi$ denotes the first
    Chern class of $\OO_E(1)$. This yields an isomorphism
    \[
      \frac{H^\ast(E,\Omega^\ast_{\PP(V)/k})}{\pi^\ast
        H^\ast(Y,\Omega^\ast_{Y/k})} \cong \bigoplus_{i=1}^{r}H^{\ast-i}(Y,\Omega^{\ast-i}_{Y/k})\xi^i.
    \]
    Let $a \in H^\ast(Z,\Omega^\ast Z)$.
    Since $i_Y^\ast\colon H^\ast(X,\Omega^\ast_{X/k}) \to
    H^\ast(Y,\Omega^\ast_{Z/k})$ is surjective,
    we can lift $a$ to
    an element $\hat a \in H^\ast(X,\Omega^\ast_{X/k})$ with $i_Y^\ast(\hat a)
    = a$. Furthermore, \cite[Theorem~1 (g)]{Srinivas} yields
    $i_E^\ast(-\hat\xi) = \xi$ for $\hat\xi := i_{E\ast}(1) \in
    H^1(\Bl_YX,\Omega^1_{\Bl_YX/k})$. Now, we have by \cite[Theorem~1 (d)]{Srinivas}
    \begin{align*}
      i_E^\ast(i_{E\ast}(\pi^\ast (a)\xi^\nu))
      &= i_E^\ast i_{E\ast}(\pi^\ast
      i_Z^\ast(\hat a) \cdot i_E^\ast((-\hat\xi)^\nu))\\
      &= i_E^\ast i_{E\ast}
      i_E^\ast(\pi^\ast(a) \cdot (-\hat\xi)^\nu)\\
      &= i_E^\ast(\pi^\ast(a) \cdot (-\hat\xi)^\nu\cdot i_{E\ast}(1))\\
      &= \pi^\ast (a)\xi^\nu \cdot i_E^\ast i_{E\ast}(1)\\
      &= -\pi^\ast (a)\xi^{\nu+1}
    \end{align*}
    for $0 \le \nu \le r-1$.
    Furthermore, we have $\pi_\ast i_{E\ast}(\pi^\ast (a) \xi^\nu) =
    i_{Z\ast}(a) \pi_\ast(\xi^\nu) = 0$ for $0 \le \nu < r$ by \cite[Theorem~1
    (b) and (d)]{Srinivas} and a dimension consideration.

    Let $x = \sum_{i=1}^ra_i \xi^i \in H^\ast(E,\Omega^\ast_{E/k})/\pi^\ast
    H^\ast(Z,\Omega^\ast_{Z/k})$ and write $x/\xi := \sum_{i=1}^ra_i
    \xi^{i-1} \in H^\ast(E,\Omega^\ast_{E/k})$. The above computation yields
    \[
      -i_E^\ast i_{E\ast}(x/\xi) = x.
    \]
    Since the assignment $x \mapsto x/\xi$ is additive, we have thus
    constructed an explicit right inverse of $i_E^\ast$ that induces the map $H^\ast(E,\Omega^\ast_{\PP(V)/k})/\pi^\ast
    H^\ast(Y,\Omega^\ast_{Y/k}) \to H^\ast(\Bl_YX,\Omega^\ast_{\Bl_YX/k})$
    from the splitting in Proposition \ref{prop:hodge-blowup}.

    Now, let $x = \sum_{i=1}^ra_i \xi^i,y = \sum_{i=1}^rb_i \xi^i \in H^\ast(E,\Omega^\ast_{E/k})/\pi^\ast
    H^\ast(Z,\Omega^\ast_{Z/k})$ and consider $x/\xi$ and $y/\xi$ as defined
    above as elements of $H^\ast(E,\Omega^\ast_{E/k})$. Hence, we have by
    \cite[Theorem~1]{Srinivas}
    \begin{align*}
      q(x,y)
      &= \Tr_{\Bl_YX/k}((-i_{E\ast}(x/\xi))\cdot (-i_{E\ast}(y/\xi)))\\
      &= \Tr_{\Bl_YX/k}(i_{E\ast}(x/\xi)\cdot i_{E\ast}(y/\xi))\\
      &= \Tr_{\Bl_YX/k}(i_{E\ast}(x/\xi\cdot i_E^\ast i_{E\ast}(y/\xi)))\\
      &= \Tr_{\Bl_YX/k}i_{E\ast}(x/\xi\cdot (-\xi) \cdot y/\xi)\\
      &= -\Tr_{E/k}(x/\xi \cdot y/\xi\cdot \xi).
    \end{align*}
    Thus, we need to show that the symmetric bilinear form
    \[
      \bigoplus_{i=0}^{r-1}H^\ast(Y,\Omega^\ast_{Y/k})\xi^i \times
      \bigoplus_{i=0}^{r-1}H^\ast(Y,\Omega^\ast_{Y/k})\xi^i \to
      H^{n+r}(E,\Omega^{n+r}_{E/k}) \xrightarrow{\Tr_{E/k}} k
    \]
    given by first sending $(\sum_{i=0}^{r-1}a_i\xi^i,\sum_{i=0}^{r-1}b_i\xi^i)$ to
    $(\sum_{i=0}^{r-1}a_i\xi^i\cdot \sum_{i=0}^{r-1}b_i\xi^i)\cdot (-\xi)$ and then
    applying $\Tr_{E/k}$ induces $\chi_G(E/k) - \chi_G(Z/k)$ in
    $\GW^G(k)$. Here, we equip $H^q(Y,\Omega^p_{Y/k})\xi^i$ with degree
    $q-p$ and then map the $\Z$-graded form to $\GW^G(k)$
    using the canonical map.

    By the proof of Proposition \ref{prop:eqec-pbf}, the induced form in $\GW^G(k)$ agrees with the
    form induced by
    \begin{equation}
      \label{eq:g-eq-blow-up:simplified-form}
      \bigoplus_{i=0}^{r-1}H^\ast(Y,\Omega^\ast_{Y/k})\xi^i \times
      \bigoplus_{i=0}^{r-1}H^\ast(Y,\Omega^\ast_{Y/k})\xi^i \to
      H^\ast(Y,\Omega^\ast_{Y/k})[\xi]/(\xi^{r+1}) \xrightarrow{\varphi_r} k,
    \end{equation}
    where the first map sends $(x,y)$ to $x\cdot y\cdot (-\xi)$, and where $\varphi_r$ is
    the composition of $\Tr_{Y/k}$ and the projection to the factor of
    $\xi^r$.
    Since projecting to the $\xi^r$ part of $xy\xi$ is the same as projecting
    to the $\xi^{r-1}$-part of $xy$, we get that
    the $\Z$-graded form \eqref{eq:g-eq-blow-up:simplified-form} agrees with the form
    \begin{equation}
      \label{eq:g-eq-blow-up:pr-form}
      \bigoplus_{i=0}^{r-1}H^\ast(Y,\Omega^\ast_{Y/k})\xi^i \times
      \bigoplus_{i=0}^{r-1}H^\ast(Y,\Omega^\ast_{Y/k})\xi^i \to
      H^\ast(Y,\Omega^\ast_{Y/k})[\xi]/(\xi^{r}) \xrightarrow{\varphi_{r-1}} k,
    \end{equation}
    where the first map sends $(x,y)$ to $xy$ and $\varphi_{r-1}$ is the
    composition of $-\Tr_{Y/k}$ and the projection to the factor of
    $\xi^{r-1}$. Again by the proof of Proposition \ref{prop:eqec-pbf}, we get
    that the form \eqref{eq:g-eq-blow-up:pr-form} induces $\langle -1\rangle
    \cdot\chi_c(\PP^{r-1}/k)\cdot \chi_G(Y/k)$. Thus, we have
    \[
      q = \langle -1\rangle
      \cdot\chi_c(\PP^{r-1}/k)\cdot \chi_G(Y/k) = (\chi_c(\PP^r/k) - \langle
      1\rangle) \cdot \chi_G(Y/k) = \chi_G(E/k) - \chi_G(Y/k)
    \]
    by Proposition \ref{prop:eqec-pbf}.
  \end{proof}
\end{proposition}

\begin{corollary}
  \label{cor:g-eq-blow-up-fixed}
  Let $X$ be a smooth, projective scheme over $k$ with a $G$-action. Let $Y \subset X$
  be a smooth, closed subscheme of constant codimension $r+1$ fixed by $G$. Then we
  have
  \[
    \chi_G(\Bl_YX/k) = \chi_G(X/k) + (\chi_c(\PP^r/k) - \langle 1\rangle)\chi_c(Y/k).
  \]

  \begin{proof}
    Combine Corollary \ref{cor:g-eq-normal-bundle-fixed} and Proposition
    \ref{prop:g-eq-blow-up}.
  \end{proof}
\end{corollary}

\subsection{Application to the $\A^1$-Euler Characteristic of Quotients}
We shall now apply the $G$-equivariant quadratic Euler characteristic
to the study of quotient schemes. That is, consider the following setup: Let $X$ be a smooth, projective scheme
over $k$ and let $G$ act on $X$. Let $Y = X/G$ be the quotient and denote the
quotient map by $\pi\colon X \to Y$.
We want to computed the $\A^1$-Euler characteristic of $Y$
from  the $G$-equivariant $\A^1$-Euler characteristic of $X$ and possibly
other information.

The following lemma already appeared in
\cite[Lemma~21]{bv2024quadratic}. Since this is one of the building blocks for
comparing Euler characteristics, we include its proof here for the convenience
of the reader.

\begin{lemma}[{\cite[Lemma~21]{bv2024quadratic}}]
  \label{lem:trace-finite-change}
  Let $f\colon X \to Y$ be a finite morphism of degree $d$ between smooth,
  projective, $n$-dimensional schemes $X$ and $Y$
  over $k$. Then the diagram
  \[
    \begin{tikzcd}
      H^n(X, \Omega^n_{X/k}) \ar[r, "\Tr_{X/k}"] & k.\\
      H^n(Y, \Omega^n_{Y/k}) \ar[ur, "d\cdot \Tr_{Y/k}"']\ar[u,
      "f^\ast"] &
    \end{tikzcd}
  \]
  commutes.

  \begin{proof}
    Consider the pushforward $f_\ast\colon H^n(X,\Omega^n_{X/k})\to H^n(Y,\Omega^n_{Y/k})$
    constructed by Srinivas \cite[Theorem 1]{Srinivas}. The pushforward of the
    respective structure maps agrees with the trace
    maps by \cite[Theorem~1~(e)]{Srinivas}. Hence, the diagram
    \[
      \begin{tikzcd}
        H^n(X, \Omega^n_{X/k}) \ar[r, "\Tr_{X/k}"]\ar[d, "f_\ast"] & k\\
        H^n(Y, \Omega^n_{Y/k}) \ar[ur, "\Tr_{Y/k}"'] &
      \end{tikzcd}
    \]
    commutes by \cite[Theorem~1~(b)]{Srinivas}. Since we have $f_\ast f^\ast = d\cdot \operatorname{id}$ by \cite[Theorem~1~(d) and
    (f)]{Srinivas}, we get the statement of the lemma.
  \end{proof}
\end{lemma}

The following proposition is well-known to the experts. We include a proof
here for the convenience of the reader.

\begin{proposition}
  \label{prop:quotient-hodge-cohomology}
  Assume that $G$ is a group of order prime to
  $\characteristic k$ acting freely on an open dense subscheme of $X$ and assume that $Y$ is smooth. Then the pullback map $\pi^\ast \colon
  H^q(Y,\Omega^p_{Y/k}) \to H^q(X,\Omega^p_{X/k})$ induces an isomorphism
  \[
    H^q(Y,\Omega^p_{Y/k}) \cong (H^q(X,\Omega^p_{X/k}))^G.
  \]
  Since this isomorphism is canonical, we shall just write $H^q(Y,\Omega^p_{Y/k}) =
  (H^q(X,\Omega^p_{X/k}))^G$ considering both groups to be the same.

  \begin{proof}
    Since $X$ and $Y$ are smooth, the quotient map $\pi$ is flat. Hence $\pi_\ast\Omega^p_{X/k}$ is a locally free coherent sheaf on $Y$. Since $|G|$ is prime to $\characteristic k$, we get that $G$ acts semi-simply on $\pi_\ast\Omega^p_{X/k}$. In particular, $(\pi_\ast\Omega^p_{X/k})^G$ is a direct summand of $\pi_\ast\Omega^p_{X/k}$.

    Furthermore, we have $\Omega^p_{Y/k} =
    (\pi_\ast\Omega^p_{X/k})^G$ by Knighten \cite[Corollary after Theorem 3 in
    Chapter II]{Knighten1973DiifQuotients}, since we can check this equality
    after passing to an algebraic closure, and the ramification is tame, since
    $|G|$ is prime to the characteristic of $k$.

    Since the $G$-representation
    $\pi_\ast\Omega^p_{X/k}$ is semisimple, we get a decomposition into
    isotypical components
    \[
      \pi_\ast\Omega^p_{X/k} = \bigoplus_\rho (\pi_\ast\Omega^p_{X/k})^\rho
    \]
    where $\rho$ ranges over all irreducible
    characters. Note that this decomposition makes sense since the
    $G$-action on $Y$ is trivial.

    Now, we have
    \[
      \bigoplus\nolimits_\rho H^q(Y,\pi_\ast\Omega^p_{X/k})^\rho =
      H^q(Y,\pi_\ast\Omega^p_{X/k}) =
      H^q(Y,\bigoplus\nolimits_\rho(\pi_\ast\Omega^p_{X/k})^\rho) =
      \bigoplus\nolimits_\rho H^q(Y,(\pi_\ast\Omega^p_{X/k})^\rho),
    \]
    where $\bigoplus_\rho H^q(Y,\pi_\ast\Omega^p_{X/k})^\rho$ is again the
    decomposition into isotypical components. Let $\mathfrak U = \{U_i\}$ be an open affine
    cover of $Y$.
    Hence, the Čech complex $\check C^\bullet(\mathfrak U, -)$
    computes sheaf cohomology since $Y$ is again projective and hence,
    separated. The decomposition into isotypical components of
    $\pi_\ast\Omega^p_{X/k}$ induces the decomposition
    \[
      \bigoplus_\rho \check C^\bullet(\mathfrak U,
      \pi_\ast\Omega^p_{X/k})^\rho = \bigoplus_\rho \check C^\bullet(\mathfrak U,
      (\pi_\ast\Omega^p_{X/k})^\rho)
    \]
    of $C^\bullet(\mathfrak U, \pi_\ast\Omega^p_{X/k})$ into isotypical
    components with $C^\bullet(\mathfrak U, \pi_\ast\Omega^p_{X/k})^\rho = \check C^\bullet(\mathfrak U,
      (\pi_\ast\Omega^p_{X/k})^\rho)$. Because the Čech complex
    computes sheaf cohomology and because taking isotypical components is
    exact, this induces an equality $H^q(Y,(\pi_\ast\Omega^p_{X/k})^\rho)=
    H^q(Y,\pi_\ast\Omega^p_{X/k})^\rho$ for every $\rho$.

    In particular, this yields $H^q(Y,\Omega^p_{Y/k}) = H^q(Y,(\pi_\ast\Omega^p_{X/k})^G) =
    H^q(Y,\pi_\ast\Omega^p_{X/k})^G$.
    The map $\pi^\ast$ induces a $G$-equivariant isomorphism
    $H^q(Y,\pi_\ast\Omega^p_{X/k})\cong H^q(X,\Omega^p_{X/k})$ because $\pi$ is
    finite. Hence
    \[
      H^q(Y,\Omega^p_{Y/k}) = H^q(Y,\pi_\ast\Omega^p_{X/k})^G\cong H^q(X,\Omega^p_{X/k})^G.\qedhere
    \]
  \end{proof}
\end{proposition}

\begin{theorem}
  \label{thm:smooth-quotient-euler-char}
  Let $X$ be a smooth, projective scheme over $k$ and let $G$ be a group of order prime to $\characteristic k$ acting on $X$. Let $Y = X/G$ be the quotient and $\pi\colon X \to Y$
  the quotient map. Assume that $G$ acts freely on an open, dense subset of $X$ and assume that $Y$ is smooth.
  Then $\chi_c(Y/k) = \langle|G|\rangle \cdot \chi_G(X/k)^G$.

  \begin{proof}
    Since the equivariant quadratic Euler characteristic is invariant under taking
    disjoint unions, we can assume that $X$ is equidimensional of dimension $n$.
    Consider the quadratic forms
    \[
      \chi(X/k)_{p,q} \colon (H^p(X,\Omega^q_{X/k})\oplus H^{n-p}(X,\Omega^{n-q}_{X/k}))
      \times (H^p(X,\Omega^q_{X/k})\oplus H^{n-p}(X,\Omega^{n-q}_{X/k}))
      \to k
    \]
    given by $\chi(X/k)_{p,q}((a,b),(a',b'))=\Tr_{X/k}(a\cup b'+(-1)^{(p+q)(n+1)}b\cup a')$
    for $0 \le p < q \le n$ or $p = q$ with $0 \le 2p < n$. If $n$ is even,
    also consider the quadratic form
    \[
      \chi(X/k)_{n/2,n/2} \colon H^{n/2}(X,\Omega^{n/2}_{X/k})
      \times H^{n/2}(X,\Omega^{n/2}_{X/k})
      \xrightarrow{\cup} H^n(X,\Omega^n_{X/k}) \xrightarrow{\Tr_{X/k}} k.
    \]
    These again carry a $G$-action. Define $\chi(Y/k)_{p,q}$ similarly.
    By Theorem \ref{thm:motivic-gauss-bonnet}, we have
    \[
      \chi_c(Y/k) = \sum_{p,q}(-1)^{q-p}\chi(Y/k)_{p,q}
    \]
    and by definition, we have
    \begin{equation}
      \label{eq:smooth-quotient-euler-char:x}
      \chi_G(X/k) = \sum_{p,q}(-1)^{q-p}\chi(X/k)_{p,q}.
    \end{equation}
    The map $\pi$ is finite of degree $|G|$ because $G$ acts freely on an
    open, dense subscheme of $X$.
    Since $\pi^\ast$ induces an isomorphism $H^p(Y,\Omega^q_{Y/k}) \cong
    H^p(X,\Omega^q_{X/k})^G$ by Proposition \ref{prop:quotient-hodge-cohomology}, we have $\chi(Y/k)_{p,q} =
    \langle|G|\rangle\cdot \chi(X/k)_{p,q}^G$ by Lemma
    \ref{lem:trace-finite-change}. If we now apply $(-)^G$ to
    \eqref{eq:smooth-quotient-euler-char:x}, we get the statement.
  \end{proof}
\end{theorem}

Let $Z \subset X$ be the fixed-point locus of the $G$-action.
By the universal property of blow-ups, the $G$-action on $X$ induces a
$G$-action on $\Bl_ZX$. Thus, $\Bl_ZX/G \to Y = X/G$ is proper and an isomorphism
over $(X\setminus Z)/G$. In a special case
of this situation, we can say the following:

\begin{proposition}
  \label{prop:euler-char-quot-blow-up}
  Let $G$ be a cyclic group of prime order $p \ne \characteristic k$ and
  assume that the induced $G$-action on the exceptional divisor $E$ of
  $\Bl_ZX$ is
  trivial. Furthermore, assume that the quotient $(\Bl_ZX)/G$ is smooth and that
  $Z$ has constant codimension $c$ in $X$. Then
  the compactly supported $\A^1$-Euler characteristic of $Y$ is given by
  \[
    \chi_c(Y/k) = \chi_G(X/k)^G\cdot \langle |G|\rangle + (\langle |G|\rangle -
    \langle 1 \rangle) \cdot \langle -1\rangle \cdot \chi_c(\PP^{c-2}/k)\cdot
    \chi_c(Z/k)\in \GW(k).
  \]

  \begin{proof}
    Since $X$, $\Bl_Z X$, and $(\Bl_ZX)/G$ are smooth and projective, we can use Theorem \ref{thm:smooth-quotient-euler-char}
    to compare the $\A^1$-Euler characteristics of $\Bl_ZX$ and $(\Bl_ZX)/G$.

    The $G$-action on $Z$ is trivial because the $G$-action on $E$ is trivial
    and $\pi\colon E\to Z$ is $G$-equivariant. By Corollary \ref{cor:g-eq-blow-up-fixed}, we have
    \begin{align*}
      \chi_G(\Bl_ZX/k)
      &= \chi_G(X/k) + (\chi_c(\PP^{c-1}/k) - \langle
        1\rangle)\chi_c(Z/k)\\
      &= \chi_G(X/k) + \langle -1\rangle \cdot \chi_c(\PP^{c-2}/k)\chi_c(Z/k).
    \end{align*}
    Now, Theorem \ref{thm:smooth-quotient-euler-char} yields
    \begin{align*}
      \chi_c(((\Bl_ZX)/G)/k)
      &= \langle |G|\rangle\cdot\chi_G(((\Bl_ZX)/G)/k)^G\\
      &=\langle |G|\rangle\cdot \chi_G(X/k)^G + \langle |G|\rangle\cdot \langle -1\rangle \cdot \chi_c(\PP^{c-2}/k)\chi_c(Z/k).
    \end{align*}
    Similar to the situation on $\Bl_ZX$, we also have since $Z/G = Z$ and
    $E/G = E$
    \[
      \chi_c(((\Bl_ZX)/G)/k) = \chi_c(Y/k) + \langle -1\rangle \cdot \chi_c(\PP^{c-2}/k)\chi_c(Z/k).
    \]
    Comparing the above to formulae, yields the statement.
  \end{proof}
\end{proposition}

\begin{corollary}
  \label{cor:euler-char-quot-blow-up-cases}
  In the notation and under the assumptions of Proposition \ref{prop:euler-char-quot-blow-up}, we have
  \[
    \chi_c(Y/k) = \chi_G(X/k)^G\cdot \langle |G|\rangle +
    \begin{cases}
      0 &\text{if } \dim Z \text{ odd or } c \text{ odd},\\
       (\langle |G|\rangle -
    \langle 1 \rangle) \cdot \langle -1\rangle \cdot \chi_c(Z/k) & \text{otherwise}.\\
    \end{cases}
  \]
  \begin{proof}
    This follows from the observation that $(\langle |G|\rangle -
    \langle1\rangle)$ multiplied with a hyperbolic form is zero, the explicit
    description of $\chi_c(\PP^{c-2}/k) = \sum_{i=0}^{c-2} \langle
    (-1)^i\rangle$ and the observation that $\chi_c(Z/k)$ is a hyperbolic form if $\dim
    Z$ is odd, by Corollary \ref{cor:motivic-gb}.
  \end{proof}
\end{corollary}

\begin{corollary}
  \label{cor:z2-action-euler-char}
  Let $G = \Z/2\Z$ act on a smooth, projective scheme $X$ over $k$. Then we have
  \[
    \chi_c((X/G)/k) = \langle 2\rangle \cdot \chi_G(X/k)^G +
    \begin{cases}
      0 &\text{if } \dim Z \text{ odd or } c \text{ odd},\\
      (\langle -2\rangle -
    \langle -1 \rangle) \cdot \chi_c(X^G/k) & \text{otherwise}.
    \end{cases}
  \]

  \begin{proof}
    By Corollary~\ref{cor:smooth-z2-blow-up-quotient}, the assumptions of
    Corollary~\ref{cor:euler-char-quot-blow-up-cases} are satisfied when we
    consider the blow-up of $X^G$ in $X$. Thus, the
    statement follows from Corollary~\ref{cor:euler-char-quot-blow-up-cases}.
  \end{proof}
\end{corollary}

\section{\texorpdfstring{Explicit Computation for $\Sym^2$}{Explicit Computation for Sym2}}
\label{sec:explicit-sym2}
We shall now compute the $\A^1$-Euler characteristic of $\Sym^2X$ for a connected, smooth, projective scheme $X$ over $k$.

\begin{proposition}
  \label{prop:fixed-points-quadratic-action}
  Let $V$ be an $n$-dimensional $k$-vector space and let $\beta \colon V\otimes
  V\to k$ be a non-degenerate symmetric bilinear form on $V$. Assume that
  $\beta$ is isometric to the diagonal form $\sum_{i=1}^n\langle \alpha_i\rangle$. Let $G
  \subset S_m$ be a subgroup and consider the non-degenerate symmetric bilinear form
  \[
    \beta^m\colon (V^{\otimes m})\otimes (V^{\otimes m})\to k;
    (v_1\otimes \dots \otimes v_m)\otimes (w_1\otimes \dots \otimes w_m)\mapsto
    \beta(v_1,w_1)\cdots \beta(v_m,w_m).
  \]
  Then the following holds:
  \begin{enumerate}
  \item Let $G$ act on $V^{\otimes m}$ by permuting the factors; that is,
    $\sigma \in G$ acts on $v_1\otimes\dots \otimes v_m \in V^{\otimes m}$ by $\sigma\cdot v_1\otimes \dots \otimes v_m =
    v_{\sigma(1)}\otimes \dots \otimes v_{\sigma(m)}$. Then $\beta^m$ is a
    non-degenerate $G$-equivariant symmetric bilinear form and we have
    \[
      (\beta^m)^{G} = \sum_{G\cdot (i_\ast)\in \{1, \dots, n\}^m/G}\langle
      |G\cdot (i_\ast)| \cdot \alpha_{i_1}\cdots \alpha_{i_m}\rangle
    \]
    where $G$ acts on $\{1, \dots, n\}^m$ by permuting the factors, $(i_\ast)$
    is shorthand for $(i_1, \dots, i_n) \in \{1, \dots, n\}^m$, and $|G\cdot (i_\ast)|$ denotes the cardinality of the orbit of $(i_\ast)$ under the
    $G$-action.
    \label{prop:fixed-points-quadratic-action:no-sign}
  \item Let $G$ act on $V^{\otimes m}$ by permuting the factors and
    introducing a sign; that is, $\sigma \in G$ acts on $v_1\otimes \dots
    \otimes v_m \in V^{\otimes m}$ by $\sigma\cdot v_1\otimes \dots \otimes v_m =
    \sgn \sigma \cdot v_{\sigma(1)}\otimes \dots \otimes v_{\sigma(m)}$.
    Then $\beta^m$ is a
    non-degenerate $G$-equivariant symmetric bilinear form and we have
    \[
      (\beta^m)^{G} = \sum_{G\cdot (i_\ast)\in (\{1, \dots, n\}^m/G)'}\langle
      |G\cdot (i_\ast)| \cdot \alpha_{i_1}\cdots \alpha_{i_m}\rangle
    \]
    where $G$ acts on $\{1, \dots, n\}^m$ by permuting the factors, and the
    sum runs over the orbits $G\cdot (i_\ast)$ where the stabiliser of some
    (or equivalently all)
    representative $(i_\ast)$ of $G\cdot (i_\ast)$ is contained $A_m$.
    \label{prop:fixed-points-quadratic-action:sign}
  \end{enumerate}

  \begin{proof}
    The condition in (\ref{prop:fixed-points-quadratic-action:sign}) that
    there is a representative $(i_\ast)$ of $G\cdot (i_\ast)$ such that the
    stabiliser $\Stab_G(i_\ast)$ of $(i_\ast)$ is contained in $A_m$ does not depend on the
    choice of representative. If $(i'_\ast)$ is another representative of
    $G\cdot (i_\ast)$, then $\Stab_G(i'_\ast)$ and $\Stab_G(i_\ast)$ are
    conjugated. Since $A_m$ is normal, we thus have $\Stab_G(i_\ast) \subset
    A_m$ if and only if $\Stab_G(i'_\ast) \subset A_m$. In particular, the
    condition in (\ref{prop:fixed-points-quadratic-action:sign}) only depends
    on the orbit.

    Note, that $\beta^m$ is non-degenerate and $G$-equivariant in the
    situation of (\ref{prop:fixed-points-quadratic-action:no-sign}) and
    (\ref{prop:fixed-points-quadratic-action:sign}).
    Let $v_1, \dots, v_n\in V$ be a basis of $V$ such that
    \[
      \beta(v_i,v_j) =
      \begin{cases}
        \alpha_i &\text{if } i=j,\\
        0 & \text{if } i\ne j.
      \end{cases}
    \]
    Therefore, $\{v_{i_1}\otimes \dots\otimes v_{i_m}\mid i_1,
    \dots,i_m \in \{1, \dots, n\} \}$ is a basis of $V^{\otimes m}$ and
    we have
    \begin{equation}
      \label{eq:fixed-points-quadratic-action}
      \beta^m(v_{i_1}\otimes \dots\otimes v_{i_m},v_{j_1}\otimes \dots\otimes
      v_{j_m}) =
      \begin{cases}
        \alpha_{i_1}\cdots \alpha_{i_m} & \text{if } (i_1, \dots, i_m) = (j_1, \dots, j_m),\\
        0 & \text{if }(i_1, \dots, i_m) \ne (j_1, \dots, j_m).
      \end{cases}
    \end{equation}
    We now prove (\ref{prop:fixed-points-quadratic-action:no-sign}).
    In this situation, $(V^{\otimes m})^G$ has a basis consisting of the elements
    \[
      v_{G\cdot (i_\ast)} = \sum_{(j_1, \dots, j_m) \in G\cdot (i_\ast)} v_{j_1}\otimes
      \dots\otimes v_{j_m}
    \]
    where $G\cdot (i_\ast)$ runs over the orbits in $\{1, \dots,
    n\}^m/G$. Thus, we get using \eqref{eq:fixed-points-quadratic-action}
    \begin{align*}
      \beta^m(v_{G\cdot (i_\ast)}, v_{G\cdot (i_\ast)})
      &= \sum_{(j_1,
        \dots, j_m) \in G\cdot (i_\ast)}\sum_{(j'_1,
        \dots, j'_m) \in G\cdot (i_\ast)} \beta^m(v_{j_1}\otimes
        \dots\otimes v_{j_m}, v_{j_1'}\otimes
        \dots\otimes v_{j_m'})\\
      &= \sum_{(j_1,\dots, j_m) \in G\cdot (i_\ast)} \alpha_{j_1}\cdots \alpha_{j_m}\\
      &= |G\cdot (i_\ast)|\cdot \alpha_{i_1}\cdots \alpha_{i_m}.
    \end{align*}
    The last equality holds because
    the group action only permutes the elements $i_1, \dots,
    i_m$. Since \eqref{eq:fixed-points-quadratic-action} also yields
    \[
      \beta^m(v_{G\cdot (i_\ast)}, v_{G\cdot (j_\ast)}) = 0
    \]
    for $G\cdot (i_\ast) \ne G\cdot (j_\ast)$,
    we get the claim.

    The proof of (\ref{prop:fixed-points-quadratic-action:sign}) works
    similarly: let $v = \sum_{(i_\ast)\in \{1,\dots,
      n\}^m}\lambda_{(i_\ast)}v_{i_1}\otimes \dots\otimes v_{i_m} \in
    (V^{\otimes m})^G$. Let $(i_\ast) \in \{1,\dots, n\}^m$ and
    $\sigma \in \Stab_{(i_\ast)}G$. Now, $\sigma v = v$ implies
    $\lambda_{(i_\ast)} = \sgn\sigma \cdot \lambda_{(i_\ast)}$. In particular, if
    there exists $\sigma \in \Stab_{(i_\ast)}G$ with $\sgn\sigma = -1$, then
    $\lambda_{(i_\ast)} = 0$, since we are assuming that $k$ has
    characteristic different from two. Thus, we can describe a basis of $(V^{\otimes
      m})^G$ as follows: let $(i^1_\ast), \dots, (i^N_\ast)$ be a system of
    representatives of the orbits in $(\{1, \dots, n\}^m/G)'$, where $(\{1, \dots,
    n\}^m/G)'$ consists of all orbits of $\{1, \dots, n\}^m/G$ whose
    stabiliser is contained in the kernel of $\sgn$. The elements
    \[
      v_{(i^\nu_\ast)} = \sum_{\sigma \in G/\Stab_{(i^\nu_\ast)}G} \sigma\cdot
      v_{i_1^\nu}\otimes\dots \otimes v_{i_m^\nu}
    \]
    for $\nu = 1, \dots, N$ form a basis of $(V^{\otimes m})^G$. Since we get
    just as in the proof of (\ref{prop:fixed-points-quadratic-action:no-sign})
    that
    \[
      \beta^m(v_{(i^\nu_\ast)}, v_{(i^\nu_\ast)})
      = |G\cdot (i^\nu_\ast)|\cdot \alpha_{i^\nu_1}\cdots \alpha_{i^\nu_m}
    \]
    and
    \[
      \beta^m(v_{(i^\nu_\ast)}, v_{(i^{\nu'}_\ast)}) = 0
    \]
    for $\nu \ne \nu'$, we get the claim.
  \end{proof}
\end{proposition}

\begin{lemma}
  \label{lem:sym2-vector-space-comp}
  Let $V^\bullet$ be a bounded finite-dimensional $\Z$-graded $k$-vector space and let $\beta\colon
  V^\bullet\times V^\bullet \to k$ be a $\Z$-graded non-degenerate symmetric bilinear form on $V^\bullet$. Let
  $\Z/2\Z$ act on $V^\bullet\otimes V^\bullet$ by $a\otimes b \mapsto
  (-1)^{\nu\nu'}b\otimes a$ for $a\in V^\nu$ and $b \in V^{\nu'}$, and
  assume that $\beta^\even$ is isometric to $\sum_{i=1}^n\langle\alpha_i\rangle$. Then we have
  \[
    (\can_\Z(\beta^2))^{\Z/2\Z} = n\langle 1\rangle + \sum_{1\le i<j\le
      n}\langle 2\alpha_i\alpha_j\rangle  + m\langle -2\rangle + (m^2 - (n+1)m)\cdot H \in \GW(k),
  \]
  where $\dim V^\odd = 2m$.

  \begin{proof}
    Note that $\can_\Z\beta = \sum_{i=1}^n\langle \alpha_i\rangle-mH$ and that
    $\dim V^\even = \rank \beta^\even = n$. Furthermore, we can find a
    diagonalising basis
    $v_1^\odd, \dots, v_{2m}^\odd$ of $V^\odd$ for $\beta^\odd$ with diagonal
    elements $\alpha'_i = (-1)^i$.

    We have $(\beta^2)^\even = ((\beta^\even)^2 + (\beta^\odd)^2)$ and
    $(\beta^2)^\odd = (\beta^\even\beta^\odd + \beta^\odd\beta^\even)$ as $\Z/2\Z$-equivariant
    bilinear forms by Remark \ref{rem:even-odd-mult}. The group
    $\Z/2\Z$ acts on $(\beta^\even)^2$ and $\beta^\even\beta^\odd +
    \beta^\odd\beta^\even$ by permuting the factors, and on
    $(\beta^\odd)^2$ by permuting the factors with a sign. Thus,
    we can use Proposition \ref{prop:fixed-points-quadratic-action} to compute
    $((\beta^2)^\even)^{\Z/2\Z}$ as
    \begin{align*}
      ((\beta^2)^\even)^{\Z/2\Z}
      &= ((\beta^\even)^2)^{\Z/2\Z} +
      ((\beta^\odd)^2)^{\Z/2\Z}\\
      &= n\langle 1\rangle + \sum_{1 \le i <
        j \le n}\langle 2\alpha_i\alpha_j\rangle + \sum_{1\le i<j\le
        2m}\langle 2(-1)^{i+j}\rangle\\
      &= n\langle 1\rangle + \sum_{1 \le i <
        j \le n}\langle 2\alpha_i\alpha_j\rangle + m\langle -2\rangle + (m^2-m)H.
    \end{align*}

    We are thus left with computing
    $((\beta^2)^\odd)^{\Z/2\Z}$.
    Let $v_1^\even, \dots, v_n^\even$ be a basis of $V^\even$
    and let $v_1^\odd, \dots, v_{2m}^\odd$ be a basis of $V^\odd$. Now,
    $\Z/2\Z$ acts on $V^\even\otimes V^\odd \oplus V^\odd \otimes V^\even$ by
    sending $v_i^\even\otimes v_j^\odd$ to $v_j^\odd\otimes v_i^\even$ and $v_j^\odd\otimes v_i^\even$ to $v_i^\even\otimes v_j^\odd$. Thus, a basis of
    $((V\otimes V)^\odd)^{\Z/2}$ is given by the elements
    \[
      w_{ij} := v_i^\even \otimes v_j^\odd + v_j^\odd\otimes v_i^\even
    \]
    for $i = 1, \dots, n$ and $j = 1, \dots, 2m$. In particular, $\dim ((V\otimes
    V)^\odd)^{\Z/2} = 2nm$. Thus, $((\beta^2)^\odd)^{\Z/2\Z}$
    contributes $-\frac 12 \dim ((V\otimes
    V)^\odd)^{\Z/2}\cdot H = -nmH$ to $(\can_\Z(\beta^2))^{\Z/2\Z}$ in $\GW(k)$.
  \end{proof}
\end{lemma}

\begin{theorem}
  \label{thm:qec-sym2}
  Let $X$ be a connected, smooth, projective scheme over $k$. Express the $\A^1$-Euler
  characteristic of $X$ as $\chi_c(X/k) = \beta -mH$ with
  $\beta = \sum_{i=1}^n\langle \alpha_i\rangle$ and $m\ge 0$.
  Then if $\dim X$ is odd, we have in $\GW(k)$
  \begin{equation}
    \label{thm:qec-sym2:odd}
    \chi_c(\Sym^2X/k) = n\langle 2\rangle + \sum_{1\le i<j\le n}\langle
    \alpha_i\alpha_j\rangle + m\langle -1\rangle + (m^2 - (n+1)m)\cdot H.
  \end{equation}
  If $\dim X$ is even, we have
  \begin{equation}
    \label{thm:qec-sym2:even}
    \begin{split}
      \chi_c(\Sym^2X/k) &= n\langle 2\rangle + \sum_{1\le i<j\le n}\langle
                          \alpha_i\alpha_j\rangle +
                          m\langle -1\rangle + (m^2 - (n+1)m)\cdot H\\
      &\quad + (\langle -2\rangle -
      \langle -1 \rangle) \cdot \chi_c(X/k)
    \end{split}
  \end{equation}
  in $\GW(k)$.

  \begin{proof}
    First, note that the formulae \eqref{thm:qec-sym2:odd} and
    \eqref{thm:qec-sym2:even} are independent of the choice of $m$ in the representation of
    $\chi_c(X/k)$ in the statement of the theorem. Indeed, since
    \eqref{thm:qec-sym2:odd} and \eqref{thm:qec-sym2:even} only differ by the
    summand $(\langle -2\rangle -
    \langle -1 \rangle) \cdot \chi_c(X/k)$, it is enough to show the
    independence for \eqref{thm:qec-sym2:odd}. If we change the above representation as
    \[
      \chi_c(X/k) = \beta -mH = (\beta + H) -(m+1)H,
    \]
    formula \eqref{thm:qec-sym2:odd} changes to
    \begin{align*}
      &(n+2)\langle 2\rangle + \sum_{1\le i< j\le n}\langle
        \alpha_i\alpha_j\rangle + \underset{=nH}{\underbrace{\sum_{i=0}^n\langle \alpha_i\rangle + \sum_{i=0}^n\langle
        -\alpha_i\rangle}} + \langle -1\rangle + (m+1)\langle -1\rangle\\
      &\quad+ ((m+1)^2 - (n+3)(m+1))\cdot H\\
      =\;& n\langle 2\rangle + \sum_{1\le i < j\le n}\langle \alpha_i\alpha_j\rangle +
         m\langle -1\rangle + nH +  \underset{= 2H}{\underbrace{2\langle 2
         \rangle + 2\langle -1\rangle}}\\
      &\quad +(m^2 + 2m + 1 - (n+1)m -n -2m - 3)H\\
      =\;& n\langle 2\rangle + \sum_{1\le i< j\le n}\langle
         \alpha_i\alpha_j\rangle + m\langle -1\rangle + (m^2 - (n+1)m)\cdot H.
    \end{align*}

    We want to use Corollary~\ref{cor:z2-action-euler-char} to compute the
    $\A^1$-Euler characteristic of $\Sym^2X$, that is, we have to compute
    $\chi_c(X\times X/k)^{\Z/2}$.

    The formulae \eqref{thm:qec-sym2:odd} and \eqref{thm:qec-sym2:even} are
    invariant under changing $m$ in the representation of $\chi_c(X/k)$ in the
    statement of the theorem. Thus, we can assume by Example
    \ref{ex:graded-euler-char} that $\beta$ is the even part of the trace form on Hodge
    cohomology
    \begin{equation}
      \label{eq:qec-sym2:trace-form}
      \left(\bigoplus_{i,j=0}^{\dim_kX}
        H^i(X,\Omega^j_{X/k})[j-i],\Tr\right)
    \end{equation}
    from Example \ref{ex:graded-euler-char} and the odd part of
    \eqref{eq:qec-sym2:trace-form} is isometric to $mH$.
    With this explicit representation of $\chi_c(X/k)$, Lemma
    \ref{lem:sym2-vector-space-comp} yields an explicit description of
    $\chi_c(X\times X/k)^{\Z/2}$.
    The result now follows from Corollary~\ref{cor:z2-action-euler-char} after
    observing that the codimension of the diagonal in $X\times X$ equals the
    dimension of $X$ and its $\A^1$-Euler characteristic is hyperbolic if $X$
    is odd-dimensional.
  \end{proof}
\end{theorem}

\begin{corollary}
  Let $X$ be an odd-dimensional connected, smooth, projective scheme over $k$ with
  $\A^1$-Euler characteristic $\chi_c(X/k) = mH$ for some $m \in \Z$. Then
  we have
  \[
    \chi_c(\Sym^2X/k) =
    \begin{cases}
      m\langle 1\rangle + m^2H &\text{if } m \ge 0,\\
      -m\langle -1\rangle + (m^2 + m)\cdot H &\text{if } m < 0.
    \end{cases}
  \]

  \begin{proof}
    First note that we can write $\chi_c(X/k) = mH$ for every odd-dimensional
    smooth, projective scheme by Corollary \ref{cor:motivic-gb}. The corollary now follows by plugging
    $mH$ into formula \eqref{thm:qec-sym2:odd} of Theorem \ref{thm:qec-sym2}
    and note that $\langle 1\rangle + \langle 1\rangle  = \langle 2\rangle
    + \langle 2\rangle$.
  \end{proof}
\end{corollary}

\section{\texorpdfstring{Explicit Computation for $\Sym^3$}{Explicit Computation for Sym3}}
\label{sec:explicit-sym3}
Assume that $\characteristic k \ne 2,3$ throughout this section and let $X$ be
an $n$-dimensional, connected, smooth quasi-projective scheme over $k$ throughout this section. When we
turn towards the actual computation of $\A^1$-Euler characteristics, we will
additionally assume that $X$ is projective.

We now want to compute the compactly supported $\A^1$-Euler characteristic of
$\Sym^3X$ in a similar manner as in the case for $\Sym^2X$. That is, we
equip $X^3$ with the action of $S_3$ permuting the factors, and
we relate the $\A^1$-Euler characteristic of $X$ with the compactly
supported $\A^1$-Euler characteristic of $\Sym^3X = X^3/S_3$.

In order to find a sequence of blow-ups whose $S_3$-quotient is smooth, we
start with an intermediary step: consider the subgroup
$A_3\subset S_3$ consisting of all permutations with sign $1$. These are the
identity and the two cyclic permutations, so $A_3\cong \Z/3\Z$. We first find a
sequence of blow-ups of $X^3$ such that the $A_3$-quotient $Y/A_3$ of the
final blow-up $Y$ is smooth. The $A_3$-quotient carries an induced action of
$S_3/A_3\cong \Z/2\Z$. Then we shall blow-up $Y/A_3$ such that the
$\Z/2\Z$-quotient is smooth by Corollary \ref{cor:smooth-z2-blow-up-quotient}. We shall relate this in turn back to $\Sym^3X = X^3/S_3$ so
that we can compute its compactly supported $\A^1$-Euler
characteristic.

We can lift the $\Z/2\Z$-action on the quotient $X^3/A_3$ to $X^3$ in three different
ways: we can permute the first two factors, the second two factors, or the first and the third
factor. We denote the $\Z/2\Z$-action fixing the $i$-th factor and permuting the
other two by $\tau_i$ for $i = 1,2,3$. Let $\Delta_i = \{(x_1,x_2,x_3) \in
X^3\mid x_j = x_k \text{ for } j \ne i \ne
k\} \subset X^3$ be one of the big diagonals. Thus, we have $\Delta_i =
(X^3)^{\tau_i}$.

Consider the blow-up $\Bl_\Delta X^3 \to X^3$ of $X^3$ along the diagonal
$\Delta \subset X^3$ and denote its exceptional divisor by $E_\Delta\subset
\Bl_\Delta X^3$. Furthermore, denote the blow-up of $\Delta_i$ along $\Delta$ by
$\Delta_i'$. The closed immersion $\Delta_i \hookrightarrow X^3$ lifts to a
closed immersion $\Delta_i' \hookrightarrow
\Bl_\Delta X^3$ identifying $\Delta_i'$ with the strict transform of
$\Delta_i$ in $\Bl_\Delta X^3$.

In what follows, we refer to a projective space bundle of the form $\PP_Y(V)
\to Y$, for $V$ a rank $n+1$ vector bundle over some scheme $Y$ as a
\emph{Zariski locally trivial $\PP^n$-bundle over $Y$}.

\begin{lemma}
  \label{lem:a3-fixed-points-bldelta}
  Let $R = k[T]/(T^2 + T + 1)$. Then the fixed-point locus $Z$ of the induced
  $A_3$-action on $\Bl_\Delta X^3$ is a Zariski locally trivial
  $\PP^{n-1}$-bundle over $\Delta
  \times \Spec R$. Moreover, $Z\subset E_\Delta$ and $\tau_i$ acts on $Z$ by
  permuting the roots of $T^2+T+1$ in $R$. Furthermore, $Z/\tau_i$ is a
  Zariski locally trivial
  $\PP^{n-1}$-bundle over $\Delta$.

  \begin{proof}
    The observation that $Z \subset E_\Delta$ follows from the observation that there are no $A_3$-fixed
    points on $X^3\setminus \Delta\cong \Bl_\Delta X^3 \setminus E_\Delta$. Therefore, we only need to consider the
    fixed points contained in $E_\Delta$, which is equivariantly isomorphic to the
    projectivisation of the normal bundle $N_\Delta X^3$ of $\Delta$ in $X^3$. Thus, we need
    to compute the eigenspaces of the action on $N_\Delta X^3$.

    Let $i\colon \Delta \to X^3$ be the closed immersion.
    We have an equivariant isomorphism $N_\Delta X^3 = i^\ast T(X^3)/T\Delta \cong
    (TX)^3/TX$ via the fibre-wise diagonal embedding of $TX$ into $TX^3$. The inclusion
    of the first two summands yields an isomorphism
    $TX^2 \cong TX^3/TX$, where the $A_3$-action restricts to the action by
    the block-matrix
    \[
      A:=
      \begin{pmatrix}
        0 & -I_n\\
        I_n & -I_n
      \end{pmatrix};
    \]
    here, we have chosen the cyclic permutation $(1\;2\;3)$ as a generator of $A_3$.
    This matrix has characteristic polynomial $(T^2 + T + 1)^n$. Thus, for every
    closed point $x\in \Delta$, we have an
    $n$-dimensional eigenspace, defined by the zero locus of $T \id -
    A$ over $R$. These eigenspaces assemble into a rank $n$ vector
    bundle over $V \to \Delta \times \Spec R$.
    The vector bundle $V\to \Delta \times \Spec R$ is naturally a sub-vector
    bundle of  the pullback $p_1^*N_\Delta X^3$ where $p_1\colon \Delta \times
    \Spec R\to \Delta$ is the projection to the first factor. The two
    roots of $T^2 + T + 1$ are the two non-trivial third roots $\zeta$ and
    $\zeta^2$ of $1$ in $\bar{k}$ and thus, for each $x\in \Delta$, the two
    eigenspaces in $N_\Delta X^3\otimes_k\bar{k}$ only intersect at
    $0$. Hence, the projection $p_{1*}(V)\subset N_\Delta X^3$ is a sub-bundle
    of $N_\Delta X^3$ whose fibre over $x\in \Delta$ is the union of the two
    $n$-dimensional eigenspaces, intersecting only at 0. If $\zeta$ is already
    in $k(x)$, then each of these eigenspaces is a $k(x)$-linear subspace of
    $N_\Delta X^3\otimes k(x)$, and if not, the two eigenspaces are $k(x,
    \zeta)$-linear subspaces of $N_\Delta\otimes k(x,\zeta)$ that are
    conjugate over $k(x)$. In either case, the map
    $\tilde{p}_1:\PP_{\Delta\times\Spec R}(V)\to \PP_\Delta(N_\Delta X^3)$
    induced by $p_1\colon \Delta \times \Spec R\to \Delta$ is a closed
    immersion. Furthermore,$\tilde{p}_1$ identifies $\PP_{\Delta\times\Spec
      R}(V)$ with the $A_3$-fixed locus $Z\subset E_\Delta$ under the
    identification of $\PP_\Delta(N_\Delta X^3)$ with the exceptional divisor
    $E_\Delta$. Thus, $Z$ is a Zariski locally trivial $\PP^{n-1}$-bundle over
    $\Delta\times\Spec R$.

    To get the statement about the $\tau_i$-action, note the following: for every
    closed point $x \in \Delta$, the $\tau_i$-action permutes the two
    $n$-dimensional eigenspaces, since conjugation by $\tau_i$ sends
    $(1\;2\;3)$ to its inverse $(3\;2\;1)$. This translates to permuting the roots of
    $T^2+T+1$ in $R$, which is the described action in the statement.
  \end{proof}
\end{lemma}

\begin{lemma}
  \label{lem:blowup-2-exceptional-action-EZ}
  Let $X$ be a smooth, quasi-projective scheme over $k$.
  Let $\Bl_Z\Bl_\Delta X^3$ and $E'_\Delta = \Bl_ZE_\Delta$ be the blow-ups of
  $Z$ in $\Bl_\Delta X^3$ and $E_\Delta$, respectively. Then the exceptional divisors $E_Z$ and
  $E_Z\cap E'_\Delta$ of $\Bl_Z\Bl_\Delta X^3$ and $E'_\Delta$, respectively,
  are equal to the $A_3$-fixed points of $\Bl_Z\Bl_\Delta X^3$ and
  $E'_\Delta$, respectively. In particular, $(\Bl_Z\Bl_\Delta X^3)/A_3$ and $E'_\Delta/A_3$ are
  smooth by Theorem \ref{thm:divisor-smooth-quotient}.

  \begin{proof}
    Because $Z$ is the fixed-point locus of the $A_3$-action on $\Bl_\Delta
    X^3$ and $E'_\Delta$, we have $E_Z \subset \Bl_Z\Bl_\Delta X^3$ $E_Z\cap
    E'_\Delta\subset (E'_\Delta)^{A_3}$.
    Since we also have $E_Z\cap E'_\Delta \subset E_Z$, it thus suffices to show that $E_Z$
    is fixed under the $A_3$-action.

    We first prove the statement for $X = \A^n$. Note that it is enough to
    prove the statement after a faithfully flat base-change, so we may assume
    that $k$ is algebraically closed. Let $\zeta \in k$ be a primitive third root of
    unity. By diagonalising the action on $(\A^n)^3$, we may rewrite
    the $A_3$-action on $(\A^n)^3$ so that the cyclic permutation $\sigma = (1\; 2\; 3)$ acts on
    $(x,y,z) \in (\A^n)^3$ by $(x,\zeta y, \zeta^2 z)$.

    Consider the blow-up $\Bl_0(\A^n)^2$ of $0$ in $(\A^n)^2$ and act on
    $(\A^n)^2$ via $\sigma\cdot
    (x,y) = (\zeta x, \zeta^2y)$ for $(x,y) \in (\A^n)^2$. Let $Z' \subset \Bl_0(\A^n)^2$ be the
    fixed-point locus of the induced $A_3$-action. Thus, $Z' \subset E_0$
    where $E_0\subset \Bl_0(\A^n)^2$ is the exceptional divisor. In order to
    prove the lemma for $\A^n$, it is enough to prove that the $A_3$-action on
    the exceptional divisor of the blow-up of $Z'$ in $\Bl_0(\A^n)^2$ is trivial by
    \cite[\href{https://stacks.math.columbia.edu/tag/0805}{Tag~0805}]{stacks-project}.

    By \cite[Chapter II, Example 7.12.1]{HartshorneAG}, we can model the
    blow-up $\Bl_0(\A^n)^2$ as the vanishing locus of $\{x_iy_j - x_jy_i\mid
    i,j = 1, \dots, 2n\}$ in $\PP^{2n-1}_A = \Proj A[Y_1, \dots, Y_{2n}]$ where $A
    = k[x_1, \dots, x_{2n}]$. Consider the distinguished open subset $D(Y_i)$ for
    some $i \in \{1, \dots, 2n\}$. Hence, $\Bl_0(\A^n)^2\cap D(Y_i)$ is
    isomorphic to the spectrum of
    \[
      \frac{A[y_1, \dots, \hat y_i,\dots, y_{2n}]}{(x_j-x_iy_j\mid j = 1, \dots
        2m; j \ne i)} = k[x_i, 
    y_1, \dots, \hat y_i,\dots, y_{2n}],
    \]
    where $\hat y_i$ indicates that we leave out the variable $y_i$.
    
    Now, let  $1 \le i \le n$. The relation $x_j-x_iy_j$ on $D(Y_i)\cap \Bl_0(\A^n)^2$, shows that $\sigma
    y_j = y_j$ for $1\le j \le n$ and $\sigma y_j = \zeta y_j$ for $n+1 \le
    j \le 2n$. Thus, $Z'\cap D(Y_i)$ is the vanishing locus of the ideal
    $I = (x_i, y_{n+1}, \dots, y_{2n})$. Since $\sigma$ acts by multiplication
    with $\zeta$ on all generators of $I$, we get that $\sigma$ acts by
    multiplication with $\zeta$ on $I/I^2$, which is isomorphic to the normal
    bundle $(N_{Z'}\Bl_0(\A^n)^2)|_{D(Y_i)}$ of $Z'\cap D(Y_i)$ in
    $\Bl_0(\A^n)^2\cap D(Y_i)$. Thus, $A_3$ acts by scalar multiplication on
    $(N_{Z'}\Bl_0(\A^n)^2)|_{D(Y_i)}$ and hence, trivially on
    $\PP(N_{Z'}\Bl_0(\A^n)^2)|_{D(Y_i)}$.

    One can show analogously that the $A_3$-action of
    $\PP(N_{Z'}\Bl_0(\A^n)^2)|_{D(Y_i)}$ is trivial for $i = n+1, \dots,
    2n$. Since the $D(Y_i)$ cover $\Bl_0(\A^n)^2$, we get that the
    $A_3$-action on $\PP(N_{Z'}\Bl_0(\A^n)^2)$, which is the exceptional
    divisor of $\Bl_{Z'}\Bl_0(\A^n)^2$, is trivial. This concludes the proof
    for $\A^n$.

    Now, let $X$ be any smooth, connected, quasi-projective scheme over $k$. It suffices
    to show that the
    $A_3$-action on $E_Z$ is trivial locally. Since $E_Z$ lies in the fibre of
    $\Bl_Z\Bl_\Delta X^3 \to X^3$ over $\Delta$, it is thus enough to show the following: for every $x
    \in X$ there exists a Zariski-open neighbourhood $U \subset X$ such that the
    lemma holds for $U$. Indeed, if we let $\pi\colon \Bl_Z\Bl_\Delta X^3 \to
    \Bl_\Delta X^3 \to X^3$ be the projection map and $x \in E_Z$, then
    such a neighbourhood $U$ for $y \in X$ with $\pi(x) = (y,y,y)$ shows
    that the $A_3$-action on the open neighbourhood $\pi^{-1}(U^3)\cap E_Z$ of
    $x$ is trivial.

    Let $x \in X$ be a point. Since $X$ is smooth, there is an open subscheme
    $U \subset X$ and an étale map $\varphi\colon U \to \A^n$. We therefore
    get an $A_3$-equivariant étale morphism $\psi\colon U^3 \to
    (\A^n)^3$. For every point $x' \in U$, we have $(N_\Delta
    U^3)_{(x',x',x')} = (N_{\varphi^{-1}(\Delta)}U^3)_{(x',x',x')}$. After
    restricting to an $A_3$-stable open neighbourhood $V \subset U^3$ of the
    diagonal $\Delta \subset U^3$, we thus get an $A_3$-equivariant isomorphism $Z \to
    (\Bl_{\psi^{-1}(\Delta)}V)^{A_3}$. By a similar argument, we get that
    $E_Z$ is $A_3$-equivariantly isomorphic to the exceptional divisor of
    $\Bl_Z\Bl_{\varphi^{-1}(\Delta)}V$. Since we already know the lemma for
    $\A^n$, we thus also know the lemma for $U$ since blowing up commutes with
    flat base-change by
    \cite[\href{https://stacks.math.columbia.edu/tag/0805}{Tag~0805}]{stacks-project}.
  \end{proof}
\end{lemma}

\begin{lemma}
  \label{lem:blowup-3-exceptional-action-FiEi}
  Let $F_i = (\Bl_\Delta X^3)^{\tau_i}$. Then $E_i := E_\Delta\cap F_i =
  (E_\Delta)^{\tau_i}$ consists of two disjoint Zariski locally trivial $\PP^{n-1}$-bundles over
  $\Delta$. Furthermore, the $A_3$-action permutes the $F_i$ and we have $F_i
  \cap F_j = \emptyset$ for $i \ne j$ and $Z\cap F_i =
  \emptyset$. Finally, one of the two components of $E_i$ is $\Delta_i'\cap
  E_\Delta$; letting $\tilde E_i$ denote the other component, $F_i$ is a
  disjoint union of $\Delta_i'$ and $\tilde E_i$.

  \begin{proof}
    Since $\Bl_\Delta X^3\setminus E_\Delta \cong X^3\setminus \Delta$ and
    $\Delta_i\cap \Delta_j = \Delta$ for $i \ne j$, it is
    enough to prove the assertion in $E_\Delta$. For this, we use as in Lemma
    \ref{lem:a3-fixed-points-bldelta} that $TX^2 \cong TX^3/T\Delta$ via the
    inclusion into the first two factors. Now, $\tau_1$, $\tau_2$, and
    $\tau_3$ act via the block matrices
    \[
      \begin{pmatrix}
        I_n & -I_n\\
        0 & -I_n
      \end{pmatrix}, \quad
      \begin{pmatrix}
        -I_n & 0\\
        -I_n & I_n
      \end{pmatrix}, \quad \text{and}\quad
      \begin{pmatrix}
        0 & I_n\\
        I_n & 0
      \end{pmatrix},
    \]
    respectively. These actions have the
    following eigenspaces:

    \begin{center}
      \begin{tabular}{c|c|c}
        Eigenvalue & $1$ & $-1$\\\hline
        $i= 1$ & $\{(v,0)\mid v\in TX\}$ & $\{(v,2v)\mid v\in TX\}$\\
        $i=2$ & $\{(0,v)\mid v\in TX\}$ & $\{(2v,v)\mid v\in TX\}$ \\
        $i=3$ & $\{(v,v)\mid v\in TX\}$ & $\{(v,-v)\mid v\in TX\}$
      \end{tabular}
    \end{center}
    
    For $i = 1,2,3$, the projectivisation of the eigenspace for the eigenvalue
    $1$ is the
    exceptional divisor of $\Delta_i'$ and the projectivisation of the
    eigenspace for the eigenvalue $-1$ is $\tilde E_i$. Since all of these
    eigenspaces pairwise intersect trivially and they intersect trivially
    with the subspace
    defining $Z$ as described in the proof of Lemma \ref{lem:a3-fixed-points-bldelta}, we get the claim.
  \end{proof}
\end{lemma}

We now have the following picture:
\begin{equation}
  \label{eq:sym3-blow-up-sequence}
  \begin{tikzcd}[column sep=tiny]
    X^3 & & \Bl_\Delta X^3\ar[ll] & & \Bl_Z\Bl_\Delta X^3\ar[ll] & &
    \Bl_{\coprod_iF_i}\Bl_Z\Bl_\Delta X^3\ar[ll]\\
    \Delta_i \ar[u, hook] & & F_i = (\Bl_\Delta X^3)^{\tau_i} \ar[u, hook] \ar[ll] & & F_i \ar[ll,
    "\sim"]\ar[u, hook]
    & & E_{F_i}\ar[ll]\ar[u, hook] \\
    \Delta \ar[u, hook] & &  E_\Delta \ar[ll] \ar[uu, bend left, hook] & Z =
    (E_\Delta)^{A_3}\ar[l, hook] & E_Z \ar[l] \ar[uu, bend left, hook] & &
  \end{tikzcd}
\end{equation}
Here $E_Z$ is the exceptional divisor of $\Bl_Z\Bl_\Delta X^3$ and $E_{F_i}$
is the preimage of $F_i$ under the map $\Bl_{\coprod_iF_i}\Bl_Z\Bl_\Delta X^3
\to \Bl_Z\Bl_\Delta X^3$. In particular $E_{F_i}\cong
\PP(N_{F_i}\Bl_Z\Bl_\Delta X^3)$ and the exceptional divisor of
$\Bl_{\coprod_iF_i}\Bl_Z\Bl_\Delta X^3$ is $\coprod_i E_{F_i}$. Furthermore,
we will consider the blow-up $\Bl_{\coprod_i E_i}E_\Delta$ where $E_i = F_i
\cap E_\Delta$. Denote the
preimage of $E_i$ under the projection $\Bl_{\coprod E_i}E_\Delta \to
E_\Delta$ by $E_{E_i}$. In particular, $E_{E_i} = \PP(N_{E_i}E_\Delta)$ and the
exceptional divisor of $\Bl_{\coprod E_i}E_\Delta$ is $\coprod_iE_{E_i}$.

We are now turning towards computing the $\A^1$-Euler characteristic of
$\Sym^3X$; so assume that $X$ is smooth and projective from now on.
We can summarise the connections between the various schemes involved in the computation in the following
proposition.

\begin{proposition}
  \label{porp:properties-sym3-comp}
  The above schemes have the following relations, where all projective bundles
  are Zariski locally trivial:
  \begin{enumerate}\renewcommand{\theenumi}{\arabic{enumi}}
  \item $\Delta_i'\cap \Delta_j' = \emptyset$ for $i \ne j$,
    \label{porp:properties-sym3-comp:delta-ij}
  \item $F_i\cap F_j = \emptyset$ for $i \ne j$,
    \label{porp:properties-sym3-comp:F-ij}
  \item $Z\cap F_i = \emptyset$,
    \label{porp:properties-sym3-comp:ZF}
  \item $E_\Delta$ is a $\PP^{2n-1}$-bundle over $\Delta$,
    \label{porp:properties-sym3-comp:E-Delta}
  \item $Z$ is a $\PP^{n-1}$-bundle over $\Delta\times \Spec R$,
    \label{porp:properties-sym3-comp:Z}
  \item $Z/\tau_i$ is a $\PP^{n-1}$-bundle over $\Delta$,
    \label{porp:properties-sym3-comp:Z-tau}
  \item $F_i$ is the disjoint union of $\Delta_i'$ and a $\PP^{n-1}$-bundle
    over $\Delta$, denoted $\tilde E_i$,
    \label{porp:properties-sym3-comp:F}
  \item $E_{F_i}$ is the disjoint union of a $\PP^{n-1}$-bundle over $\Delta'_i$
    and a $\PP^n$-bundle over $\tilde E_i$,
    \label{porp:properties-sym3-comp:EF}
  \item $\Delta_i'\cap E_\Delta$ is a $\PP^{n-1}$-bundle over $\Delta$,
    \label{porp:properties-sym3-comp:Delta-E-Delta}
  \item $E_i$ is a disjoint union of two $\PP^{n-1}$-bundles over $\Delta$,
    \label{porp:properties-sym3-comp:E}
  \item The $A_3$-action permutes the $E_i$, and the action by $\tau_i$ on $E_i$
    is trivial,
    \label{porp:properties-sym3-comp:E-action}
  \item The $A_3$-action on $E_Z$ and $Z$ is trivial and the $\tau_i$-action
    on $E_Z$ and $Z$ is free,
    \label{porp:properties-sym3-comp:EZ-action}
  \item $E_{E_i}$ is a $\PP^{n-1}$-bundle over $E_i$,
    \label{porp:properties-sym3-comp:EE}
  \item The $\tau_i$-action on $E_{F_i}$ and $E_{E_i}$ is trivial.
    \label{porp:properties-sym3-comp:EE-action}
  \item $E_Z$ is a $\PP^n$-bundle over $Z$,
    \label{porp:properties-sym3-comp:EZ}
  \item $E_Z/\tau_i$ is a $\PP^n$-bundle over $Z/\tau_i$,
    \label{porp:properties-sym3-comp:EZ-tau}
  \item $E_Z\cap E_\Delta'$ is a $\PP^{n-1}$-bundle over $Z$ where $E_\Delta'
    = \Bl_ZE_\Delta$, and
    \label{porp:properties-sym3-comp:EZ-EDelta}
  \item $(E_Z\cap E_\Delta')/\tau_i$ is a $\PP^{n-1}$-bundle over $Z/\tau_i$;
    \label{porp:properties-sym3-comp:EZ-EDelta-tau}
  \end{enumerate}
  and the following $\A^1$-Euler characteristics
  \begin{enumerate}\renewcommand{\theenumi}{\alph{enumi}}
  \item $\chi_c(\Delta/k) = \chi_c(X/k)\in \GW(k)$,
    \label{porp:properties-sym3-comp:chi-delta}
  \item $\chi_c(\Delta_i/k) = \chi_c(X/k)^2 \in \GW(k)$,
    \label{porp:properties-sym3-comp:chi-delta-i}
  \item $\chi_{S_3}(Z/k) = \chi_c(\PP^{n-1}/k)\chi_{S_3}(\Spec R/k)\chi_c(X/k)\in \GW^{S_3}(k)$,
    \label{porp:properties-sym3-comp:chi-Z}
  \item $\chi_{S_3}(Z/k)^{\tau_i} = \langle 2\rangle \cdot
    \chi_c(\PP^{n-1}/k)\chi_c(X/k)\in \GW(k)$,
    \label{porp:properties-sym3-comp:chi-Z-tau}
  \item $\chi_{S_3}(E_Z/k) = \chi_c(\PP^n/k)\chi_c(\PP^{n-1}/k)\chi_{S_3}(\Spec R/k)\chi_c(X/k) =
    \frac{n(n+1)}{2}\chi_{S_3}(\Spec R/k)\cdot \chi_c(X/k)\cdot H\in \GW^{S_3}(k)$,
    \label{porp:properties-sym3-comp:chi-EZ}
  \item $\chi_{S_3}(E_Z/k)^{\tau_i} = \langle 2\rangle \cdot \chi_c(\PP^n/k)\chi_c(\PP^{n-1}/k)\chi_c(X/k)
    = \frac{n(n+1)}{2}\cdot \chi_c(X/k)\cdot H\in \GW(k)$,
    \label{porp:properties-sym3-comp:chi-EZ-tau}
  \item $\chi_{S_3}(E_Z\cap E_\Delta'/k) = \chi_c(\PP^{n-1}/k)^2\chi_{S_3}(\Spec
    R/k)\chi_c(X/k)\in \GW^{S_3}(k)$,
    \label{porp:properties-sym3-comp:chi-EZ-EDelta}
  \item $\chi_{S_3}(E_Z\cap E_\Delta'/k)^{\tau_i} = \langle 2\rangle \cdot
    \chi_c(\PP^{n-1}/k)^2\chi_c(X/k)\in \GW(k)$,
    \label{porp:properties-sym3-comp:chi-EZ-EDelta-tau}
  \item $\chi_{S_3}(E_\Delta/k) = \chi_c(\PP^{2n-1}/k)\cdot \chi_c(X/k) = n\cdot
    \chi_c(X/k)\cdot H \in \GW^{S_3}(k)$,
    \label{porp:properties-sym3-comp:chi-EDelta}
  \item $\chi_{\tau_i}(E_i/k) = 2\chi_c(\PP^{n-1}/k)\chi_c(X/k)\in \GW^{\Z/2\Z}(k)$,
    \label{porp:properties-sym3-comp:chi-Ei}
  \item $\chi_{\tau_i}(E_{E_i}/k) = \chi_c(\PP^{n-1}/k)\chi_c(E_i/k) =
    2\chi_c(\PP^{n-1}/k)^2\chi_c(X/k) \in \GW^{\Z/2\Z}(k)$,
    \label{porp:properties-sym3-comp:chi-EEi}
  \item $\chi_c(\Delta_i'/k) = (\chi_c(X/k) + \chi_c(\PP^{n-1}/k) - \langle
    1\rangle)\chi_c(X/k)\in \GW(k)$,
    \label{porp:properties-sym3-comp:chi-EDeltai}
  \item $\chi_c(F_i/k) = \chi_c(X/k)^2 - \chi_c(X/k) + 2\chi_c(\PP^{n-1}/k)\cdot \chi_c(X/k) =
    (\chi_c(X/k)-\langle 1\rangle +2\chi_c(\PP^{n-1}/k))\cdot \chi_c(X/k)$ in $\GW(k)$, and
    \label{porp:properties-sym3-comp:chi-Fi}
  \item $\chi_c(E_{F_i}/k) = (\langle (-1)^n\rangle - \langle 1\rangle +
    \chi_c(X/k) + 2\chi_c(\PP^{n-1}/k))\chi_c(\PP^{n-1}/k)\chi_c(X/k) =
    \langle(-1)^n\rangle \cdot \chi_c(\PP^{n-1}/k)\chi_c(X/k) +
    \chi_c(\PP^{n-1}/k)\chi_c(F_i/k)\in\GW(k)$.
    \label{porp:properties-sym3-comp:chi-EFi}
  \end{enumerate}

  \begin{proof}
    Properties (\ref{porp:properties-sym3-comp:delta-ij}) --
    (\ref{porp:properties-sym3-comp:ZF}), (\ref{porp:properties-sym3-comp:F}),
    (\ref{porp:properties-sym3-comp:E}), and (\ref{porp:properties-sym3-comp:E-action}) are a restatement of Lemma
    \ref{lem:blowup-3-exceptional-action-FiEi}. Properties~(\ref{porp:properties-sym3-comp:E-Delta}),
    (\ref{porp:properties-sym3-comp:EF}),
    (\ref{porp:properties-sym3-comp:Delta-E-Delta}),
    (\ref{porp:properties-sym3-comp:EE}), and
    (\ref{porp:properties-sym3-comp:EZ-EDelta}) hold by their respective
    definitions and general properties of blow-ups.
    Properties~(\ref{porp:properties-sym3-comp:Z}),
    (\ref{porp:properties-sym3-comp:Z-tau}), (\ref{porp:properties-sym3-comp:EZ}),
    and (\ref{porp:properties-sym3-comp:EZ-tau}) are consequences of Lemma~\ref{lem:a3-fixed-points-bldelta}. Property
    (\ref{porp:properties-sym3-comp:EZ-action}) is a consequence of Lemma~\ref{lem:blowup-2-exceptional-action-EZ} and Lemma~\ref{lem:blowup-3-exceptional-action-FiEi}. Property~(\ref{porp:properties-sym3-comp:EZ-EDelta-tau}) is a consequence of Property
    (\ref{porp:properties-sym3-comp:EZ-action}). Property~\ref{porp:properties-sym3-comp:EE-action}
    follows from Lemma~\ref{lem:z2-normal-action}. The assertion that the
    $\PP^n$-bundle in Property~\ref{porp:properties-sym3-comp:EZ-tau} and the
    $\PP^{n-1}$-bundle in
    Property~\ref{porp:properties-sym3-comp:EZ-EDelta-tau} are Zariski locally
    trivial is deduced using Hilbert's Theorem 90 to descend the relevant
    vector bundles, as in the proof of Lemma~\ref{lem:a3-fixed-points-bldelta}.

    For Properties (\ref{porp:properties-sym3-comp:chi-delta}) and
    (\ref{porp:properties-sym3-comp:chi-delta-i}), note that $\Delta \cong X$ and
    $\Delta_i \cong X^2$.

    For Properties (\ref{porp:properties-sym3-comp:chi-Z}) --
    (\ref{porp:properties-sym3-comp:chi-EFi}), we use Proposition
    \ref{prop:eqec-pbf} and Corollary
    \ref{cor:g-eq-normal-bundle-fixed}. For
    Properties~(\ref{porp:properties-sym3-comp:chi-Z-tau}),
    (\ref{porp:properties-sym3-comp:chi-EZ-tau}), and (\ref{porp:properties-sym3-comp:chi-EZ-EDelta-tau}), we additionally use
    Theorem~\ref{thm:smooth-quotient-euler-char} to compute it from the Euler
    characteristic of the quotient. For Property~(\ref{porp:properties-sym3-comp:chi-EDeltai}), we additionally use
    Proposition~\ref{prop:g-eq-blow-up}.
  \end{proof}
\end{proposition}

\begin{lemma}
  \label{lem:edelta-euler-char}
  We have
  \[
    \chi_c((E_\Delta/S_3)/k) = n\cdot \chi_c(X/k)\cdot H \in \GW(k).
  \]

  \begin{proof}
    Consider the following diagram of blow-ups
    \[
      \begin{tikzcd}
        E_\Delta & E_\Delta' = \Bl_ZE_\Delta\ar[l] & E_\Delta'' = \Bl_{\coprod_i
          E_i}E_\Delta' \ar[l]\\
        Z\ar[u, phantom, "\rotatebox{90}{$\subseteq$}"] & E_Z\cap E_\Delta'; \coprod_iE_i \ar[l]
        \ar[u,"\rotatebox{90}{$\subseteq$}", phantom] & \coprod_i E_{E_i}\ar[l]\ar[u, "\rotatebox{90}{$\subseteq$}",phantom]
      \end{tikzcd}
    \]
    where the second row describes the exceptional divisors. We get the same
    sequence of abstract blow-ups when taking the quotient with respect to the
    $S_3$-action. Therefore, we get
    \begin{align*}
      \chi_c((E_\Delta/S_3)/k)
      &= \chi_c((E_\Delta'/S_3)/k) - \chi_c((E_Z\cap
        E_\Delta'/S_3)/k) + \chi_c((Z/S_3)/k)\\
      &= \chi_c((E_\Delta''/S_3)/k)
        - \chi_c(((\coprod\nolimits_i E_{E_i})/S_3)/k) + \chi_c(((\coprod\nolimits_iE_i)/S_3)/k)\\
      &\quad - \chi_c((E_Z\cap E_\Delta'/S_3)/k) + \chi_c((Z/S_3)/k)\\
    \end{align*}

    By Proposition \ref{porp:properties-sym3-comp}
    (\ref{porp:properties-sym3-comp:E-action}), the
    $A_3$-action on $(\coprod_iE_i)$ only permutes the $E_i$, and thus, we
    have $(\coprod_iE_i)/S_3 = E_i/\tau_i = E_i$, since the $\tau_i$-action on
    $E_i$ is trivial. Similarly, we get $(\coprod_i E_{E_i})/S_3 =
    E_{E_i}$ by Proposition \ref{porp:properties-sym3-comp} (\ref{porp:properties-sym3-comp:EE-action}).

    Since by Proposition \ref{porp:properties-sym3-comp} (\ref{porp:properties-sym3-comp:EZ-action}) the action of $A_3$ on $E_Z\cap
    E_\Delta'$ is trivial and the action of $\tau_i$ is free, we get $(E_Z\cap
    E_\Delta')/S_3 = (E_Z\cap E_\Delta')/\tau_i$ and the quotient is smooth and
    projective. Thus, Theorem \ref{thm:smooth-quotient-euler-char} yields $\chi_c(((E_Z\cap E_\Delta')/S_3)/k) = \chi_c(((E_Z\cap
    E_\Delta')/\tau_i)/k) = \chi_{\tau_i}(E_Z\cap E_\Delta'/k)^{\tau_i} \cdot \langle
    2\rangle$. By Proposition \ref{porp:properties-sym3-comp} (\ref{porp:properties-sym3-comp:chi-EZ-EDelta-tau}), we have $\chi_{\tau_i}(E_Z\cap
    E_\Delta'/k)^{\tau_i} = \langle 2\rangle \cdot \chi_c(\PP^{n-1}/k)^2\chi_c(X/k)$, and thus,
    $\chi_c(((E_Z\cap E_\Delta')/S_3)/k) = \chi_c(\PP^{n-1}/k)^2\chi_c(X/k)$. Similarly, we get $\chi_c((Z/S_3)/k) =
    \chi_c(\PP^{n-1}/k)\chi_c(X/k)$.

    Thus, we have using Proposition \ref{porp:properties-sym3-comp}
    (\ref{porp:properties-sym3-comp:chi-Ei}) and (\ref{porp:properties-sym3-comp:chi-EEi})
    \begin{equation}
      \label{eq:edelta-euler-char:edelta}
      \begin{split}
        \chi_c((E_\Delta/S_3)/k)
        &= \chi_c((E_\Delta''/S_3)/k)
          - \chi_c(E_{E_i}/k) + \chi_c(E_i/k)\\
        &\quad - \chi_c(\PP^{n-1}/k)^2\chi_c(X/k) +
          \chi_c(\PP^{n-1}/k)\chi_c(X/k)\\
        &= \chi_c((E_\Delta''/S_3)/k)
          + (\langle 1\rangle - \chi_c(\PP^{n-1}/k))\chi_c(E_i/k)\\
        &\quad  + (\langle 1\rangle - \chi_c(\PP^{n-1}/k))\chi_c(\PP^{n-1}/k)\chi_c(X/k)\\
        &= \chi_c((E_\Delta''/S_3)/k)
          + 3(\langle 1\rangle - \chi_c(\PP^{n-1}/k))\chi_c(\PP^{n-1}/k)\chi_c(X/k).
      \end{split}
    \end{equation}

    We know that $E_\Delta''/S_3 =
    (\Bl_{E_i}(E_\Delta'/A_3))/\tau_i$. Furthermore, the action of $\tau_i$ on
    $E_{E_i}$ is trivial by Proposition
    \ref{porp:properties-sym3-comp} (\ref{porp:properties-sym3-comp:EE-action}) and
    $E_\Delta'/A_3$ and $\Bl_{E_i}(E_\Delta'/A_3)$ are smooth.
    Thus, we have using Theorem \ref{thm:smooth-quotient-euler-char}
    \begin{equation}
      \label{eq:edelta-euler-char:edeltapp}
      \begin{split}
        \chi_c((E_\Delta''/S_3)/k)
        &= \chi_c(((\Bl_{E_i}(E_\Delta'/A_3))/\tau_i)/k)\\
        &= \langle 2\rangle \chi_{\tau_i}((\Bl_{E_i}(E_\Delta'/A_3))/k)^{\tau_i}\\
        &= \langle 2\rangle (\chi_{\tau_i}((E_\Delta'/A_3)/k) - \chi_{\tau_i}(E_i/k) +
          \chi_{\tau_i}(E_{E_i}/k))^{\tau_i}\\
        &= \langle 6\rangle \chi_{S_3}(E_\Delta'/k)^{S_3} - \langle 2\rangle \cdot (\chi_{\tau_i}(E_i/k) -
          \chi_{\tau_i}(E_{E_i}/k))^{\tau_i}.
      \end{split}
    \end{equation}
    Here, we also used that the exceptional divisor of $\Bl_{E_i}(E_\Delta'/A_3)$
    is isomorphic to $E_{E_i}$ because the $A_3$ action only permutes the
    $E_i$ and the blow-up formula from Proposition
    \ref{prop:g-eq-blow-up}. By construction, the $\tau_i$-action on $E_i$ is
    trivial and moreover the $\tau_i$-action on $E_{E_i}$ is trivial. Thus,
    \begin{equation}
      \label{eq:edelta-euler-char:ei}
      \chi_{\tau_i}(E_i/k)^{\tau_i} = \chi_c(E_i/k) =
      2\chi_c(\PP^{n-1}/k)\chi_c(X/k)
    \end{equation}
    and
    \begin{equation}
      \label{eq:edelta-euler-char:eei}
      \chi_{\tau_i}(E_{E_i}/k)^{\tau_i} =
      \chi_c(E_{E_i}/k) = 2\chi_c(\PP^{n-1}/k)^2\chi_c(X/k).
    \end{equation}
    Now, we have
    \[
      \chi_{S_3}(E_\Delta'/k) = \chi_{S_3}(E_\Delta/k) - \chi_{S_3}(Z/k) + \chi_{S_3}(E_Z\cap E_\Delta'/k)
    \]
    by Proposition \ref{prop:g-eq-blow-up}.
    By construction and Proposition
    \ref{porp:properties-sym3-comp} (\ref{porp:properties-sym3-comp:EZ-action}), the $A_3$-actions on $Z$ 
    and $E_Z\cap E_\Delta'$ are trivial. Thus, we
    get $\chi_{S_3}(Z/k)^{S_3} = \chi_{S_3}(Z/k)^{\tau_i} = \langle 2\rangle \cdot \chi_c(\PP^{n-1}/k)\chi_c(X/k)$
    and $\chi_{S_3}(E_Z\cap E_\Delta'/k)^{S_3} = \chi_{S_3}(E_Z\cap E_\Delta'/k)^{\tau_i}
    = \langle 2\rangle \cdot \chi_c(\PP^{n-1}/k)^2\chi_c(X/k)$ using Proposition
    \ref{porp:properties-sym3-comp}.

    Corollary \ref{cor:g-eq-normal-bundle-fixed} yields $\chi_{S_3}(E_\Delta/k)^{S_3} = \chi_c(E_\Delta/k) = n\cdot
    \chi_c(X/k)\cdot H$. Putting things together, we get
    \begin{equation}
      \label{eq:edelta-euler-char:edeltap}
      \begin{split}
        \chi_{S_3}(E_\Delta'/k)^{S_3}
        &= (\chi_{S_3}(E_\Delta/k) - \chi_{S_3}(Z/k) + \chi_{S_3}(E_Z\cap E_\Delta'/k))^{S_3}\\
        &= n\cdot\chi_c(X/k)\cdot H - \langle 2\rangle \cdot (\langle 1\rangle
          -\chi_c(\PP^{n-1}/k))\chi_c(\PP^{n-1}/k)\chi_c(X/k).
      \end{split}
    \end{equation}
    Thus, we get
    \begin{align*}
      \chi_c((E_\Delta/S_3)/k)
      &= \chi_c((E_\Delta''/S_3)/k)
        + 3(\langle 1\rangle - \chi_c(\PP^{n-1}/k))\chi_c(\PP^{n-1}/k)\chi_c(X/k)\\
      &= \langle 6\rangle \chi_{S_3}(E_\Delta'/k)^{S_3} - \langle 2\rangle \cdot (\chi_{\tau_i}(E_i/k) -
        \chi_{\tau_i}(E_{E_i}/k))^{\tau_i}\\
      &\quad + 3(\langle 1\rangle -
        \chi_c(\PP^{n-1}/k))\chi_c(\PP^{n-1}/k)\chi_c(X/k)\\
      &= \langle 6\rangle \cdot (n\cdot\chi_c(X/k)\cdot H - \langle 2\rangle
        \cdot (\langle 1\rangle
        -\chi_c(\PP^{n-1}/k))\chi_c(\PP^{n-1}/k)\chi_c(X/k))\\
      &\quad - 2\cdot \langle 2\rangle
        \cdot (\langle 1\rangle - \chi_c(\PP^{n-1}/k))\chi_c(\PP^{n-1}/k)\chi_c(X/k)\\
      &\quad + 3(\langle 1\rangle -
        \chi_c(\PP^{n-1}/k))\chi_c(\PP^{n-1}/k)\chi_c(X/k)\\
      &= n\cdot \chi_c(X/k)\cdot H\\
        &\quad + (\langle 1\rangle -\langle 3\rangle) \cdot (\langle 1\rangle -
        \chi_c(\PP^{n-1}/k))\chi_c(\PP^{n-1}/k)\chi_c(X/k)\\
        &\quad + 2(\langle 1\rangle -\langle 2\rangle) \cdot (\langle 1\rangle -
        \chi_c(\PP^{n-1}/k))\chi_c(\PP^{n-1}/k)\chi_c(X/k),
    \end{align*}
    where we use \eqref{eq:edelta-euler-char:edelta} for the first equality;
    \eqref{eq:edelta-euler-char:edeltapp} for the second equality; and
    \eqref{eq:edelta-euler-char:ei}, \eqref{eq:edelta-euler-char:eei}, and \eqref{eq:edelta-euler-char:edeltap}
    for the third equality.
    Since $(\langle 1\rangle - \chi_c(\PP^{n-1}/k))\chi_c(\PP^{n-1}/k)$ is
    hyperbolic, the additional summands vanish and we get $\chi_c((E_\Delta/S_3)/k)
    = n\cdot \chi_c(X/k)\cdot H$ as desired.
  \end{proof}
\end{lemma}

\begin{lemma}
  \label{lem:sym3-invariants-blow-ups}
  We have
  \[
    \chi_{S_3}(\Bl_Z\Bl_\Delta X^3/k)^{S_3} = \chi_{S_3}(X^3/k)^{S_3} - \left(\langle 1\rangle
      + \langle 2\rangle \cdot \chi_c(\PP^{n-1}/k) -
      \frac{n(n+3)}{2}\cdot H\right)\cdot \chi_c(X/k)
  \]
  in $\GW(k)$.

  \begin{proof}
    We have using Proposition \ref{prop:g-eq-blow-up}
    \begin{align*}
      \chi_{S_3}(\Bl_Z\Bl_\Delta X^3/k)
      &= \chi_{S_3}(\Bl_\Delta X^3/k) - \chi_{S_3}(Z/k) + \chi_{S_3}(E_Z/k)\\
      &= \chi_{S_3}(X^3/k) - \chi_{S_3}(\Delta/k) + \chi_{S_3}(E_\Delta/k) - \chi_{S_3}(Z/k) + \chi_{S_3}(E_Z/k).
    \end{align*}
    Since taking invariants is additive, we
    can consider each term separately. Since the $S_3$-action on $\Delta$ is
    trivial, we have $\chi_{S_3}(\Delta/k)^{S_3} = \chi_c(\Delta/k) =
    \chi_c(X/k)$. By Corollary \ref{cor:g-eq-normal-bundle-fixed}, we
    get $\chi_{S_3}(E_\Delta/k)^{S_3} = \chi_c(E_\Delta/k) = n\cdot \chi_c(X/k)\cdot H$.

    Proposition \ref{porp:properties-sym3-comp}
    (\ref{porp:properties-sym3-comp:EZ-action}) yields that the $A_3$-action
    on $Z$ and $E_Z$ is trivial. Thus, we only need to
    consider the $\tau_i$-action. For this, we computed in Proposition
    \ref{porp:properties-sym3-comp} that $\chi_{S_3}(Z/k)^{\tau_i} =
    \langle 2\rangle \cdot \chi_c(\PP^{n-1}/k)\chi_c(X/k)$ and $\chi_{S_3}(E_Z/k)^{\tau_i} =
    \frac{n(n+1)}{2}\cdot \chi_c(X/k)\cdot H$.

    Putting this together, we obtain
    \begin{align*}
      \chi_{S_3}(\Bl_Z\Bl_\Delta X^3/k)^{S_3}
      &= \chi_{S_3}(X^3/k)^{S_3} - \chi_{S_3}(\Delta/k)^{S_3} +
        \chi_{S_3}(E_\Delta/k)^{S_3}\\
      &\quad - \chi_{S_3}(Z/k)^{S_3} +
        \chi_{S_3}(E_Z/k)^{S_3}\\
      &= \chi_{S_3}(X^3/k)^{S_3} - \chi_c(X/k) + n\cdot \chi_c(X/k)\cdot H\\
      &\quad -\langle 2\rangle\cdot \chi_c(\PP^{n-1}/k)\chi_c(X/k)
        + \frac{n(n+1)}{2}\cdot \chi_c(X/k)\cdot H\\
      &= \chi_{S_3}(X^3/k)^{S_3} - \left(\langle 1\rangle + \langle 2\rangle \cdot \chi_c(\PP^{n-1}/k) -
        \frac{n(n+3)}{2}\cdot H\right)\cdot \chi_c(X/k).
    \end{align*}
  \end{proof}
\end{lemma}

\begin{lemma}
  \label{lem:sym3-quotient-blow-ups}
  We have  
  \begin{align*}
    \chi_c(\Sym^3X/k)
    &= \chi_c(((\Bl_{\coprod_iF_i}\Bl_Z\Bl_\Delta X^3)/S_3)/k)\\
    &\quad + \left(\langle 1\rangle + (\langle 1\rangle - \langle
      (-1)^n\rangle)\chi_c(\PP^{n-1}/k) - \frac{n(n+3)}{2}\cdot H\right) \cdot
      \chi_c(X/k) \\
    &\quad + (\langle 1\rangle - \chi_c(\PP^{n-1}/k)\chi_c(F_i/k)
  \end{align*}
  in $\GW(k)$.

  \begin{proof}
    We consider the blow-up sequence \eqref{eq:sym3-blow-up-sequence} after
    quotienting out the $S_3$-action. Thus, we get
    \begin{align*}
      \chi_c(\Sym^3X/k)
      &= \chi_c(((\Bl_\Delta X^3)/S_3)/k) +\chi_c((\Delta/S_3)/k) -
        \chi_c((E_\Delta/S_3)/k)\\
      &= \chi_c(((\Bl_Z\Bl_\Delta X^3)/S_3)/k) + \chi_c((Z/S_3)/k) -
        \chi_c((E_Z/S_3)/k)\\
      &\quad+\chi_c((\Delta/S_3)/k) -
        \chi_c((E_\Delta/S_3)/k)\\
      &= \chi_c(((\Bl_{\coprod_iF_i}\Bl_Z\Bl_\Delta X^3)/S_3)/k)
        + \chi_c((\coprod\nolimits_iF_i/S_3)/k) - \chi_c((\coprod\nolimits_iE_{F_i}/S_3)/k)\\
      &\quad + \chi_c((Z/S_3)/k) - \chi_c((E_Z/S_3)/k) +\chi_c((\Delta/S_3)/k) -
        \chi_c((E_\Delta/S_3)/k).
    \end{align*}
    The $A_3$-action only permutes the $F_i$, and the $\tau_i$-action on $F_i$ is
    trivial. Thus, we have $(\coprod_iF_i)/S_3 = F_i/\tau_i = F_i$ and
    similarly $(\coprod_i E_{F_i})/S_3 = E_{F_i}$. Hence, Proposition
    \ref{porp:properties-sym3-comp} yields
    \[
      \chi_c((\coprod\nolimits_iF_i/S_3)/k) - \chi_c((\coprod\nolimits_iE_{F_i}/S_3)/k) =
      (\langle 1\rangle -\chi_c(\PP^{n-1}/k))\chi_c(F_i/k)
      -\langle (-1)^n\rangle \chi_c(\PP^{n-1}/k)\chi_c(X/k).
    \]

    The $A_3$-action on $Z$ and $E_Z$ is trivial and the $\tau_i$-action on
    $Z$ and $E_Z$ has no fixed points. Therefore, the quotients by the $A_3$-
    and the $\tau_i$-action
    are smooth by Theorem \ref{thm:divisor-smooth-quotient}. Thus, we have $\chi_c((Z/S_3)/k) = \langle
    2\rangle \cdot \chi_{\tau_i}(Z/k)^{\tau_i} = \chi_c(\PP^{n-1}/k)\chi_c(X/k)$ and
    $\chi_c((E_Z/S_3)/k) = \frac{n(n+1)}{2}\cdot \chi_c(X/k)\cdot H$ by Theorem
    \ref{thm:smooth-quotient-euler-char} and
    Proposition \ref{porp:properties-sym3-comp}.

    The action of $S_3$ on $\Delta$ is trivial, so $\chi_c((\Delta/S_3)/k) =
    \chi_c(X/k)$. Putting things together, we get using Lemma
    \ref{lem:edelta-euler-char}
    \begin{align*}
      \chi_c(\Sym^3X/k)
      &= \chi_c(((\Bl_{\coprod_iF_i}\Bl_Z\Bl_\Delta X^3)/S_3)/k)\\
      &\quad + \chi_c((\coprod\nolimits_iF_i/S_3)/k) - \chi_c((\coprod\nolimits_iE_{F_i}/S_3)/k)\\
      &\quad + \chi_c((Z/S_3)/k) - \chi_c((E_Z/S_3)/k)\\
      &\quad +\chi_c((\Delta/S_3)/k) - \chi_c((E_\Delta/S_3)/k)\\
      &= \chi_c(((\Bl_{\coprod_iF_i}\Bl_Z\Bl_\Delta X^3)/S_3)/k)\\
      &\quad-\langle (-1)^n\rangle \chi_c(\PP^{n-1}/k)\chi_c(X/k) + (\langle
        1\rangle -\chi_c(\PP^{n-1}/k))\chi_c(F_i/k)\\
      &\quad + \chi_c(\PP^{n-1}/k)\chi_c(X/k) - \frac{n(n+1)}{2}\cdot \chi_c(X/k)\cdot H\\
      &\quad +\chi_c(X/k) - n\cdot \chi_c(X/k)\cdot H\\
      &= \chi_c(((\Bl_{\coprod_iF_i}\Bl_Z\Bl_\Delta X^3)/S_3)/k)\\
      &\quad + \left(\langle 1\rangle + (\langle 1\rangle - \langle
        (-1)^n\rangle)\chi_c(\PP^{n-1}/k) - \frac{n(n+3)}{2}\cdot H\right) \cdot
        \chi_c(X/k)\\
      &\quad + (\langle 1\rangle -\chi_c(\PP^{n-1}/k))\chi_c(F_i/k)
    \end{align*}
    as desired.
  \end{proof}
\end{lemma}

\begin{theorem}
  \label{thm:euler-char-sym3}
  We have
  \[
    \chi_c(\Sym^3X/k) = \langle 6\rangle
    \cdot \chi_{S_3}(X^3/k)^{S_3}
    + (\langle 1\rangle - \langle 6\rangle)\cdot
    \chi_c(X/k) + (\langle 1\rangle -\langle 2\rangle)(\chi_c(X/k) - \langle
    1\rangle)\chi_c(X/k)
  \]
  in $\GW(k)$.

  \begin{proof}
    By Lemma \ref{lem:sym3-quotient-blow-ups}, we know
    \begin{align*}
      \chi_c(\Sym^3X/k)
      &= \chi_c((\Bl_{\coprod_iF_i}\Bl_Z\Bl_\Delta X^3)/S_3)\\
      &\quad + \left(\langle 1\rangle + (\langle 1\rangle - \langle
      (-1)^n\rangle)\chi_c(\PP^{n-1}/k) - \frac{n(n+3)}{2}\cdot H\right) \cdot
        \chi_c(X/k)\\
      &\quad + (\langle 1\rangle - \chi_c(\PP^{n-1}/k))\chi_c(F_i/k).
    \end{align*}
    Note that $(\Bl_{\coprod_iF_i}\Bl_Z\Bl_\Delta X^3)/S_3 =
    (\Bl_{F_i}((\Bl_Z\Bl_\Delta X^3)/A_3))/\tau_i$ is smooth by
    Theorem~\ref{thm:divisor-smooth-quotient} and the actions of $A_3$ on
    $(\Bl_Z\Bl_\Delta X^3)$ and of $\tau_i$ on $\Bl_{F_i}((\Bl_Z\Bl_\Delta
    X^3)/A_3)$ are free on the complement of the fixed-point locus. Thus, we
    have by Theorem~\ref{thm:smooth-quotient-euler-char} and Proposition~\ref{prop:g-eq-blow-up}
    \begin{align*}
      \chi_c(((\Bl_{\coprod_iF_i}\Bl_Z\Bl_\Delta X^3)/S_3)/k)
      &= \langle 2\rangle
        \cdot \chi_{\tau_i}((\Bl_{F_i}((\Bl_Z\Bl_\Delta X^3)/A_3))/k)^{\tau_i}\\
      &= \langle 2\rangle \cdot (\chi_{\tau_i}(((\Bl_Z\Bl_\Delta X^3)/A_3)/k) -
        \chi_{\tau_i}(F_i/k) + \chi_{\tau_i}(E_{F_i}/k))^{\tau_i}.
    \end{align*}
    The $\tau_i$-action on $F_i$ and $E_{F_i}$ is trivial and thus, we have
    $\chi_{\tau_i}(F_i/k)^{\tau_i} = \chi_c(F_i/k)$ and $\chi_{\tau_i}(E_{F_i}/k)^{\tau_i} =
    \chi_c(E_{F_i}/k) = \langle(-1)^n\rangle \cdot \chi_c(\PP^{n-1}/k)\chi_c(X/k) +
    \chi_c(\PP^{n-1}/k)\chi_c(F_i/k)$ by Proposition~\ref{porp:properties-sym3-comp}. Furthermore, the quotient
    $(\Bl_Z\Bl_\Delta
    X^3)/A_3$ is smooth and thus, we again get using Theorem~\ref{thm:smooth-quotient-euler-char} that
    \[
      \chi_{\tau_i}(((\Bl_Z\Bl_\Delta X^3)/A_3)/k) = \langle 3\rangle \cdot
      \chi_{S_3}((\Bl_Z\Bl_\Delta X^3)/k)^{A_3}.
    \]
    Putting these statements together, we obtain
    \begin{align*}
      \chi_c(((\Bl_{\coprod_iF_i}\Bl_Z\Bl_\Delta X^3)/S_3)/k)
      &= \langle 6\rangle
        \cdot \chi_{S_3}((\Bl_Z\Bl_\Delta X^3)/k)^{S_3}\\
      &\quad + \langle 2\rangle \cdot \langle(-1)^n\rangle \cdot
        \chi_c(\PP^{n-1}/k)\chi_c(X/k)\\
      &\quad - \langle 2\rangle (\langle 1\rangle -\chi_c(\PP^{n-1}/k))\chi_c(F_i/k).
    \end{align*}
    Thus, we get using Lemma \ref{lem:sym3-invariants-blow-ups} and the observation
    that multiplying a rank zero form by a hyperbolic form yields zero that
    \begin{align*}
      \chi_c(\Sym^3X/k)
      &= \langle 6\rangle
      \cdot \chi_{S_3}((\Bl_Z\Bl_\Delta X^3)/k)^{S_3} + \langle 2\rangle \cdot \langle(-1)^n\rangle \cdot
      \chi_c(\PP^{n-1}/k)\chi_c(X/k)\\
      &\quad + \left(\langle 1\rangle + (\langle 1\rangle - \langle
      (-1)^n\rangle)\chi_c(\PP^{n-1}/k) - \frac{n(n+3)}{2}\cdot H\right) \cdot
        \chi_c(X/k)\\
      &\quad + (\langle 1\rangle -\langle 2\rangle)(\langle 1\rangle -\chi_c(\PP^{n-1}/k))\chi_c(F_i/k)\\
      &= \langle 6\rangle
      \cdot \chi_{S_3}((\Bl_Z\Bl_\Delta X^3)/k)^{S_3}\\
      &\quad + \left(\langle 1\rangle + (\langle 1\rangle + (\langle 2\rangle -\langle
        1\rangle)\cdot  \langle
      (-1)^n\rangle)\chi_c(\PP^{n-1}/k) - \frac{n(n+3)}{2}\cdot H\right) \cdot
        \chi_c(X/k)\\
      &\quad + (\langle 1\rangle -\langle 2\rangle)(\langle 1\rangle -\chi_c(\PP^{n-1}/k))\chi_c(F_i/k)\\
      &= \langle 6\rangle
      \cdot \left(\chi_{S_3}(X^3/k)^{S_3} - \left(\langle 1\rangle + \langle 2\rangle \cdot \chi_c(\PP^{n-1}/k) -
        \frac{n(n+3)}{2}\cdot H\right)\cdot \chi_c(X/k)\right)\\
      &\quad + \left(\langle 1\rangle + (\langle 1\rangle + (\langle 2\rangle -\langle
        1\rangle)\cdot  \langle
      (-1)^n\rangle)\chi_c(\PP^{n-1}/k) - \frac{n(n+3)}{2}\cdot H\right) \cdot
        \chi_c(X/k)\\
      &\quad + (\langle 1\rangle -\langle 2\rangle)(\langle 1\rangle -\chi_c(\PP^{n-1}/k))\chi_c(F_i/k)\\
      &= \langle 6\rangle
      \cdot \chi_{S_3}(X^3/k)^{S_3}\\
      &\quad + (\langle 1\rangle - \langle 6\rangle)(\langle 1\rangle - \frac{n(n+3)}{2}\cdot H) \cdot
        \chi_c(X/k)\\
      &\quad + (\langle 1\rangle - \langle 3 \rangle)\cdot \underset{\text{hyperbolic}}{\underbrace{\chi_c(\PP^{n-1}/k)\chi_c(X/k)}}\\
      &\quad +(\langle 2\rangle -\langle
        1\rangle)\cdot  \langle
        (-1)^n\rangle\underset{\text{hyperbolic}}{\underbrace{\chi_c(\PP^{n-1}/k)\chi_c(X/k)}}\\
      &\quad + (\langle 1\rangle -\langle 2\rangle)(\langle 1\rangle -\chi_c(\PP^{n-1}/k))\chi_c(F_i/k)\\
      &= \langle 6\rangle
      \cdot \chi_{S_3}(X^3/k)^{S_3}
        + (\langle 1\rangle - \langle 6\rangle)\cdot
        \chi_c(X/k)\\
      &\quad + (\langle 1\rangle -\langle 2\rangle)(\langle 1\rangle -\chi_c(\PP^{n-1}/k))\chi_c(F_i/k).\\
    \end{align*}
    By Proposition \ref{porp:properties-sym3-comp}, we have
    \[
      (\langle 1\rangle -\chi_c(\PP^{n-1}/k))\chi_c(F_i/k) = (\langle 1\rangle
      -\chi_c(\PP^{n-1}/k))(\chi_c(X/k) - \langle 1\rangle + 2\chi_c(\PP^{n-1}/k))\chi_c(X/k).
    \]
    Since at least one of $\chi_c(\PP^{n-1}/k)$ and $\chi_c(X/k)$ is hyperbolic
    for dimension reasons and multiplication with $\langle 1\rangle - \langle
    2\rangle$ kills hyperbolic forms, we get
    \[
      (\langle 1\rangle -\langle 2\rangle)(\langle 1\rangle
      -\chi_c(\PP^{n-1}/k))\chi_c(F_i/k) = (\langle 1\rangle - \langle 2
      \rangle)(\chi_c(X/k) - \langle 1\rangle) \chi_c(X/k).
    \]
    If we combine this with our already done computation of $\chi_c(\Sym^3X/k)$,
    we get the formula in the statement.
  \end{proof}
\end{theorem}

\begin{lemma}
  \label{lem:quadratic-form-s3-fixed}
  In the notation and under the assumptions of Proposition \ref{prop:fixed-points-quadratic-action} (\ref{prop:fixed-points-quadratic-action:no-sign}), we
  have for $G = S_3$ that
  \[
    (\beta^3)^{S_3} = (\langle 1\rangle + (n-1)\langle 3\rangle)\beta +
    \sum_{1\le i<j<k\le n}\langle 6a_ia_ja_k\rangle.
  \]

  \begin{proof}
    We expand the formula from Proposition
    \ref{prop:fixed-points-quadratic-action} (\ref{prop:fixed-points-quadratic-action:no-sign}). Thus, we need to classify the
    orbits of elements of $\{1, \dots, n\}^3$ under the $S_3$-action. The
    elements $(1, 1, 1), \dots, (n,n,n)$ are fixed under the action.
    
    Elements representing the orbits  with an orbit size of three are $(i,j,j)$ for $i\ne j \in
    \{1, \dots, n\}$. Representatives for the orbits with an orbit size of six are $(i,j,k)$ for
    $1 \le i<j<k\le n$.
    Since this covers all elements of $\{1, \dots, n\}^3$, there are no
    elements with an orbit size of two.
    Thus, we get the desired formula from Proposition
    \ref{prop:fixed-points-quadratic-action} (\ref{prop:fixed-points-quadratic-action:no-sign})
    after noting that $\langle \alpha_i\alpha_j^2\rangle = \langle \alpha_i\rangle$ and
    $\langle \alpha_i^3\rangle = \langle \alpha_i\rangle$.
  \end{proof}
\end{lemma}

\begin{lemma}
  \label{lemma:virtual-quadratic-form-s3-fixed}
  Let $V^\bullet$ be a $\Z$-graded finite-dimensional vector space
  over $k$. Let $\beta^\bullet$ be a $\Z$-graded
  non-degenerate symmetric bilinear form on $V^\bullet$. Let $S_3$ act on
  $\beta^3$ by permuting the factors.
  Then we have
  \[
    (\can_\Z(\beta^3))^{S_3} = ((\beta^\even)^3)^{S_3} + l\langle -6\rangle \beta^\even +
    ml(l-1)H - \left(m+\binom m2\right)lH - \frac 12 \binom{2l}{3}H
  \]
  in $\GW(k)$,
  where $l = \frac 12 \dim_kV^\odd = \frac 12 \rank\beta^\odd$ and $m =\dim_kV^\even = \rank
  \beta^\even$.

  \begin{proof}
    Note that since $\beta^\odd$ is alternating, it corresponds with $lH
    \in GW(k)$. We have
    \begin{align*}
      V^{\otimes 3}
      &= (V^\even \oplus V^\odd)^{\otimes 3}\\
      &= \underset{= (V^{\otimes 3})^\even}{\underbrace{(V^\even)^{\otimes 3} \oplus
        V^\even\otimes V^\odd\otimes V^\odd \oplus V^\odd\otimes V^\even\otimes V^\odd \oplus
        V^\odd\otimes V^\odd\otimes V^\even}}\\
      &\quad \oplus
        \underset{=(V^{\otimes 3})^\odd}{\underbrace{V^\even\otimes V^\even\otimes V^\odd \oplus V^\even\otimes V^\odd\otimes V^\even \oplus V^\odd\otimes V^\even\otimes V^\even
        \oplus (V^\odd)^{\otimes 3}}}.
    \end{align*}
    The $S_3$-action factors into actions on $(V^\even)^{\otimes 3}$,
    \begin{gather*}
      V_1 =
      V^\even\otimes V^\odd\otimes V^\odd \oplus V^\odd\otimes V^\even\otimes V^\odd \oplus 
      V^\odd\otimes V^\odd\otimes V^\even,\\
      V_2 = V^\even\otimes V^\even\otimes V^\odd \oplus
      V^\even\otimes V^\odd\otimes V^\even \oplus V^\odd\otimes V^\even\otimes V^\even,
    \end{gather*}
    and
    $(V^\odd)^{\otimes 3}$. Thus,
    \[
      (V^{\otimes 3})^{S_3} = \underset{= ((V^{\otimes 3})^\even)^{S_3}}{\underbrace{((V^\even)^{\otimes 3})^{S_3} \oplus
          V_1^{S_3}}} \oplus \underset{((V^{\otimes
          3})^\odd)^{S_3}}{\underbrace{V_2^{S_3} \oplus ((V^\odd)^{\otimes 3})^{S_3}}}.
    \]
    The summand $((V^\even)^{\otimes 3})^{S_3}$ contributes
    $+((\beta^\even)^3)^{S_3}$ to $(\can_{\Z}(\beta^3))^{S_3}$.

    Consider the cyclic permutation $\sigma = (1\; 2\; 3)\in S_3$.
    Let $w_1^\odd, \dots, w_N^\odd$ be a
    basis for $((V^\odd)^{\otimes 2})^{S_2}$, where $S_2$ acts on $(V^\odd)^{\otimes 2}$
    by swapping the factors with a sign, and let $v_1^\even, \dots,
    v_m^\even$ be a basis of $V^\even$. Therefore,
    \[
      \left\{\sum_{\nu=0}^2 \sigma^\nu \cdot v_i^\even\otimes w_j^\odd\mid i=1,
      \dots, m; j = 1, \dots, N\right\}
    \]
    is a basis of
    $V_1^{S_3}$. Similarly, if $w_1^\even, \dots, w_{N'}^\even$ is a
    basis of
    $((V^\even)^{\otimes 2})^{S_2}$, where $S_2$ acts on $(V^\even)^{\otimes 2}$ by swapping the factors
    without introducing a sign, and $v^\odd_1,\dots, v^{\odd}_{2l}$ is a
    basis of $V^\odd$, then
    \[
      \left\{\sum_{\nu=0}^2 \sigma^\nu \cdot v_i^\odd\otimes w_j^\even\mid i=1,
      \dots, 2l; j = 1, \dots, N'\right\}
    \]
    is a basis of $V_2^{S_3}$. In particular by Proposition
    \ref{prop:fixed-points-quadratic-action}
    (\ref{prop:fixed-points-quadratic-action:no-sign}), we have
    \[\dim V_2^{S_3} =
      2l\cdot \dim ((V^\even)^{\otimes
        2})^{S_2} = 2l(m + \binom m2).
    \]
    Thus, the restriction of
    $\can_{\Z}(\beta^3)$ to $V_2^{S_3}$ contributes
    $-l(m+\binom m2)H$ to $(\can_\Z(\beta^3))^{S_3}$.
    Furthermore, we have
    \[
      \beta^3\left(\sum_{\nu=0}^2 \sigma^\nu v_i^\even\otimes w_j^\odd, \sum_{\nu=0}^2 \sigma
      v_i^\even\otimes w_j^\odd\right)= 3\beta^\even(v_i^\even, v_i^\even) \cdot (\beta^\odd)^2(w_j^\odd,w_j^\odd).
    \]
    Thus, using Lemma \ref{lem:sym2-vector-space-comp} for $\beta^\odd$ as a
    quadratic form on $0 \oplus V^\odd$, we get that $V_1^{S_3}$ contributes
    \[
      \langle 3\rangle \beta^\even(l\langle -2\rangle + l(l-1)H) = l\langle
      -6\rangle \beta^\even + ml(l-1)H
    \]
    to $(\can_\Z(\beta^3))^{S_3}$.

    Thus, we are left with computing the dimension of $((V^-)^{\otimes
      3})^{S_3}$. Proposition~\ref{prop:fixed-points-quadratic-action}~(\ref{prop:fixed-points-quadratic-action:sign})
    and the proof of Lemma \ref{lem:quadratic-form-s3-fixed}
    yield $\dim ((V^\odd)^{\otimes 3})^{S_3} = \binom{2l}{3}$, and thus, this
    summand contributes $-\frac 12\binom{2l}{3}$ to $(\beta^3)^{S_3}$. If we
    combine these contributions, we get the desired formula.
  \end{proof}
\end{lemma}

\begin{theorem}\label{thm:sym3}
  Let $X$ be a connected, smooth, projective scheme over $k$. Express the $\A^1$-Euler characteristic of $X$ as
  $\chi_c(X/k) = \sum_{i=1}^m\langle \alpha_i\rangle -lH \in \GW(k)$ such that
  $l = \frac{1}{2}\cdot \sum_{p-q \text{ odd}}\dim_k H^q(X,\Omega^p_{X/k})$ and write $\beta = \sum_{i=1}^m \langle
  \alpha_i\rangle$. Then
  \begin{align*}
    \chi_c(\Sym^3 X/k)
    &= \langle 6\rangle
      \cdot \Bigg((\langle 1\rangle + (m-1)\langle 3\rangle)\beta +
      \sum_{1\le i<j<k\le m}\langle 6\alpha_i\alpha_j\alpha_k\rangle\\
    &\qquad  + l\langle -6\rangle \beta +
      ml(l-1)H - \left(m+\binom m2\right)lH - \frac 12 \binom{2l}{3}H\Bigg)\\
    &\quad + (\langle 1\rangle - \langle 6\rangle)\cdot
      \chi_c(X/k) + (\langle 1\rangle -\langle 2\rangle)(\chi_c(X/k) - \langle
      1\rangle)\chi_c(X/k)\\
  \end{align*}
  in $\GW(k)$.

  \begin{proof}
    Example
    \ref{ex:graded-euler-char} yields that $\beta$ is the even part of the trace form on Hodge
    cohomology
    \begin{equation}
      \label{eq:sym3:trace-form}
      \left(\bigoplus_{i,j=0}^{\dim_kX}
        H^i(X,\Omega^j_{X/k})[j-i],\Tr\right)
    \end{equation}
    from Example \ref{ex:graded-euler-char} and the odd part of
    \eqref{eq:sym3:trace-form} is isometric to $lH$.
    The theorem now follows by combining Theorem \ref{thm:euler-char-sym3}, Lemma
    \ref{lem:quadratic-form-s3-fixed}, and Lemma
    \ref{lemma:virtual-quadratic-form-s3-fixed}.
  \end{proof}
\end{theorem}

\begin{rem}
  \label{rem:sym3-formula-indep}
  We shall see in Theorem \ref{thm:sym3-comparison} that $\chi_c(\Sym^3X/k) =
  a_3^\PPal(\chi_c(X/k))$ by explicitly calculating $a_3^\PPal(\chi_c(X/k))$ and
  comparing the result with the formula in Theorem \ref{thm:sym3}. Since the
  $a_3^\PPal(\chi_c(X/k))$ is independent of the choice of $l\ge 0$ in Theorem~\ref{thm:sym3}, we get that the assumption that $l = \frac{1}{2}\cdot
  \sum_{p-q \text{ odd}}\dim_k H^q(X,\Omega^p_{X/k})$ can be replaced with
  the assumption that $l \ge 0$ is an integer.
\end{rem}

\section{Comparison with the Power Structure Prediction}
\label{sec:comparison-sym23}
We shall now compare the computation of $\chi_c(\Sym^2X/k)$ and
$\chi_c(\Sym^3X/k)$ from Theorem~\ref{thm:qec-sym2} and Theorem~\ref{thm:sym3} to the prediction from the power structure $a_\ast^\PPal$ of
Pajwani-Pál.

\begin{proposition}
  \label{prop:q-form-power-structure-prediction}
  Let $\beta = \langle \alpha_1\rangle + \dots + \langle \alpha_m\rangle \in \GW(k)$ be a
  quadratic form. Then
  \[
    a^\PPal_2(\beta) = m\cdot \langle 2\rangle + \sum_{0 \le i < j \le m}\langle \alpha_i\alpha_j\rangle + (\langle
    -2\rangle - \langle -1\rangle)\beta\in \GW(k).
  \]

  \begin{proof}
    We have $a^\PPal_0(\langle \alpha_i\rangle) = \langle 1\rangle$, $a^\PPal_1(\langle
    \alpha_i\rangle) = \langle \alpha_i\rangle$, and $a_2^\PPal(\langle \alpha_i\rangle) = \langle
    \alpha_i^2\rangle + \langle 2\rangle + \langle \alpha_i\rangle - \langle 1\rangle -
    \langle 2\alpha_i\rangle$.
    Therefore, we get by Lemma \ref{lem:power-structre-large-sum}
    \begin{align*}
      a_2^\PPal(\beta)
      &= \sum_{i_1+ \dots + i_m =2}a^\PPal_{i_1}(\langle \alpha_1\rangle) \cdots
        a^\PPal_{i_m}(\langle \alpha_m\rangle)\\
      &= \sum_{i=1}^m (\langle \alpha_i^2\rangle + \langle 2\rangle +\langle \alpha_i\rangle -
        \langle 1\rangle - \langle 2\alpha_i\rangle) + \sum_{1 \le i < j \le m}
        \langle \alpha_i\alpha_j\rangle\\
      &= \sum_{1\le i< j\le m}\langle \alpha_i\alpha_j\rangle +
        \underset{=(\langle -2\rangle - \langle -1\rangle)\beta}{\underbrace{\beta -\langle 2\rangle
        \cdot \beta}} + m\cdot \langle 2\rangle\\
      &= m\cdot \langle 2\rangle +  \sum_{0 \le i < j \le m}\langle \alpha_i\alpha_j\rangle + (\langle
        -2\rangle - \langle -1\rangle)\beta.\qedhere
    \end{align*}
  \end{proof}
\end{proposition}

\begin{proposition}
  \label{prop:compuation-a2-general}
  Let $\beta = \langle \alpha_1\rangle + \dots + \langle \alpha_n\rangle - mH \in \GW(k)$
  with $m \ge 0$ be a virtual
  quadratic form. Then
  \[
    a_2^\PPal(\beta) = n\langle 2\rangle \sum_{0 \le i < j \le n}\langle \alpha_i\alpha_j\rangle +
    m\langle -1\rangle + (m^2-(n+1)m)H + (\langle
    -2\rangle - \langle -1\rangle)\beta\in\GW(k).
  \]
  \begin{proof}
    Note $(\langle -2\rangle - \langle -1\rangle)\cdot  H = 0$.
    For $\beta' = \langle \alpha_1\rangle + \dots + \langle \alpha_n\rangle$, we have by
    Proposition \ref{prop:q-form-power-structure-prediction}
    \[
      a_2^\PPal(\beta')
      =n\langle 2\rangle + \sum_{0\le i < j\le n} \langle \alpha_i\alpha_j\rangle
        +(\langle -2\rangle - \langle
        -1\rangle)\cdot\beta
    \]
    Now, we have $a_2^\PPal(-mH) = m\langle -1\rangle + (m^2-m)H$ for $m \in
    \Z_{\ge 0}$ by \cite[Lemma~23]{bv2024quadratic}, and thus,
    \begin{align*}
      a_2^\PPal(\beta) &= a_2^\PPal(\beta') + \underset{=-nmH}{\underbrace{a_1^\PPal(-mH)a_1^\PPal(\beta')}} +
               a_2^\PPal(-mH)\\
      &= n\langle 2\rangle + \sum_{0\le i < j\le n} \langle \alpha_i\alpha_j\rangle
        +(\langle -2\rangle - \langle
        -1\rangle)\cdot
        \beta -nmH + m\langle -1\rangle + (m^2-m)H\\
      &=n\langle 2\rangle+ \sum_{0 \le i < j \le n}\langle \alpha_i\alpha_j\rangle +
    m\langle -1\rangle + (m^2-(n+1)m)H + (\langle
    -2\rangle - \langle -1\rangle)\beta.\qedhere
    \end{align*}
  \end{proof}
\end{proposition}

\begin{theorem}
  \label{thm:compatibility-a2-smooth-projective}
  For a smooth, projective scheme $X$ over $k$ that is not-necessarily connected, we have
  \[
    \chi_c(\Sym^2X/k) =  a_2^\PPal(\chi_c(X/k))
  \]
  in $\GW(k)$.
  \begin{proof}
    Let $X = \coprod_{i=1}^N X_i$ with $X_i$ connected be a decomposition into
    connected components. Since $a_\ast^{\Sym}$ and $a_\ast^\PPal$ are power
    structures, we have
    \[
      \Sym^2X = \sum_{i_1+ \dots+ i_N =2}\Sym^{i_1}X_1\cdots \Sym^{i_N}(X_N)
    \]
    in $K_0^\rs(\Var_k)$ and
    \[
      a_2^\PPal(\chi_c(X/k)) = \sum_{i_1+ \dots+ i_N
        =2}a_{i_1}^\PPal(\chi_c(X_1/k))\cdots a_{i_N}^\PPal(\chi_c(X_N/k))
    \]
    in $\GW(k)$. Because we already know $\chi_c(\Sym^iX/k) =
    a_i^\PPal(\chi_c(X/k))$ for $i = 0,1$, we can therefore assume that $X$ is
    connected.
    Now, compare the formulae in Theorem \ref{thm:qec-sym2} with the formula
    in Proposition \ref{prop:compuation-a2-general} and note that $\chi_c(X/k)$
    is hyperbolic if $X$ is odd-dimensional.
  \end{proof}
\end{theorem}

\begin{lemma}
  \label{lem:computation-a3}
  For a quadratic form $\beta =  \sum_{i=1}^n\langle \alpha_i\rangle \in \GW(k)$,
  we have
  \[
    a_3^\PPal(\beta) = n\cdot \langle 2 \rangle\cdot \beta
        +(\langle 1\rangle - \langle 2\rangle)\cdot \beta^2  
        +\sum_{1\le i <j<k\le n}\langle \alpha_i\alpha_j\alpha_k\rangle\in \GW(k).
  \]

  \begin{proof}
    Lemma \ref{lem:power-structre-large-sum} yields
    \begin{align*}
      a_3^\PPal(\beta)
      &= a_3^\PPal\left(\sum_{i=1}^n\langle \alpha_i\rangle\right)\\
      &= \sum_{i_1+ \dots + i_n = 3}a_{i_1}^\PPal(\langle \alpha_1\rangle ) \cdots
        a^\PPal_{i_n}(\langle \alpha_n\rangle)\\
      &= \sum_{i=1}^na_3^\PPal(\langle \alpha_i\rangle) + \sum_{\substack{0 \le i,j\le
        n\\i \ne j}}\underset{=\langle
      \alpha_i\rangle}{\underbrace{a_1^\PPal(\langle
      \alpha_i\rangle)}}a_2^\PPal(\langle \alpha_j\rangle)\\
      &\quad +\sum_{0\le i <j<k\le n}\underset{=\langle \alpha_i\alpha_j\alpha_k\rangle}{\underbrace{a_1^\PPal(\langle \alpha_i\rangle)a_1^\PPal(\langle
      \alpha_j\rangle)a_1^\PPal(\langle \alpha_k\rangle)}}\\
      &= \sum_{i=1}^n\underset{=\langle \alpha_i^3\rangle +
      \langle 2\rangle +\langle \alpha_i\rangle - \langle 1\rangle -\langle
      2\alpha_i\rangle}{\underbrace{a_3^\PPal(\langle \alpha_i\rangle)}} + \sum_{\substack{0 \le i,j\le
        n\\i \ne j}}\langle \alpha_i\rangle\cdot \underset{=\langle \alpha_j^2\rangle +
      \langle 2\rangle +\langle \alpha_j\rangle - \langle 1\rangle -\langle
      2\alpha_j\rangle}{\underbrace{a_2^\PPal(\langle \alpha_j\rangle)}}\\
      &\quad + \sum_{0\le i <j<k\le n}\langle \alpha_i\alpha_j\alpha_k\rangle\\
      &= \sum_{i=1}^n(\underset{=2\langle 2\rangle}{\underbrace{2\langle 1\rangle}} - \langle 2\rangle)\langle \alpha_i\rangle
        + n\cdot (\langle 2\rangle - \langle 1\rangle)
        + \sum_{\substack{0 \le i,j\le
        n\\i \ne j}}\langle \alpha_i\rangle\cdot (\langle 2\rangle + \langle
      \alpha_j\rangle - \langle 2\alpha_j\rangle)\\
      &\quad +\sum_{0\le i <j<k\le n}\langle \alpha_i\alpha_j\alpha_k\rangle\\
      &= \langle 2 \rangle\cdot \beta
        + n(\langle 2\rangle - \langle 1\rangle)
        + \underset{=(n-1)\cdot \langle 2\rangle\cdot\beta}{\underbrace{\sum_{\substack{0 \le i,j\le
        n\\i \ne j}}\langle 2\alpha_i\rangle}}
      + \underset{=(\langle 2\rangle -\langle 1\rangle)(n\langle 1\rangle - \beta^2)}{\underbrace{\sum_{\substack{0 \le i,j\le
        n\\i \ne j}} (\langle
      1\rangle - \langle 2\rangle)\langle \alpha_i\alpha_j\rangle}}\\
      &\quad +\sum_{0\le i <j<k\le n}\langle \alpha_i\alpha_j\alpha_k\rangle\\
      &= n\cdot \langle 2 \rangle\cdot \beta
        + n\cdot \underset{=0}{\underbrace{2(\langle 2\rangle - \langle 1\rangle)}} +(\langle 1\rangle - \langle 2\rangle)\beta^2  
        +\sum_{0\le i <j<k\le n}\langle \alpha_i\alpha_j\alpha_k\rangle\\
      &= n\cdot \langle 2 \rangle\cdot \beta
        +(\langle 1\rangle - \langle 2\rangle)\cdot \beta^2  
        +\sum_{0\le i <j<k\le n}\langle \alpha_i\alpha_j\alpha_k\rangle.\qedhere
    \end{align*}
  \end{proof}
\end{lemma}

\begin{lemma}
  \label{lem:computation-a3-general}
  For a quadratic form $\beta =  \sum_{i=1}^n\langle \alpha_i\rangle \in \GW(k)$
  and $l \in \Z_{\ge 0}$,
  we have
  \begin{align*}
    a_3^\PPal(\beta - lH)
    &= (n\cdot \langle 2 \rangle -l\cdot \langle 1\rangle)\cdot \beta
        +(\langle 1\rangle - \langle 2\rangle)\cdot \beta^2  
        +\sum_{1\le i <j<k\le n}\langle \alpha_i\alpha_j\alpha_k\rangle\\
    &\quad-\left(nl(1-l) + l\binom n2+\frac 12\cdot \binom{2l}{3}\right)H.
  \end{align*}

  \begin{proof}
    We have
    \[
      a_3^\PPal(\beta - lH) = a_3^\PPal(\beta) + a_2^\PPal(\beta)\cdot a_1^\PPal(-lH) +
      a_1^\PPal(\beta)\cdot a_2^\PPal(-lH) + a_3^\PPal(-lH).
    \]
    By Lemma \ref{lem:computation-a3}, we have
    \[
      a_3^\PPal(\beta) = n\cdot \langle 2 \rangle\cdot \beta
      +(\langle 1\rangle - \langle 2\rangle)\cdot \beta^2
      +\sum_{1\le i <j<k\le n}\langle \alpha_i\alpha_j\alpha_k\rangle.
    \]
    Furthermore, we have $a_1^\PPal(-lH) = -lH$ and $a_1^\PPal(\beta) = \beta$. By \cite[Lemma~23]{bv2024quadratic}, we have
    \[
      a_2^\PPal(-lH) = l(l-1)\langle 1\rangle
      + l^2\langle -1\rangle = l^2H - l\langle 1\rangle
    \]
    and
    \begin{align*}
      a_3^\PPal(-lH)
      &= (-1)^3\cdot \sum_{i=0}^3\binom{l}{i}\binom{l}{3-i}\langle
        (-1)^i\rangle\\
      &= -\binom l0\binom l3 \langle 1\rangle -\binom l1\binom l2
        \langle -1\rangle - \binom l2\binom l1 \langle 1\rangle - \binom
        l3\binom l0 \langle -1\rangle\\
      &= -\left(\binom l3+l\binom l2\right)H \\
      &= -\frac 12\cdot \binom{2l}{3} \cdot H.
    \end{align*}
    By Proposition \ref{prop:q-form-power-structure-prediction}, we have
    \[
      a_2^\PPal(\beta) = n\langle 2\rangle + \sum_{1\le i < j\le n} \langle \alpha_i\alpha_j\rangle
        +(\langle -2\rangle - \langle
        -1\rangle)\cdot\beta.
    \]
    Thus, we get
    \begin{align*}
      a_3^\PPal(\beta - lH)
      &= a_3^\PPal(\beta) + a_2^\PPal(\beta)\cdot a_1^\PPal(-lH) +
        a_1^\PPal(\beta)\cdot a_2^\PPal(-lH) + a_3^\PPal(-lH)\\
      &= n\cdot \langle 2 \rangle\cdot \beta
        +(\langle 1\rangle - \langle 2\rangle)\cdot \beta^2  
        +\sum_{1\le i <j<k\le n}\langle \alpha_i\alpha_j\alpha_k\rangle\\
      &\quad-  (n\langle 2\rangle + \sum_{1\le i < j\le n} \langle \alpha_i\alpha_j\rangle
        +(\langle -2\rangle - \langle
        -1\rangle)\cdot\beta)\cdot lH\\
      &\quad+ \beta \cdot (l^2H - l\langle 1\rangle) - \frac
        12\cdot \binom{2l}{3} \cdot H\\
      &=n\cdot \langle 2 \rangle\cdot \beta
        +(\langle 1\rangle - \langle 2\rangle)\cdot \beta^2  
        +\sum_{1\le i <j<k\le n}\langle \alpha_i\alpha_j\alpha_k\rangle\\
      &\quad- \left(n+\binom n2\right)\cdot lH + nl^2H-l\beta- \frac
        12\cdot \binom{2l}{3} \cdot H\\
      &=(n\cdot \langle 2 \rangle -l\cdot \langle 1\rangle)\cdot \beta
        +(\langle 1\rangle - \langle 2\rangle)\cdot \beta^2  
        +\sum_{1\le i <j<k\le n}\langle \alpha_i\alpha_j\alpha_k\rangle\\
      &\quad -\left(nl(1-l) + l\binom n2+\frac 12\cdot \binom{2l}{3}\right)H.\qedhere
    \end{align*}
  \end{proof}
\end{lemma}

\begin{theorem}
  \label{thm:sym3-comparison}
  Let $X$ be a smooth, projective scheme over $k$ that is not-necessarily
  connected. Assume that $\characteristic k \ne 2,3$. Then we have
  \[
    \chi_c(\Sym^3X/k) = a_3^\PPal(\chi_c(X/k))
  \]
  in $\GW(k)$.

  \begin{proof}
    Let $X = \coprod_{i=1}^N X_i$ with $X_i$ connected be a decomposition into
    connected components. Since $a_\ast^{\Sym}$ and $a_\ast^\PPal$ are power
    structures, we have
    \[
      \Sym^3X = \sum_{i_1+ \dots+ i_N =3}\Sym^{i_1}X_1\cdots \Sym^{i_N}(X_N)
    \]
    in $K_0^\rs(\Var_k)$ and
    \[
      a_3^\PPal(\chi_c(X/k)) = \sum_{i_1+ \dots+ i_N
        =3}a_{i_1}^\PPal(\chi_c(X_1/k))\cdots a_{i_N}^\PPal(\chi_c(X_N/k))
    \]
    in $\GW(k)$. Because we already know $\chi_c(\Sym^iX/k) =
    a_i^\PPal(\chi_c(X/k))$ for $i = 0,1,2$ by Theorem \ref{thm:compatibility-a2-smooth-projective}, we can therefore assume that $X$ is
    connected.

    Write $\chi_c(X/k) = \sum_{i=1}^m\langle
    \alpha_i\rangle - lH \in \GW(k)$ with $l$ as in Theorem \ref{thm:sym3} and $\beta = \sum_{i=1}^m\langle
    \alpha_i\rangle$.
    By Theorem \ref{thm:sym3}, we have
    \begin{align*}
      \chi_c(\Sym^3 X/k)
      &= \langle 6\rangle
        \cdot \Bigg((\langle 1\rangle + (m-1)\langle 3\rangle)\beta +
        \sum_{1\le i<j<k\le m}\langle 6\alpha_i\alpha_j\alpha_k\rangle\\
      &\qquad  + l\langle -6\rangle \beta +
        ml(l-1)H - \left(m+\binom m2\right)lH - \frac 12 \binom{2l}{3}H\Bigg)\\
      &\quad + (\langle 1\rangle - \langle 6\rangle)\cdot
        \chi_c(X/k) + (\langle 1\rangle -\langle 2\rangle)(\chi_c(X/k) - \langle
        1\rangle)\chi_c(X/k)\\
      &= (m\langle 2\rangle + l\langle -1\rangle)\beta +
        \sum_{1\le i<j<k\le m}\langle \alpha_i\alpha_j\alpha_k\rangle + (\langle 1\rangle -
        \langle 2\rangle)\beta^2\\
      &\quad + \left(ml(l-2) -l\binom m2 -\frac 12\binom{2l}{3}\right)H.
    \end{align*}
    If we compare this with Lemma \ref{lem:computation-a3-general}, we get
    \begin{align*}
      \chi_c(\Sym^3 X/k) - a_3^\PPal(\chi_c(X/k))
      &= l\langle -1\rangle \beta - mlH +l\beta\\
      &= l(\langle 1\rangle + \langle -1\rangle)\beta - mlH = 0.
    \end{align*}
    Thus, $\chi_c(\Sym^3 X/k) = a_3^\PPal(\chi_c(x/k))$.
  \end{proof}
\end{theorem}

The following proposition is a direct consequence of the proof of \cite[Lemma
2.9]{Pajwani-PalPS}. We include a proof here for the convenience
of the reader.

\begin{proposition}
  \label{prop:negative-additive-ps-compatibility}
  Let $a_\ast$ be a power structure on a ring $R$ and let $a'_\ast$ be a power structure on
  a ring $R'$. Let $\varphi\colon R \to R'$ be a ring
  homomorphism and $r,s \in R$ and $n \ge 0$. If we have $\varphi(a_i(r)) =
  a'_i(\varphi(r))$ and $\varphi(a_i(s)) =
  a'_i(\varphi(s))$ for all $0 \le i \le n$, then we have $\varphi(a_i(-r)) =
  a'_i(\varphi(-r))$ and $\varphi(a_i(r+s)) = a'_i(\varphi(r+s))$ for all $0 \le i \le n$.

  \begin{proof}
    The statement for $r+s$ follows directly from the sum decomposition
    $a_i(r+s) = \sum_{\nu=0}^ia_\nu(r)a_{i-\nu}(s)$. So it only remains to
    prove the statement for $-r$.
    
    We prove this by induction on $n$. For $n = 0$, there is nothing to prove
    since $\varphi(a_0(-r)) = \varphi(1) = 1 = a'_0(\varphi(-r))$.

    Let $n \ge 1$. We get compatibility for $i \le n-1$ by the induction
    hypothesis. So, we only need to prove $\varphi(a_n(-r)) = a'_n(\varphi(-r))$. Now, we have $0 = a_n(r-r) = \sum_{i=0}^na_i(r)a_{n-i}(-r)$
    and thus, $a_n(-r) = -\sum_{i=1}^na_i(r)a_{n-i}(-r)$. This implies
    \[
      \varphi(a_n(-r)) = -\sum_{i=1}^n\varphi(a_i(r))\varphi(a_{n-i}(-r)) =
      -\sum_{i=1}^na'_i(\varphi(r))a'_{n-i}(\varphi(-r)) = a'_n(\varphi(-r)),
    \]
    where we use the induction hypothesis for the second equality.
  \end{proof}
\end{proposition}

\begin{corollary}
  \label{cor:compatibility-ps-generator-check}
  Let $a_\ast$ be a power structure on a ring $R$ and let $a'_\ast$ be a power structure on
  a ring $R'$. Let $\varphi\colon R \to R'$ be a ring
  homomorphism and $n \ge 0$. Let $S \subset R$ be a subset generating $R$ as a group. If we have $\varphi(a_i(r)) =
  a'_i(\varphi(r))$ for all $0 \le i \le n$ and all $r \in S$, then we have $\varphi(a_i(r)) =
  a'_i(\varphi(r))$ for all $0 \le i \le n$ and $r\in R$.

  \begin{proof}
    Let $r \in R$. By Proposition \ref{prop:negative-additive-ps-compatibility}, we can
    assume $S = -S$. Since $S$ generates $R$ as a group, we can find $s_i \in
    S$ such that $r = \sum_{i=1}^ms_{i}$. Again by Proposition
    \ref{prop:negative-additive-ps-compatibility}, we have 
    \[
      a'_n(\varphi(r)) = a'_n\left(\sum_{i=1}^m\varphi(s_{i})\right)
      = \varphi\left(a_n\left(\sum_{i=1}^m s_i\right)\right) = \varphi(a_n(r)).\qedhere
    \]
  \end{proof}
\end{corollary}

\begin{theorem}
  \label{thm:sym23-compatibility-char-zero}
  Let $X$ be a variety over $k$ and suppose $k$ has characteristic zero. Then we
  have
  \[
    \chi_c(\Sym^2X/k) = a_2^\PPal(\chi_c(X/k)) \quad \text{and}\quad \chi_c(\Sym^3X/k) = a_3^\PPal(\chi_c(X/k)).
  \]

  \begin{proof}
    Combine Bittner's Theorem \cite[Theorem~3.1]{bittner_universal_2004} with
    Theorem~\ref{thm:compatibility-a2-smooth-projective},
    Theorem~\ref{thm:sym3-comparison}, and
    Corollary~\ref{cor:compatibility-ps-generator-check}.
  \end{proof}
\end{theorem}

\printbibliography[heading=bibintoc]

~\\
	
\newpage 
\noindent Louisa F. Br\"oring \\
Universit\"at Duisburg-Essen \\
Fakult\"at f\"ur Mathematik, Thea-Leymann-Str. 9, 45127 Essen, Germany \\
E-Mail: \href{mailto:louisa.broering@uni-due.de}{louisa.broering@uni-due.de}\\ \\

\noindent Keywords: motivic homotopy theory, refined enumerative geometry, other fields\\
Mathematics Subject Classification: 14G27, 14N10, 14F42
\end{document}